\documentclass[12pt,leqno]{article}

\usepackage{amsmath,amsfonts,amssymb,epsfig,theorem}
\usepackage{yfonts}
\usepackage[usenames]{color}

\usepackage{indentfirst}
\usepackage{xspace}
\usepackage{multirow}
\usepackage{hyperref}
\usepackage{xcolor}
\hypersetup{colorlinks=true,urlcolor=blue,linkcolor=blue,citecolor=blue}

\usepackage[all]{xy}
\usepackage{psfrag}

\DeclareFontFamily{U}{rsf}{}
\DeclareFontShape{U}{rsf}{m}{n}{
  <5> <6> rsfs5 <7> <8> <9> rsfs7 <10->  rsfs10}{}
\DeclareMathAlphabet{\mathscr}{U}{rsf}{m}{n}
\newcommand{\mycal}[1]{\mathscr{#1}}

\newcommand{\op}[1]{\operatorname{#1}}
\newcommand{\sone}[1]{{}^{\dagger}#1}
\newcommand{\stwo}[1]{{}^{\ddagger}#1}

\newcommand{\bD}{\boldsymbol{D}}
\newcommand{\bO}{\boldsymbol{O}}

\newcommand{\nc}{{\sf\bfseries nc}}

\newcommand{\bw}{{\sf w}}

\newcommand{\sGr}{\text{\sf\bfseries Gr}}

\newcommand{\MF}{\text{\sf\bfseries MF}}
\newcommand{\bDelta}{\boldsymbol{\Delta}}
\newcommand{\vol}{\text{{\sf vol}}}
\newcommand{\Ybar}{\text{\sf Z}}
\newcommand{\wbar}{\text{\sf f}}
\newcommand{\hatwbar}{\widehat{\op{{\sf f}}}}

\newcommand{\bmid}{\boldsymbol{\mid}}
\newcommand{\Fuk}{\text{{\sf\bfseries Fuk}}}
\newcommand{\FS}{\text{{\sf\bfseries FS}}}
\newcommand{\sD}{\text{{\sf\bfseries D}}}
\newcommand{\Perf}{\text{\sf\bfseries Perf}}
\newcommand{\hor}[1]{#1^{\text{\sf h}}}
\newcommand{\ver}[1]{#1^{\text{\sf v}}}
\newcommand{\con}{{\boldsymbol{\iota}}}
\newcommand{\dtg}{{\boldsymbol{\mathfrak{g}}}}
\newcommand{\pvg}{{\boldsymbol{\mathfrak{G}}}}
\newcommand{\sg}{\text{{\sf g}}}
\newcommand{\ddiv}{\op{{\sf div}}}
\newcommand{\bk}{\boldsymbol{\Bbbk}}
\newcommand{\Fr}{\op{{\sf Fr}}}
\newcommand{\car}{\boldsymbol{\gamma}}
\newcommand{\for}{\op{{\sf for}}}
\newcommand{\lFr}{\op{\mathfrak{Fr}}}
\newcommand{\bz}{\boldsymbol{z}}
\newcommand{\beps}{\boldsymbol{\epsilon}}
\newcommand{\bZ}{\boldsymbol{{\mathcal Z}}}
\newcommand{\bGamma}{\boldsymbol{\Gamma}}
\newcommand{\bY}{\boldsymbol{{\mathcal Y}}}
\newcommand{\bE}{\text{\bfseries\sf E}}
\newcommand{\bphi}{\boldsymbol{\phi}}
\newcommand{\Anc}[1]{{}^{\text{{\sf A}}}#1}
\newcommand{\Bnc}[1]{{}^{\text{{\sf B}}}#1}
\newcommand{\Amod}[1]{{}^{\boldsymbol{\mathfrak{a}}}#1}
\newcommand{\Bmod}[1]{{}^{\boldsymbol{\mathfrak{b}}}#1}
\newcommand{\AAmod}[1]{{}^{\boldsymbol{\textgoth{A}}}#1}
\newcommand{\BBmod}[1]{{}^{\boldsymbol{\textgoth{B}}}#1}

\newcommand{\sM}{\op{{\sf M}}}
\newcommand{\bG}{\text{{\sf\bfseries G}}}

\def\punkt{\refstepcounter{subsubsection}
           \noindent{\bf \thesubsubsection.\ }}

\theoremstyle{plain}

\newtheorem{thm}{Theorem}
\newtheorem{theo}{Theorem}[section]
\newtheorem{lemma}[theo]{Lemma}

\newtheorem{cor}[theo]{Corollary}

\newtheorem{conn}[theo]{Conjecture}
\newtheorem{defi}[theo]{Definition}
\newtheorem{prop}[theo]{Proposition}

\newtheorem{claim}[theo]{Claim}
\newtheorem{prop-defi}[theo]{Proposition-Definition}
\newtheorem{lemma-defi}[theo]{Lemma-Definition}
{\theorembodyfont{\rmfamily} \newtheorem{rem}[theo]{Remark}}
{\theorembodyfont{\rmfamily} \newtheorem{ex}[theo]{Example}}
{\theorembodyfont{\rmfamily} }

\numberwithin{equation}{subsection}

\begin{document}

\title{Bogomolov-Tian-Todorov theorems for Landau-Ginzburg models}
\author{L.Katzarkov \and M.Kontsevich \and T.Pantev}
\date{ }
\maketitle

\centerline{\large\it to the memory of Andrey Todorov}

\

\begin{abstract}
  In this paper we prove the smoothness of the moduli space of
  Landau-Ginzburg models. We formulate and prove a
  Bogomolov-Tian-Todorov theorem for the deformations of
  Landau-Ginzburg models, develop the necessary Hodge theory for
  varieties with potentials, and prove a double degeneration statement
  needed for the unobstructedness result. We discuss the various
  definitions of Hodge numbers for non-commutative Hodge structures of
  Landau-Ginzburg type and the role they play in mirror symmetry. We
  also interpret the resulting families of de Rham complexes attacted
  to a potential in terms of mirror symmetry for one parameter
  families of symplectic Fano manifolds and argue that modulo a
  natural triviality property the moduli spaces of Landau-Ginzburg
  models posses canonical special coordinates.
\end{abstract}

\tableofcontents

\section{Introduction} \label{s-introduction}

In this paper we study the local structure of the moduli space of
complex Landau-Ginzburg models. Such a Landau-Ginzburg model is
determined by a pair $(Y,\bw)$, where $Y$ is a complex
quasi-projective variety, and $\bw : Y \to \mathbb{A}^{1}$ is a
holomorphic function on $Y$. Our main objective is to prove the
unobstructedness of the deformations of $(Y,\bw)$ in the case when $Y$
has a trivial canonical class $K_{Y} \cong \mathcal{O}_{Y}$. 
At a
first glance, such a statement is not likely to hold since the
non-compactness of $Y$ will often cause the moduli space of the pair
$(Y,\bw)$ to be infinite dimensional and to have a complicated and
unwieldy local behavior in general.

Before we adress this difficulty it is useful to look at the model
example provided by the classical unobstructedness statement for the
deformations of compact Calabi-Yau varieties. This statement was
proven by different methods by Bogomolov \cite{fedya-ihes}, Tian
\cite{tian}, and Todorov \cite{andrey}.  Recall from
\cite{fedya-ihes}, \cite{tian} and \cite{andrey}, that if $X$ is a
smooth compact Calabi-Yau manifold of dimension $\dim_{\mathbb{C}} X =
d$, then the (formal) versal deformation space $\mycal{M}_{X}$ of $X$
is smooth and of dimension $h^{d-1,1}(X)$. Moreover a choice of a
splitting of the Hodge filtration on $H^{d}_{DR}(X,\mathbb{C})$
defines an analytic affine structure on $\mycal{M}_{X}$. This theorem
has many variants establishing the unobstructedness of deformations of
log Calabi-Yau varieties or Deligne-Mumford stacks, or of weak Fano
varieties or Deligne-Mumford stacks, see e.g.
\cite{ran,kawamata,manetti,kkp,iacono-manetti,iacono,sano}. The log
version of the Bogomolov-Tian-Todorov theorem suggests that if we want
to attain a good control of the deformations of a Landau-Ginzburg
model $(Y,\bw)$, we should look at a nice, e.g. log Calabi-Yau,
compactification $\Ybar$ of $Y$ and consider only deformations that
fix the boundary divisor $D_{\Ybar} = \Ybar - Y$. To streamline this
discussion it will be convenient to distinguish notationally the
varying and the fixed parts in any given deformation problem. Our
convention in that regard will be that when the deformations of some
collection of geometric data are studied, the moving part of the data
will be listed in parentheses, while the part of the data that is kept
fixed will be listed in a subscript. Thus when we say that we are
analyzing the deformations of $(\Ybar,\wbar)_{D_{\Ybar}}$, we mean
that we consider deformations of the pair $(\Ybar,\wbar)$ together
with compatible {\bfseries\it trivial} deformations of the divisor
$D_{\Ybar}$.

In this framework 
prove the following unobstructedness result

\begin{thm} \label{thm:main.TT} Let $\Ybar$ be a smooth projective
  variety, $\wbar : \Ybar \to \mathbb{P}^{1}$ a flat morphism, and
  $D_{\Ybar} \subset \Ybar$ a reduced anti-canonical divisor with
  strict normal crossings.  Assume moreover that $\op{crit}(\wbar)$
  does not intersect the horizontal part of $D_{\Ybar}$, and that the
  vertical part of $D_{\Ybar}$ coincides with the scheme theoretic
  fiber $\wbar^{-1}(\infty)$ of $\wbar$ over $\infty \in
  \mathbb{P}^{1}$.  Then the versal deformation space
  $\mycal{M}_{(\Ybar,\wbar)_{D_{\Ybar}}}$ of
  $(\Ybar,\wbar)_{D_{\Ybar}}$ is smooth.
\end{thm}

\

\noindent
This theorem can be viewed as an unobstructedness result for the
Calabi-Yau Landau-Ginzburg model $(Y,\bw)$ where $Y = \Ybar -
D_{\Ybar}$, $\bw = \wbar_{|Y} : Y \to \mathbb{A}^{1}$. Indeed, the
theorem asserts that if $(Y,\bw)$ admits a compactification
$(\Ybar,\wbar)$ with normal crossing boundary $D_{\Ybar}$, then the
deformations of $(Y,\bw)$ that are ``anchored at infinity'', i.e. the
deformations of the compactification that keep the boundary fixed, are
unobstructed. To prove Theorem~\ref{thm:main.TT} we identify the
$L_{\infty}$-algebra that controls the deformation theory of
$(\Ybar,\wbar)_{D_{\Ybar}}$ and show in Theorem~\ref{theo:ha} that
this $L_{\infty}$-algebra is homotopy abelian. We argue
that, as in the case of compact Calabi-Yau manifolds, the latter
statement can be reduced to a Hodge theoretic property: the double
degeneration property for the Hodge-to-De Rham spectral sequence
associated with the complex of $\wbar$-adapted logarithmic forms (see
Definition~\ref{defi:main}).  This double degeneration is then established in
Theorem~\ref{theo:double.degeneration}.

The setup and conclusion of Theorem~\ref{thm:main.TT} are natural from
the point of view of mirror symmetry. To elaborate on this, note first
that a Landau-Ginzburg pair $(Y,\bw)$ as above will typically arise as
the mirror of a symplective manifold $(X,\omega_{X})$ underlying a
projective Fano variety. Now the homological mirror symmetry
conjecture predicts that the Fukaya category $\Fuk(X,\omega_{X})$ of
$(X,\omega_{X})$ will be equivalent to the category $\MF(Y,\bw)$ of
matrix factorizations of the potential $\bw : Y \to
\mathbb{A}^{1}$. In particular the deformation theories of the Fukaya
category and of the category of matrix factorizations will be
identified. The heuristics motivating Theorem~\ref{thm:main.TT} comes
from the comparison of the corresponding moduli spaces. The versal
deformation space of the Fukaya category is manifestly smooth since it
is an open cone in the space of harmonic 2-forms on $X$. Thus mirror
symmetry predicts that the versal deformation space of the category of
matrix factorizations will also be smooth. Next recall
\cite{orlov,orlov-lagrange,orlov-global} that $\MF(Y,\bw)$ is the
coproduct $\coprod_{\lambda \in \op{crit} \bw}
\sD^{b}_{\op{sing}}\left(Y_{\lambda}\right)$ of the categories of
singularities of the singular fibers of $\bw$. This interpretation
indicates that flat deformations of the geometric data $(Y,\bw)$ will
not necessarily give rise to flat deformations of
$\MF(Y,\bw)$. Indeed, when we deform $(Y,\bw)$ geometrically, the
singularities of fibers of $\bw$ can coalesce and more importantly can
run away to infinity. This will happen for instance if we deform a
compactification of $(Y,\bw)$ so that some interior singular fiber
gets absorbed in the fiber at infinity.  Because of this phenomenon we
will have flat families of Landau-Ginzburg models which give us
families of categories of matrix factorizations whose periodic cyclic
homologies jump. This suggests that we should only consider geometric
deformations of $(Y,\bw)$ that are anchored at infinity. Indeed, if
$((\Ybar,\wbar),D_{\Ybar})$ is a compactification of $(Y,\bw)$, then
the deformations of $(\Ybar,\wbar)$ that fix the boundary divisor
$D_{\Ybar}$ will give deformations of $(Y,\bw)$ without jumps in the
global vanishing cohomology. In this setting the corresponding
categories of matrix factorizations will move in a flat family and we
expect that the deformations of a compactification with a fixed
boundary will provide enough parameters to cover the full versal
deformation space $\mycal{M}_{\MF(Y,\bw)}$ of the category
$\MF(Y,\bw)$. In fact, in the process of proving the double
degeneration property Theorem~\ref{theo:double.degeneration} we will
check that, under the hypothesis of Theorem~\ref{thm:main.TT}, the
natural map of versal deformation spaces
$\mycal{M}_{(\Ybar,\wbar)_{D_{\Ybar}}} \to \mycal{M}_{\MF(Y,\bw)}$ is
\'{e}tale. More precisely, it is not hard to see that the composition
of the isomorphism constructed in
Lemma~\ref{lemma:cohomology.of.the.pair} with Efimov's comparison
isomorphism \cite{efimov-cyclic} can be identified with the
differential of the map \linebreak $\mycal{M}_{(\Ybar,\wbar)_{D_{\Ybar}}} \to
\mycal{M}_{\MF(Y,\bw)}$ at the closed point. In particular this
differential is an isomorphism and so the map is \'{e}tale. Altogether
this heuristic reasoning explains why the unobstructedness of the
deformation theory of $((\Ybar,\wbar),D_{\Ybar})$ is indeed the
expected behavior.

Since the compactified Landau-Ginzburg model
$((\Ybar,\wbar),D_{\Ybar})$ plays a central role in the above
heuristics it is natural to expect that this compactification should
also have a mirror interpretation. A closer look at the associated
Hodge/de Rham data and the double degeneration property of
$((\Ybar,\wbar),D_{\Ybar})$ suggests that the mirror of
$((\Ybar,\wbar),D_{\Ybar})$ is an anticanonical pencil on the symplectic
manifold $(X,\omega_{X})$. In Section~\ref{ssec:mirrorsZf} we
discuss this mirror picture in detail and compare the Hodge theoretic
data appearing on the two sides of this extended mirror
correspondence. We use this analysis to explain how the commutative
pure Hodge structure of the compactified Landau-Ginzburg model arising
from the double degeneration property can be reconstructed from the
non-commutative Hodge structure of the original Landau-Ginzburg model
$(Y,\bw)$. Through the extended mirror symmetry picture we rewrite
this reconstruction process for the Fano mirror and use the resulting
structure to propose a definition of Hodge numbers for $X$ which is
formulated entirely in symplectic terms.

\

\medskip

\noindent{\bfseries Acknowledgments:} We are very grateful to H\'{e}l\`{e}ne
Esnault, Claude Sabbah, Jeng-Daw Yu, and Morihiko Saito for sharing
with us preliminary versions of \cite{esy,saito.appendixE,sabbah.yu}
and for many insigthful letters and conversations on the various
aspects of the irregular Hodge filtration and its asymptotic
properties. Special thanks go to Denis Auroux  and Paul Seidel for
explaining to us various delicate technical points of the variational
theory of Fukaya and wrapped Fukaya categories, and for several illuminating
discussions on Landau-Ginzburg mirror symmetry. 

During the preparation of this work Ludmil Katzarkov was partially
supported by the FRG grant DMS-0854977, and the research grants DMS-
0854977, DMS-0901330, DMS-1265230, DMS-1201475, and by an OISE-1242272
PASI grant from the National Science Foundation, by the FWF grant
P24572-N25, and by an ERC GEMIS grant.  Tony Pantev was partially
supported by NSF Research Training Group Grant DMS-0636606, and NSF
research grants DMS-1001693 and DMS-1302242.

\section{Moduli of Landau-Ginzburg models} 
\label{sec:moduliLG}

In this section we study variations of pure \nc \ Hodge
structures that arise from universal families of Landau-Ginzburg
models. We focus on the components of the universal variation that
encode geometric properties of the Landau-Ginzburg moduli space and
investigate the Hodge theoretic input in the Landau-Ginzburg
deformation theory. 

The relevant class of Landau-Ginzburg models appears naturally in the
context of mirror symmetry for Fano manifolds. Since this context is a
primary source of examples for us, we recall
it next. 

\subsection{Mirrors of Fano manifolds} \label{ss-msFano}

Mirror symmetry is a duality that identifies seemingly different two
dimensional supersymmetric quantum field theories. Geometrically such
theories arise as linear or non-linear sigma models with K\"{a}hler
targets, or as Landau-Ginzburg models with targets given by K\"{a}hler
manifolds equipped with holomorphic superpotentials. The mirror map
matches the target geometries that produce mirror symmetric models
into mirror pairs. Typically a sigma model or a Landau-Ginzburg model
with a given target geometry admits two topological twists - the {\sf
  A} and {\sf B} twists - each of which gives rise to a category of
boundary field theories or $D$-branes (see
e.g. \cite{mirrorbook,mirrorbookD}). According to the homological
mirror conjecture from \cite{kontsevich-icm}, the mirror
correspondence can be generalized to an identification of the
categories of boundary field theories. Specifically the homological
mirror conjecture predicts that in a mirror pair, the category of {\sf
  A}-branes for one side of the pair must be equivalent to the
category of {\sf B}-branes for the other side. Such an equivalence
induces non-obvious isomorphisms between the various invariants that
one can extract from the categories. In particular we get a
conjectural matching of the cohomology of the two categories; matching
of the \nc \ Hodge structures on the cohomology of the two categories;
matching of the deformation spaces of the two categories; and matching
of the natural variations of \nc \ motives over these deformation
spaces. We will exploit these conjectural identifications to deduce
interesting predictions for the properties of the moduli spaces and
the Hodge theory of the requisite geometric backgrounds and will
eventually prove these predictions directly. To set things up we begin
by recalling the basic geometric framework for the mirror
correspondence.

We will indicate that two geometries $(X, \cdots)$ and $(Y, \cdots)$
are mirror equivalent by writing $(X, \cdots)\bmid (Y, \cdots)$.
Mirror pairs of geometries fall naturally into three classes: mirror
pairs of Calabi-Yau, Fano, and general type.  Here we will discuss in
detail only the mirror pairs of Fano type.

\

\noindent
By definition a {\it\bfseries mirror pair of Fano
  type} is a pair 
\[
\left(X,\omega_{X},s_{X}\right) \bmid
\left((Y,\bw),\omega_{Y},\vol_{Y}\right) 
\]
where:

\quad $\bullet$ \  $X$ is a projective Fano manifold;

\quad $\bullet$ \   $(Y,\bw)$ is a holomorphic
Landau-Ginzburg model consisting of a quasi-projective Calabi-Yau
manifold $Y$ with $\dim_{\mathbb{C}} Y = \dim_{\mathbb{C}} X = n$, and
a surjective algebraic function $\bw : Y \to \mathbb{A}^{1}$ with a compact
critical locus $\op{crit}(\bw) \subset Y$;

\quad $\bullet$ \     $\omega_{X} \in
A^{2}_{\mathbb{C}}(X)$ and $\omega_{Y} \in A^{2}_{\mathbb{C}}(Y)$ are
(complexified) K\"{a}hler forms on $X$ and $Y$;

\quad $\bullet$ \ $s_{X} \in H^{0}(X,K_{X}^{-1})$ is an anti-canonical
section of $X$, and $\vol_{Y} \in H^{0}(Y,K_{Y})$ is a trivialization
of the canonical bundle of $Y$, i.e. a holomorphic volume form on $Y$;

The anti-canonical section $s_{X} \in H^{0}(X,K_{X}^{-1})$ defines a 
Calabi-Yau hypersurface  $D_{X} = \op{divisor}(s_{X}) \subset X$ and a
nowhere vanishing section $s_{X|X-D_{X}} \in H^{0}(X-D_{X},
K_{X}^{-1})$. We will write $\vol_{X-D_{X}} = 1/s_{X}$ for the
corresponding holomorphic volume form on $X-D_{X}$.

\begin{rem} \label{rem:msFano.refine} Mirror pairs of Fano type can be
  qualified/refined in different ways:
\begin{itemize}
\item[(i)] Requiring that $D_{X}$ is smooth is mirrored by the requirement
  that $\bw$ is proper.
\item[(ii)] Requiring that $D_{X}$ has strict normal crossings is mirrored
  by the requirement that the fibers of $\bw$ are Zariski open subsets
  in projective $(n-1)$-dimensional Calabi-Yau manifolds.
\end{itemize}
\end{rem}

\noindent
It is helpful to examine the shape of the geometry of a mirror pair in
examples. Many explicit and detailed descriptions of mirror pairs of
Fano type are discussed in e.g.
\cite{givental,hori-vafa,
auroux-katzarkov-orlov,auroux-katzarkov-orlov2,
abouzaid-toric}. Here we just briefly recall
Givental's picture \cite{givental} of mirrors of projective spaces.

\begin{ex} \label{ex:mirrorPn}  The first instance of a Fano type
  mirror pair was described by Givental \cite{givental}. In the most basic
  setting $X = \mathbb{P}^{n}$ is a projective space with homogeneous
  coordinates $u_{0}$, \ldots, $u_{n}$, $\omega_{X}$ is the Fubini-Studi form,
  $D_{X} \subset \mathbb{P}^{n}$ is the union of the $(n+1)$
  coordinate hyperplanes, and $s_{X}$ is given by the product of the
  homogeneous coordinate functions. On the
  mirror side $Y = \left(\mathbb{C}^{\times}\right)^{n}$ is an
  $n$-dimensional affine torus with coordinates $z_{1}$, \ldots,
  $z_{n}$, the potential  $\bw : Y \to \mathbb{A}^{1}$ is given by
\[
\bw(z_{1},\ldots,z_{n}) = \sum_{i = 1}^{n} z_{i} +
\frac{1}{z_{1}\cdots z_{n}},
\]
the symplectic form is 
\[
\omega_{Y} = \sum_{i = 1}^{n} \frac{1}{|z_{i}|^{2}} dz_{i}\wedge d\bar{z}_{i},
\]
the point $a$ is the the point at infinity, i.e. $a = \infty$, and the
holomorphic volume form is 
\[
\vol_{Y} = \bigwedge_{i = 1}^{n} \frac{dz_{i}}{z_{i}}. 
\]
If we change the setting so that on the left hand side of the pair
$D_{X}$ is not the 
toric divisor of $\mathbb{P}^{n}$ but rather is a smooth Calabi-Yau
hypersurface, the mirror $Y$ is a partial compactification of the
torus so that $\bw$ becomes a proper map with $n-1$ dimensional
Calabi-Yau fibers. Accordingly the symplectic form $\omega_{Y}$ and
holomorphic volume form $\vol_{Y}$ have to be extended to
the compactification.
\end{ex}

\

\noindent
The mirror correspondence gives a non-trivial matching
\cite{mirrorbook} of the various ingredients of the mirror pair. The
complexified K\"{a}hler structure $\omega_{X}$ is identified with a
combination of the complex structure on $Y$, the potential $\bw$, and
the volume form $\vol_{Y}$. In the other direction $\omega_{Y}$ is
identified with a combination of the complex structure on $X$ and the
section $s_{X}$ \cite{hori-vafa,mirrorbook}.

A Fano type mirror pair  gives rise to a pair of mirror  non-compact
Calabi-Yau manifolds:
\[
\left( X-D_{X}, \omega_{X|X-D_{X}}, \vol_{X-D_{X}}\right) \bmid
\left( Y, \omega_{Y}, \vol_{Y}\right). 
\]
Under a convergence assumption on the quantum product on $X$ the
category of {\sf A}-branes for the background $(X,\omega_{X},s_{X})$
can be identified with the Fukaya category $\Fuk(X,\omega_{X})$ of the
symplectic manifold underlying the Fano variety $X$
\cite{mirrorbookD}.  $\Fuk(X,\omega_{X})$ is a ${\mathbb C}$-linear
$A_{\infty}$ category which is only ${\mathbb Z}/2$-graded
\cite{fooo,fooo2}.  The category of ${\sf B}$-branes associated with
$(X,\omega_{X},s_{X})$ is identified with a dg enhancement of the
bounded derived category $D^{b}(X)$ of coherent sheaves on $X$. We
will write $\sD^{b}(X)$ for this $\mathbb{Z}$-graded ${\mathbb
  C}$-linear dg category. There are many choices for $\sD^{b}(X)$,
e.g. the homotopy category of complexes of injective
$\mathcal{O}_{X}$-modules with coherent cohomology, or Block's
category of graded $C^{\infty}$ complex vector bundles on $X$ with
$(0,\bullet)$ superconnections \cite{block}. By a theorem of Lunts and
Orlov \cite{lunts-orlov} all dg enhancements of $D^{b}(X)$ are
quasi-equivalent so one can work with any of those enhancements. By
definition $\sD^{b}(X)$ depends only on the complex structure on $X$
and is independent of the complexified K\"{a}hler structure
$\omega_{X}$ or the section $s_{X}$.  Both $D_{X}$ and $s_{X}$ are of
course essential for defining the associated Calabi-Yau pair and its
categories of branes. 

On the right hand side of the mirror pair the definition of the
categories of branes is modified to incorporate the potential
$\bw$. The category of {\sf A}-branes associated with the background
$\left( (Y,\bw), \omega_{Y}, \vol_{Y}\right)$ is the Fukaya-Seidel
category $\FS\left( (Y,\bw), \omega_{Y}, \vol_{Y}\right)$
\cite{seidel-book} and the category of {\sf B}-branes for $\left(
  (Y,\bw), \omega_{Y}, \vol_{Y}\right)$ is defined as the category
$\MF(Y,\bw)$ of matrix factorizations of the holomorphic function $\bw
: X \to \mathbb{A}^{1}$
\cite{orlov,kkp,orlov-global,lin-pomerleano,preygel,positselski-mf}. By
construction $\FS\left( (Y,\bw), \omega_{Y}, \vol_{Y}\right)$ is a
$\mathbb{C}$-linear $\mathbb{Z}$-graded $A_{\infty}$ category. Again
the $\mathbb{Z}/2$-folding of $\FS\left( (Y,\bw), \omega_{Y},
  \vol_{Y}\right)$ depends only on the $C^{\infty}$ manifold
underlying $Y$, on the function $\bw$, and on the complexified
symplectic structure $\omega_{Y}$ while the $\mathbb{Z}$-graded
version $\FS\left( (Y,\bw), \omega_{Y}, \vol_{Y}\right)$ depends also
on $\vol_{Y}$ viewed as a $C^{\infty}$ form on $Y$. Similarly
$\MF(Y,\bw)$ is a d$(\mathbb{Z}/2)$graded ${\mathbb C}$-linear
category which depends only on the complex structure of $Y$ and on the
holomorphic function $\bw$.

Homological mirror symmetry now predicts several conjectural
equivalences of categories of branes for the Fano mirror pair and for
the associated Calabi-Yau mirror pair. These equivalences are
summarized in Table~\ref{table:hms.fano}. In this table $\sD^{b}_{c}$
denotes a dg enhancements of the derived categories of coherent
sheaves with compact support and $\Fuk^{\op{wr}}$ denotes the wrapped
version of the Fukaya category \cite{abouzaid.seidel-viterbo,
  abouzaid-wrapped}. Additionally, our convention is that whenever the
notation for a Fukaya or a Fukaya-Seidel category includes a
holomorphic volume form, the objects of this category are graded spin
Lagrangians or Lagrangian thimbles, and so the category is
$\mathbb{Z}$-graded. In particular, aside from $\Fuk(X,\omega_{X})$
and $\MF(Y,\bw)$ all categories appearing in
Table~\ref{table:hms.fano} are $\mathbb{Z}$-graded.

\begin{table}[ht] 
\begin{center}
\begin{tabular}{|c|c|}
\hline 
{\sf A}-branes & {\sf B}-branes \\ \hline \hline
\begin{minipage}[c]{0.8in}

\

$\Fuk(X,\omega_{X})$ 

\

\end{minipage}
& $\MF(Y,\bw)$ \\
\hline
\begin{minipage}[c]{2.05in}

\

$\Fuk^{\op{wr}}(X-D_{X},\omega_{X},\vol_{X-D_{X}})$ 

\

\end{minipage}
&
$\sD^{b}(Y)$ \\
\hline
\begin{minipage}[c]{2.05in}

\

$\Fuk(X-D_{X},\omega_{X},\vol_{X-D_{X}})$ 

\

\end{minipage}
&
$\sD^{b}_{c}(Y)$ \\
\hline
\multicolumn{2}{c}{$\xymatrix@1@C+2pc{ \ar@{<->}@/_1pc/[rrr]_{\op{HMS}} &
 &   & }$}
\end{tabular}
\quad
\begin{tabular}{|c|c|}
\hline 
{\sf B}-branes & {\sf A}-branes \\ \hline \hline
$\sD^{b}(X)$  & 
\begin{minipage}[c]{1.35in}

\

$\FS((Y,\bw),\omega_{Y},\vol_{Y})$ 

\

\end{minipage}
\\ \hline
$\sD^{b}(X-D_{X})$ & 
\begin{minipage}[c]{1.35in}

\

$\Fuk^{\op{wr}}(Y,\omega_{Y},\vol_{Y})$ 

\

\end{minipage}
\\ \hline
$\sD^{b}_{c}(X-D_{X})$ & 
\begin{minipage}[c]{1.35in}

\

$\Fuk(Y,\omega_{Y},\vol_{Y})$ 

\

\end{minipage}
\\ \hline
\multicolumn{2}{c}{$\xymatrix@1@C+2pc{ \ar@{<->}@/_1pc/[rr]_{\op{HMS}} &
    & }$}
\end{tabular}
\end{center}
\caption[]{ 
\begin{minipage}[t]{3.5in}
\addtolength{\baselineskip}{-2pc}
Homological mirror symmetry for a mirror pair \\
\centerline{$\left(X,\omega_{X},s_{X}\right) \bmid
  \left((Y,\bw),\omega_{Y},\vol_{Y}\right)$} \\
of Fano type
\end{minipage}
}
\label{table:hms.fano}
\end{table} 

\

\begin{rem} \label{rem:moreFano.hms}  {\bfseries (i)}  \ Homological
  mirror symmetry predicts that the equivalence 
\[
\Fuk(X,\omega_{X}) \cong \MF(Y,\bw)
\]
of ${\mathbb Z}/2$-graded $A_{\infty}$ categories in
Table~\ref{table:hms.fano} will respect the
natural additional structures on these categories of branes. In
particular mirror symmetry will respect the natural decompositions of
these categories. It is known from the work of Orlov
\cite{orlov,orlov-lagrange,orlov-global} that the category of matrix
factorizations decomposes 
\[
\MF(Y,\bw) = \coprod_{\substack{\lambda \in \mathbb{A}^{1} \\[+0.3pc]
    \text{a crtitical value} \\[+0.3pc] \text{of } \bw}}
\sD^{b}_{\op{sing}}(Y_{\lambda})
\]
into a sum of  categories of singularities of the
singular fibers of $\bw$.

Similarly (see \cite{kkp}) 
the Fukaya category of the Fano manifold $X$ decomposes 
\[
\Fuk(X,\omega_{X}) = \coprod_{\substack{\lambda \in \mathbb{C}\\[+0.3pc]
    \text{an eigenvalue} \\[+0.3pc] \text{of }
    c_{1}(T_{X})*_{1}(\bullet)}} \Fuk(X,\omega_{X})_{\lambda}
\]
corresponding to the eigenvalues of quantum multiplication\footnote{By
  assumption we are working here with a convergent version of the
  quantum product $*_{q}$. The value of $q$ corresponding to the
  particular complex structure on $Y$ under the mirrror map is
  normalized in the flat coordinates to be $q = 1$. This is why we use
  the $*_{1}$ quantum product in the decomposition above.} with
$c_{1}(T_{X})$ on $H^{\bullet}(X,\mathbb{C})$.

\

\noindent {\bfseries (ii)} \ The equivalences of the categories of
branes listed in Table~\ref{table:hms.fano} induce respective mirror
identifications of cohomology groups. For future reference we
collect these identifications in Table~\ref{table:hms.cohomology} below.
\begin{table}[ht] 
\begin{center}
\begin{tabular}{|c|c|}
\hline 
\begin{minipage}[c]{0.8in}
\addtolength{\baselineskip}{-0.8pc}

\

\

{\sf A}-brane

charges

\

\end{minipage}
 & 
\begin{minipage}[c]{0.8in}
\addtolength{\baselineskip}{-0.8pc}

\

{\sf B}-brane

charges

\

\end{minipage}
 \\ \hline \hline
\begin{minipage}[c]{0.8in}

\

$H^{\bullet}\left(X,\mathbb{C}\right)$ 

\

\end{minipage}
& $H^{\bullet}\left(Y,Y_{\op{sm}};\mathbb{C}\right)$ \\
\hline
\begin{minipage}[c]{1.2in}

\

$H^{\bullet}\left(X-D_{X},\mathbb{C}\right)$ 

\

\end{minipage}
&
$H^{\bullet}(Y,\mathbb{C})$ \\
\hline
\begin{minipage}[c]{1.2in}

\

$H^{\bullet}_{c}\left(X-D_{X},\mathbb{C}\right)$ 

\

\end{minipage}
&
$H^{\bullet}_{c}(Y,\mathbb{C})$ \\
\hline
\multicolumn{2}{c}{$\xymatrix@1@C+2pc{ \ar@{<->}@/_1pc/[rrr]_{\op{HMS}} &
 &   & }$}
\end{tabular}
\quad
\begin{tabular}{|c|c|}
\hline 
\begin{minipage}[c]{0.8in}
\addtolength{\baselineskip}{-0.8pc}

\

\

{\sf B}-brane

charges

\

\end{minipage}
 & 
\begin{minipage}[c]{0.8in}
\addtolength{\baselineskip}{-0.8pc}

\

\

{\sf A}-brane

charges

\

\end{minipage}
 \\ \hline \hline
$H^{\bullet}(X,\mathbb{C})$  & 
\begin{minipage}[c]{1in}

\

$H^{\bullet}(Y,Y_{\op{sm}};\mathbb{C})$

\

\end{minipage}
\\ \hline
$H^{\bullet}(X-D_{X},\mathbb{C})$ & 
\begin{minipage}[c]{1in}

\

$H^{\bullet}(Y,\mathbb{C})$ 

\

\end{minipage}
\\ \hline
$H^{\bullet}_{c}(X-D_{X},\mathbb{C})$ & 
\begin{minipage}[c]{1in}

\

$H^{\bullet}_{c}(Y,\mathbb{C})$ 

\

\end{minipage}
\\ \hline
\multicolumn{2}{c}{$\xymatrix@1@C+2pc{ \ar@{<->}@/_1pc/[rr]_{\op{HMS}} &
    & }$}
\end{tabular}
\end{center}
\caption[]{ 
\begin{minipage}[t]{3.5in}
\addtolength{\baselineskip}{-2pc}
Matching of cohomology for a mirror pair \\
\centerline{$\left(X,\omega_{X},s_{X}\right) \bmid
  \left((Y,\bw),\omega_{Y},\vol_{Y}\right)$} \\
of Fano type
\end{minipage}
}
\label{table:hms.cohomology}
\end{table} 
In this table $Y_{\op{sm}} \subset Y$ denotes a smooth fiber of
$\bw : Y \to \mathbb{A}^{1}$ taken ``near infinity'' as explained in 
\cite[Section~4.5.2(2)]{kkp}. 

\

\noindent {\bfseries (iii)} \ The {\sf B}-to-{\sf A} homological
mirror correspondence in Table~\ref{table:hms.fano} can be extended to
one more case. Let $Y_{-\infty}$ denote the fiber $\bw^{-1}(z)$ over
$z \in \mathbb{C}$ with $\op{Re} z \ll 0$. We will also write
$\omega_{-\infty}$ for the restriction $\omega_{Y|Y_{-\infty}}$ of the
symplectic form, and $\vol_{-\infty}$ for the induced holomorphic
volume form on the fiber. The parallel transport for the Erhesmann
symplectic connection on $\bw : Y \to \mathbb{A}^{1}$.  identifies
symplectically all fibers of $\bw$ over points $z \in \mathbb{A}^{1}$
with  $\op{Re} z \ll 0$. So the dg category
$\Fuk(Y_{-\infty},\omega_{-\infty},\vol_{-\infty})$ is well defined up
to quasi-equivalence. Now, we can supplement
Table~\ref{table:hms.fano} by the statement that the category
 of perfect complexes (=  the category of topological {\sf
  B} branes) on the Calabi-Yau variety $D_{X}$ is quasi-equivalent 
to the Fukaya
 category  ( = the category o f {\sf
  A}-branes) on the fiber $Y_{-\infty}$:

\

\begin{center}
\begin{tabular}{|c|c|}
\hline 
{\sf B}-branes & {\sf A}-branes \\ \hline \hline
$\Perf\left(D_{X}\right)$  & 
\begin{minipage}[c]{1.6in}

\

$\Fuk\left(Y_{-\infty},\omega_{-\infty},\vol_{-\infty}\right)$ 

\

\end{minipage}
\\ \hline
\multicolumn{2}{c}{$\xymatrix@1@C+2pc{ \ar@{<->}@/_1pc/[rr]_{\op{HMS}} &
    & }$}
\end{tabular}
\end{center}

\

\noindent
Again this induces an identification of the associated brane charges,
i.e. of the periodic cyclic homologies of the two categories. Since
$HP_{\bullet}\left( \Perf(D_{X})\right) \cong H^{\bullet}(D_{X},\mathbb{C})$
and conjecturally $HP_{\bullet}\left( 
\Fuk\left(Y_{-\infty},\omega_{-\infty},\vol_{-\infty}\right) \right)
\cong H^{\bullet}_{c}(Y_{-\infty}, \mathbb{C})$ we get a mirror
identification

\

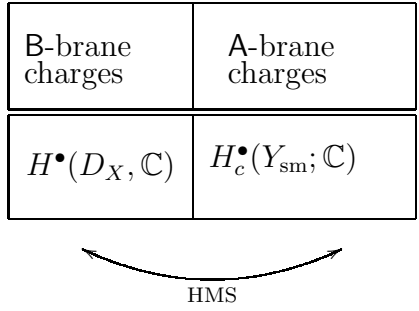
\begin{center}
\begin{tabular}{|c|c|}
\hline 
\begin{minipage}[c]{0.8in}
\addtolength{\baselineskip}{-0.8pc}

\

\

{\sf B}-brane

charges

\

\end{minipage}
 & 
\begin{minipage}[c]{0.8in}
\addtolength{\baselineskip}{-0.8pc}

\

\

{\sf A}-brane

charges

\

\end{minipage}
 \\ \hline \hline
$H^{\bullet}(D_{X},\mathbb{C})$  & 
\begin{minipage}[c]{1in}

\

$H^{\bullet}_{c}(Y_{\op{sm}};\mathbb{C})$

\

\end{minipage}
\\ \hline
\multicolumn{2}{c}{$\xymatrix@1@C+2pc{ \ar@{<->}@/_1pc/[rr]_{\op{HMS}} &
    & }$}
\end{tabular}
\end{center}

\

\noindent
Note that it is not clear how to extend the {\sf A}-to-{\sf B}
homological mirror correspondence in a similar manner. If on the {\sf
  A}-side we consider the symplectic data $\left(D_{X},
\omega_{X|D_{X}}, \vol_{D_{X}} \right)$, there is no obvious complex
fiber $Y_{c}$ of $\bw : Y \to \mathbb{A}^{1}$ for which we can hope to
get a quasi-equivalence $\Fuk\left(D_{X}, \omega_{X|D_{X}},
\vol_{D_{X}} \right) \cong \sD^{b}(Y_{c})$. The problem is that the
mirror of $\left(D_{X}, \omega_{X|D_{X}}, \vol_{D_{X}} \right)$ is
normally understood in terms of the large volume degenration, so for
an {\sf A}-to-{\sf B} mirror statemnt we will need to understand the  large
complex structure degeneration of $\bw : Y \to \mathbb{A}^{1}$ . 
\end{rem}

\subsection{Families of Landau-Ginzburg models} 
\label{ss-familiesLG}

As we explained in Section~\ref{ss-msFano} the mirror of a symplectic
manifold underlying a Fano variety is a quasi-projective
Landau-Ginzburg model $((Y,\bw),\vol_{Y})$ equipped with a holomorphic
volume form. Such Landau-Ginzburg models admit a natural class of
compactifications. 

\begin{defi} \label{defi:compactLG}
A {\bfseries compactified Landau-Ginzburg model} is  the datum
$\left((\Ybar, \wbar), D_{\Ybar}, \vol_{\Ybar}\right)$, where: 
\begin{itemize}
\item[(a)]  $\Ybar$ is a smooth projective variety and  $\wbar : \Ybar \to
\mathbb{P}^{1}$ is a flat projective morphism. 
\item[(b)] $D_{\Ybar} = \left(\cup_{i} \hor{D}_{i}\right)\cup \left(\cup_{j}
    \ver{D}_{j}\right)  \subset
  \Ybar$ is a reduced normal crossings divisor, such that
\begin{itemize}
\item $\ver{D} = \cup_{j} \ver{D}_{j}$ is the reduced pole divisor of
  $\wbar$, i.e. $\left(\wbar^{-1}(\infty)\right)_{\op{red}} = \cup_{j}
  \ver{D}_{j}$;
\item each component $\hor{D}_{i}$ of $\hor{D} = \cup_{i} \hor{D}_{i}$
  is a smooth divisor which is horizontal
for $\bw$, i.e. $\bw_{|\hor{D}_{i}}$ is a
flat projective morphism;
\item the critical locus $\op{crit}(\wbar) \subset \Ybar$ does not
  intersect $\hor{D}$.
\end{itemize}
\item[(c)] $\vol_{\Ybar}$ is a meromorphic section of $K_{\Ybar}$ with no
  zeroes and with poles at most at $D_{\Ybar}$, i.e. $\vol_{\Ybar} \in
  H^{0}\left(\Ybar, K_{\Ybar}(*D_{\Ybar})\right)$.
\end{itemize}
\end{defi}

\

\noindent
With every $\left((\Ybar,\wbar),D_{\Ybar},\vol_{\Ybar}\right)$ we 
associate its 'open part' $\left((Y,\bw),\vol_{Y}\right)$ where
$Y := \Ybar - D_{\Ybar}$, $\bw : Y \to \mathbb{A}^{1}$ is defined to
be the restriction $\bw := \wbar_{|Y}$ of $\wbar$ to $Y$, and
$\vol_{Y} := \vol_{\Ybar|Y}$ is the restriction of $\vol_{\Ybar}$ to $Y$.
The condition on $\vol_{\Ybar}$ ensures that $Y$ is a quasi-projective
variety with a trivial canonical class and that $\vol_{Y}$ is a
holomorphic volume form on $Y$. The condition that the critical locus
of $\wbar$ does not intersect the horizontal part of $D_{\Ybar}$
in particular implies that $\op{crit}(\bw)$ is proper and so
$\left((Y,\bw),\vol_{Y}\right)$ is exactly the type of Landau-Ginzburg
model that we considered in the previous section. 

\

\noindent
In addition we will often require that the datum
$\left((\Ybar,\wbar),D_{\Ybar},\vol_{\Ybar}\right)$ satisfies the
following tameness assumption which bounds the orders of poles of
$\vol_{\Ybar}$ and $\wbar$ along $D_{\Ybar}$:
\[ \label{eq:tame}
\fbox{
\tag{{\sf T}} \hspace{-1pc}
\text{\begin{tabular}{rcrrcl}
$\op{ord}_{\hor{D}_{i}}\left(\vol_{\Ybar}\right)$ & $=$ & $-1$, \qquad\qquad & 
$\op{ord}_{\ver{D}_{j}}\left(\vol_{\Ybar}\right)$ & $=$ & $-1$ \\
$\op{ord}_{\hor{D}_{i}}\left(\wbar\right)$ & $=$ & $0$, \qquad\qquad &
  $\op{ord}_{\ver{D}_{j}}\left(\wbar\right)$ & $=$ & $-1$ 
\end{tabular}}
\qquad \text{ for all $i$ and $j$}.}
\]
\

\

\medskip

\begin{rem} \label{rem:tameness}
{\bfseries (i)} \ The assumption that $\vol_{\Ybar}$ has poles of
order exactly $1$ along all the components of $D_{\Ybar}$ in
particular implies that the reduced divisor $D_{\Ybar}$ is an
anti-canonical divisor on $\Ybar$ and so $(\Ybar,D_{\Ybar})$ is a log
Calabi-Yau pair. 

Note also that if we start with a Calabi-Yau quasi-projective
Landau-Ginzburg model $\left((Y,\bw),\vol_{Y}\right)$ and we choose a
smooth normal-crossing compactification $\wbar : \Ybar \to
\mathbb{P}^{1}$ of \linebreak $\bw : Y \to \mathbb{A}^{1}$, then the
condition that the holomorphic volume form has a first order pole
along the divisor at infinity is a tight constraint which is rather
unnatural since it is not invariant under semi-stable
reduction. Nevertheless this condition is often satisfied in mirror
symmetry examples and so it is not unreasonable to impose.

\

\noindent {\bfseries (ii)} The assumption that $\bw$ has first order
poles at the vertical boundary divivsor $\ver{D}$, i.e. that the
scheme theoretic fiber $\wbar^{-1}(\infty)$ is actually reduced, is
more natural and can be justified by mirror symmetry consderations. 

Indeed, if $\wbar^{-1}(\infty) = \sum_{j} m_{j} \ver{D}_{j}$, then by
Landman's theorem \cite{landman} we know that the least common
multiple $m$ of the $m_{j}$'s is the order of the semi-simple part of
the local monodromy transformation around infinity. Concretely, choose
a small disk $\Delta \subset \mathbb{P}^{1}$ centered at $\infty \in
\mathbb{P}^{1}$ and such that $\infty$ is the only critical value of
$\wbar$ in $\Delta$.  Fix a base point $c_{0} \in \partial \Delta$,
and orient $\partial \Delta$ with the orientation on $\Delta$.
Consider the monodromy transformation $\op{{\sf mon}}_{c_{0}} :
H^{\bullet}(\Ybar_{c_{0}},\mathbb{C}) \to
H^{\bullet}(\Ybar_{c_{0}},\mathbb{C})$ corresponding to going around
$\partial \Delta$ once in the positive direction . By Landman's
theorem $\op{{\sf mon}}_{c_{0}}$ is a quasi-unipotent operator and the
minimal power of $\op{{\sf mon}}_{c_{0}}$ which is unipotent is
$m$. In other words $\left( \op{{\sf mon}}_{c_{0}}^{m} -
  \op{id}\right)^{n-1} = 0$, where $n = \dim_{\mathbb{C}} Y$, and $m$
is the minimal number with this property. Similarly, going once around
$\partial \Delta$ gives a monodromy transformation $T :
H^{\bullet}(Y,Y_{c_{0}};\mathbb{C}) \to
H^{\bullet}(Y,Y_{c_{0}};\mathbb{C})$. From our asumption on
equi-singularity of $\hor{D}$ and from the compatibility of the long
exact sequence of the pair with the action of monodromy we get that
$T$ will also be quasi-unipotent with eigenvalues which are $m$-th
roots of unity with at least one eigenvalue being a primitive $m$-th
root of unity. Next observe that the cycle class map that assigns a
relative cohomology class in $H^{\bullet}(Y,Y_{c_{0}};\mathbb{C})$ to
each Lefschetz thimble will identify the periodic cyclic homology
$HP_{\bullet}(\FS((Y,\bw),\omega_{Y},\vol_{Y}))$ of the Fukaya-Seidel
category with $H^{\bullet}(Y,Y_{c_{0}};\mathbb{C})$. From this point
of view the operator $T$ is induced from the inverse of the monodromy
auto-equivalence of $\FS((Y,\bw),\omega_{Y},\vol_{Y})$. When
$((Y,\bw),\omega_{Y},\vol_{Y}))$ is the mirror of a Fano datum
$(X,\omega_{X},s_{X})$, the mirror equivalence (see
Table~\ref{table:hms.fano}) $\FS((Y,\bw),\omega_{Y},\vol_{Y}) \cong
\sD^{b}(X)$ identifies the monodromy auto-equivalence with the Serre
functor $\otimes K_{X}[n] : \sD^{b}(X) \to \sD^{b}(X)$. But on
cohomology $H^{\bullet}(Y,Y_{c_{0}};\mathbb{C}) \cong
H^{\bullet}(X,\mathbb{C})$ the Serre functor induces multiplication
with $\exp\left((-1)^{n}c_{1}(K_{X})\right)$. In other words $T$ is a
unipotent operator, and so we must have $m = 1$. For future reference
note that this mirror symmetry description also predicts that under
the identification $H^{\bullet}(X,\mathbb{C}) \cong
H^{\bullet}(Y,Y_{c_{0}};\mathbb{C})$ the nilpotent endomorphism 
\[
\left( (-1)^{n}c_{1}(K_{X}) \right)\cup(\bullet) :
H^{\bullet}(X,\mathbb{C}) \to H^{\bullet}(X,\mathbb{C})
\]
becomes identified with the logarithm of monodromy:
\[
- \log T :  H^{\bullet}(Y,Y_{c_{0}};\mathbb{C}) \to
H^{\bullet}(Y,Y_{c_{0}};\mathbb{C}).
\]
\item[(iii)] It is possible and useful to allow for $\ver{D}$ to be
  a smooth divisor, i.e. for $\infty$ not to be a critical value of
  $\wbar$. Such Landau-Ginzburg models arise naturally as mirrors of
  quasi-Fano varieties and can be studied in the same manner. 
\end{rem}

\

\noindent
Our goal is to understand the moduli spaces of compactified complex
Landau-Ginzburg models satisfying the tameness asumption. If such a
model $((\Ybar,\wbar),D_{\Ybar},\vol_{\Ybar})$ compactifies the mirror
of a Fano datum $(X,\omega_{X},s_{X})$, then its moduli space will be
identified with the symplectic moduli of $(X,\omega_{X})$ and in fact
will look like a conical  open subset in $H^{2}(X,\mathbb{C})\oplus
H^{0}(X,\mathbb{C})$. In particular, when
$((\Ybar,\wbar),D_{\Ybar},\vol_{\Ybar})$ arises from a mirror
situation, we expect its moduli space to be smooth. This motivates the
following purely algebraic-geometric statement

\begin{theo} \label{theo:mainTT} Let
  $((\Ybar,\wbar),D_{\Ybar},\vol_{\Ybar})$ be a compactified
  Landau-Ginzburg model satisfying the tameness assumption
  \eqref{eq:tame} and the assumption
  $H^{1}(\Ybar,\mathbb{Q}) = 0$. Then the deformation theory of
  $((\Ybar,\wbar)_{D_{\Ybar}},\vol_{\Ybar})$ is unobstructed.
\end{theo}

\

\begin{rem} \label{rem:teich} The requirement that
  $H^{1}(\Ybar,\mathbb{Q}) = 0$ is a technical requirement that
  simplifies the Teichm\"{u}ller theory of $\Ybar$. It is very likely
  unnecessary  but we will not pursue this here. 
\end{rem}

\

\noindent
Theorem~\ref{theo:mainTT} extends the classical unobstructedness results of
Bogomolov \cite{fedya-unobstructed,fedya-ihes}, Tian \cite{tian},
and Todorov \cite{andrey} to the setup of log Calabi-Yau
varieties with potentials. It also gives the following more direct
generalization of Bogomolov-Tian-Todorov unobstructedness:

\begin{cor} \label{cor:rat.curves}
Let $\Ybar$ be a smooth projective variety of dimension $n$ satisfying
$H^{1}(\Ybar,\mathbb{Q}) = 0$  and such that the anti-canonical linear
system on $\Ybar$ gives a flat projective morphism \linebreak $\wbar :
\Ybar \to \mathbb{P}^{1}$. Then the deformation theory of $\Ybar$ is
unobstructed.  
\end{cor}
{\bfseries Proof.} By assumption the variety $\Ybar$ determines the
morphism $\wbar : \Ybar \to \mathbb{P}^{1}$. Let $D_{\Ybar} \subset
\Ybar$ be a smooth anti-canonical divisor in $\Ybar$.  and let
$\vol_{\Ybar}$ be a trivialization of $K_{\Ybar}(D_{\Ybar})$. Then the
datum $((\Ybar,\wbar),D_{\Ybar},\vol_{\Ybar})$ is a tame compactified
Landau-Ginzburg model satisfying the hypotheses of
Theorem~\ref{theo:mainTT}. Locally in the analytic topology the versal
deformation space of the datum
$((\Ybar,\wbar),D_{\Ybar},\vol_{\Ybar})$ is the product of the versal
deformation space of $\Ybar$ and the moduli of pairs
$(D_{\Ybar},\vol_{\Ybar})$ for a fixed $\Ybar$. But the moduli of such
pairs is isomorphic to $(\mathbb{P}^{1}-\op{crit}(\wbar))\times
\mathbb{C}^{\times}$ and is therefore smooth. Combined with the
unobstructedness of Theorem~\ref{theo:mainTT} this implies that the
deformations of $\Ybar$ are unobstructed. \ \hfill $\Box$

\

\begin{rem}
If $\Ybar$ satisfies the conditions of Corollary~\ref{cor:rat.curves},
then the family $\wbar : \Ybar \to \mathbb{P}^{1}$ is classified by a
holomorphic map from $\mathbb{P}^{1}$ to the compactified moduli space
of $(n-1)$-dimensional projective Calabi-Yau varieties. By the
classical Bogomolov-Tian-Todorov theorem \cite{tian,andrey} we know that the
moduli space $M$ of $(n-1)$-dimensional projective Calabi-Yau
varieties is smooth. By Corollary~\ref{cor:rat.curves} we know that a
certain component $M_{1}$ of the moduli space of rational curves in
$\overline{M}$ is also smooth. Considerations of mirrors of hybrid
Landau-Ginzburg models suggest that this process can be iterated: a
component $M_{2}$ of the moduli of rational curves in a
compactification $\overline{M}_{1}$ will be smooth, and so on. It will
be very interesting to analyze this problem from purely
algebraic-geometric point of view and to construct iteratively the
sequence of $L_{\infty}$ algebras controlling the corresponding
deformation problems. 
\end{rem}

\

\noindent
Before we proceed with the proof of Theorem~\ref{theo:mainTT} we will
need to establish some general facts about the deformation theory of
varieties with potentials. 

\subsection{Deformations of compactified Landau-Ginzbug
  models} \label{ssec:defosLG} 

Let $((\Ybar,\wbar),D_{\Ybar},\vol_{\Ybar})$ be a compactified
Landau-Ginzburg model satisfying $H^{1}(\Ybar,\mathbb{Q}) = 0$ and the
tameness condition \eqref{eq:tame}. Since $H^{1}(\Ybar,\mathbb{Q}) =
0$ implies $\op{Pic}^{0}(\Ybar) = 0$ and since the Neron-Severi class
$[D_{\Ybar}] \in H^{2}(\Ybar,\mathbb{Z})$ is preserved under small
deformations of $\Ybar$, it follows that the condition $D_{\Ybar} \in
|K_{\Ybar}^{-1}|$ is also preserved under small deformations of the
pair $(\Ybar,D_{\Ybar})$. By \eqref{eq:tame} the meromorphic volume
form $\vol_{\Ybar}$ is a trivialization of the line bundle
$K_{\Ybar}(D_{\Ybar})$ and so the versal deformation space of
$((\Ybar,\wbar)_{D_{\Ybar}},\vol_{\Ybar})$ is a principal
$\mathbb{C}^{\times}$-bundle over the versal deformation space of
$(\Ybar,\wbar)_{D_{\Ybar}}$. Therefore it suffices to prove the
unobstructedness of the deformation theory of
$(\Ybar,\wbar)_{D_{\Ybar}}$.

As usual the deformation theory of $(\Ybar,\wbar)_{D_{\Ybar}}$ is
controlled by an $L_{\infty}$ algebra
\cite{ks-defos,hinich,ks-ncgeometry,lurie-defos}. By standard
Kodaira-Spenser theory the deformations of the map $\wbar : \Ybar \to
\mathbb{P}^{1}$ are computed \cite{illusie1,illusie2,horikawa,sernesi}
by the sheaf of dg Lie algebras
\[
[
\xymatrix@R-2pc@C+1pc{T_{\Ybar} \ar[r]^-{d\wbar} & \wbar^{*}
  T_{\mathbb{P}^{1}} \\
0 & 1
}
].
\]
Here $\wbar^{*}  T_{\mathbb{P}^{1}}$ denotes the $\mathcal{O}$-module
pullback, and this is an complex of locally free coherent sheaves with
an $\mathcal{O}_{\Ybar}$-linear differential and $\mathbb{C}$-bilinear
(graded) Lie bracket. 

Similarly the deformations of $\wbar : \Ybar \to \mathbb{P}^{1}$ which
preserve the boundary divisor $D_{\Ybar}$ are computed by the sheaf of
dg Lie algebras  
\[
\dtg^{\bullet} := 
[\xymatrix@R-1pc@C+1pc{
T_{\Ybar,D_{\Ybar}} \ar[r]^-{d\wbar} \ar@{=}[d] & \wbar^{*}
  T_{\mathbb{P}^{1},\infty} \ar@{=}[d] \\
\dtg^{0} & \dtg^{1} 
}],
\]
where for any smooth variety $M$ and any closed reduced subscheme 
$S \subset M$ we write $T_{M,S}$ for the coherent sheaf of vector
fields on $M$ that are tangent to $S$ at the points of $S$. Since
$D_{\Ybar} \subset \Ybar$ and $\{\infty\} \subset \mathbb{P}^{1}$ are
reduced normal crossings divisors, it follows that
$T_{\Ybar,D_{\Ybar}} \subset T_{\Ybar}$ and $T_{\mathbb{P}^{1},\infty}
\subset T_{\mathbb{P}^{1}}$ are locally free subsheaves. 

Recall (see e.g. \cite{goldman-millson,ks-defos,manetti}) that the
unobstructedness of the deformation theory defined by an $L_{\infty}$
algebra follows from the stronger property that this $L_{\infty}$
algebra is homotopy abelian. Therefore Theorem~\ref{thm:main.TT} and
Theorem~\ref{theo:mainTT} will follow immediately from the following

\begin{theo} \label{theo:ha} Suppose
  $((\Ybar,\wbar),D_{\Ybar},\vol_{\Ybar})$ is a compactified
  Landau-Ginzburg model  satisfying 
  the tameness condition \eqref{eq:tame}. Then the $L_{\infty}$
  algebra 
\[
R\Gamma(\Ybar,\dtg^{\bullet}) = 
R\Gamma\left(\Ybar,
\left[\xymatrix@1@C+1pc{
T_{\Ybar,D_{\Ybar}} \ar[r]^-{d\wbar}  & \wbar^{*}
  T_{\mathbb{P}^{1},\infty}
}
  \right]\right) 
\]
is homotopy abelian.
\end{theo}

As in the compact Calabi-Yau case we will deduce
Theorem~\ref{theo:ha} from a Hodge theoretic statement - the
degeneration of a ``Hodge-to-de Rham spectral sequence'' associated
with the divisor $D_{\Ybar}$ and the potential $\wbar$. Our main tool
here is a new complex of logarithmic forms adapted to $\wbar$:

\begin{defi} \label{defi:main} Let $((\Ybar,\wbar),D_{\Ybar})$ be a
  compactified Landau-Ginzburg model satisfying the conditions
  Definition~\ref{defi:compactLG}(a) and
  Definition~\ref{defi:compactLG}(b). For any $a \geq 0$ we define the
  {\bfseries sheaf $\Omega_{\Ybar}^{a}\left(\log
    D_{\Ybar},\wbar\right) $ of $\wbar$-adapted logarithmic forms on
    $(\Ybar,D_{\Ybar})$} as the subsheaf of logarithmic $a$-forms that
  stay logarithmic after multiplication by $d\wbar$. Thus
\[
\Omega_{\Ybar}^{a}\left(\log D_{\Ybar},\wbar\right) := \left\{ \alpha \in
\Omega_{\Ybar}^{a}\left(\log D_{\Ybar}\right) \, \left| \, d\wbar\wedge
\alpha \in \Omega_{\Ybar}^{a+1}\left(\log D_{\Ybar}\right)\right.\right\}
\subset \Omega_{\Ybar}^{a}\left(\log D_{\Ybar}\right),
\]
where $\wbar$ is viewed as a meromorphic function on $\Ybar$ and
$d\wbar$ is viewed as a meromorphic one form. 
\end{defi}

\

\noindent
The sheaves $\Omega_{\Ybar}^{a}\left(\log D_{\Ybar},\wbar\right)$ have
several interesting properties. As a first remark we have the
following

\begin{lemma} \label{lemma:locally.free} 
\begin{itemize}
\item[{\normalfont (a)}]
The sheaf
  $\Omega_{\Ybar}^{a}\left(\log D_{\Ybar},\wbar\right)$ of
  $\wbar$-adapted logarithmic forms is a coherent \linebreak
  $\mathcal{O}_{\Ybar}$-module which is locally free of rank equal to
  $\op{rank} \, \Omega_{\Ybar}^{a} = \binom{n}{a}$.
\item[{\normalfont (b)}] Suppose $\varepsilon : \widehat{\Ybar} \to
  \Ybar$ is a blow-up of $\Ybar$ with smooth center contained in
  $\ver{D}$ and cleanly intersecting each component of $\ver{D}$. Let
  $\widehat{D}_{\widehat{\Ybar}} = \varepsilon^{*}D_{\Ybar}$ and
  $\hat{\wbar} = \varepsilon^{*}\wbar$ denote the pullbacks of the
  divisor and potential to $\widehat{\Ybar}$. Then
  $R\varepsilon_{*}\Omega_{\widehat{\Ybar}}^{a}\left(\log
    \widehat{D}_{\widehat{\Ybar}},\hat{\wbar}\right) =
  \Omega_{\Ybar}^{a}\left(\log 
    D_{\Ybar},\wbar\right)$.
\end{itemize}
\end{lemma}
{\bfseries Proof.}  Indeed, let $\jmath_{Y} : Y \hookrightarrow \Ybar$
denote the inclusion of $Y$ in $\Ybar$. By definition
$\Omega_{\Ybar}^{a}\left(\log D_{\Ybar},\wbar\right)$ is the preimage
of the coherent $\mathcal{O}_{\Ybar}$-submodule
$\Omega_{\Ybar}^{a+1}\left(\log D_{\Ybar}\right) \subset
\jmath_{Y*}\Omega_{Y}^{a+1}$ under the $\mathcal{O}_{\Ybar}$-linear
map $d\wbar \wedge : \Omega_{\Ybar}^{a}\left(\log D_{\Ybar}\right) \to
\jmath_{Y*}\Omega_{Y}^{a+1}$. Thus $\Omega_{\Ybar}^{a}\left(\log
D_{\Ybar},\wbar\right)$ is a torsion-free coherent submodule in
$\Omega_{\Ybar}^{a}\left(\log D_{\Ybar}\right)$ of maximal rank.

The fact that $\Omega_{\Ybar}^{a}\left(\log D_{\Ybar},\wbar\right)$ is
locally free can be checked locally analytically on $\Ybar$.  

On the open set $Y \subset \Ybar$ we have by definition
$\Omega_{\Ybar}^{a}\left(\log D_{\Ybar},\wbar\right)_{|Y} =
\Omega_{Y}^{a}$ and so it is locally free. Furthermore, if $p \in
\hor{D} - \ver{D}$, then $d\wbar$ is holomorphic in a neighborhood of
$p$. This implies that near $p$ we have that
$\Omega_{\Ybar}^{a}\left(\log D_{\Ybar},\wbar\right)$ is isomorphic to
$\Omega_{\Ybar}^{a}\left(\log D_{\Ybar}\right)$ and so is locally
free.

Suppose next $p \in \ver{D}$.  We can find local analytic coordinates
  $z_{1}, \ldots, z_{n}$ centered at $p$ so that 
in a neighborhood of $p$:

\begin{itemize} 
\item the divisor $\ver{D}$ is given by
  $\prod_{i=1}^{k} z_{i}= 0$, the divisor $\hor{D}$ is given by
  $\prod_{1 = k+1}^{k+l} z_{i} = 0$;
\item the potential $\wbar$ is given by
${\displaystyle \wbar(z_{1},\ldots, z_{n}) =
  \frac{1}{z_{1}^{m_{1}}\cdots z_{k}^{m_{k}}}}$ for some $m_{i} \geq
  1$.
\end{itemize}

\

\noindent
 Now for any $a$ we have $\Omega_{\Ybar}^{a}\left(\log
 D_{\Ybar}\right) = \oplus_{p = 0}^{a} \wedge^{p} V \otimes
 \wedge^{a-p} R$, where $V \subset \Omega_{\Ybar}^{1}\left(\log
 D_{\Ybar}\right)$ is the sub $\mathcal{O}$-module spanned by $\{
 d\log z_{i}\}_{i = 1}^{k}$, while $R \subset
 \Omega_{\Ybar}^{1}\left(\log D_{\Ybar}\right)$ is the sub
 $\mathcal{O}$-module spanned by $\{d\log z_{i}\}_{i = k+1}^{k+l}$ and
 $\{d z_{i}\}_{i = k+l+1}^{n}$.

Since $d\wbar$ only has poles at the components of $\ver{D}$, the
condition that a logarithmic form $\alpha = \sum_{p} \nu_{p}\otimes
\rho_{a-p} \in \Omega_{\Ybar}^{a}\left(\log D_{\Ybar}\right)=
\oplus_{p = 0}^{a} \wedge^{p} V \otimes \wedge^{a-p} R$ is
$\wbar$-adapted will only impose constraints on the pieces $\nu_{p}
\in \wedge^{p}V$. Thus it is enough to understand which local sections
of $\wedge^{p}V$ are $\wbar$-adapted. 

Write $W \subset V$ for the sub $\mathcal{O}$-module
spanned by $\{ d\log z_{i}\}_{i = 1}^{k-1}$. In particular 
\[
V = W \oplus \mathcal{O}\cdot d\log z_{k}, \text{ and } \wedge^{p} V =
\wedge^{p}W \oplus \left(\wedge^{p-1}W\wedge d\log z_{k}\right), 
\]
and so given any $\nu \in \wedge^{p} V$, we can write $\nu$  and
$d\wbar\wedge \nu$ uniquely as
\[
\begin{split}
\nu & = \eta + \beta\wedge d\log z_{k}, \text{ with } \eta \in
\wedge^{p}W, \text{ and } \beta \in \wedge^{p-1}W, \\[+0.5pc]
d\wbar \wedge \nu & = \varphi + \psi\wedge d\log z_{k}, \text{ with }
\varphi \in 
(\wedge^{p+1}W)(*\ver{D}), \text{ and } \psi \in (\wedge^{p}W)(*\ver{D}).
\end{split}
\]
We have $d\wbar = \wbar\cdot d\log\wbar$. The logarithmic $1$-form 
$d\log\wbar$ also decomposes as \linebreak $d\log\wbar = \omega - m_{k}d\log
z_{k}$  where $\omega = - \sum_{i = 1}^{k-1} m_{i}d\log z_{i}$ is its
$W$-component. This gives
\[
\begin{split}
\varphi & = \wbar \cdot \omega\wedge\eta, \\[+0.5pc] 
\psi & = \wbar\cdot
\left( \omega\wedge\beta - (-1)^{p}m_{k}\eta\right). 
\end{split}
\]
In particular we can solve for $\eta$ in terms of $\psi$ and $\beta$. 
The condition that $\nu$ is $\wbar$-adapted is simply the condition
that $\varphi \in \wedge^{p+1}W$ and $\psi \in \wedge^{p}W$. But
for any
$\psi \in \wedge^{p}W$ and any $\beta \in \wedge^{p-1}W$ 
the form 
\[
\eta = \frac{1}{(-1)^{p}m_{k}}\cdot\left(\omega\wedge\beta -
\frac{1}{\wbar}\cdot \psi\right) 
\] 
automatically satisfies $\eta \in \wedge^{p}W$ and 
\[
\wbar \cdot \omega\wedge\eta = - \frac{1}{(-1)^{p}m_{k}} \cdot \omega\wedge
\psi \in \wedge^{p+1}W.
\]
In other words $\nu$ is $\wbar$-adapted if and only if we can find a
form $\psi \in \wedge^{p}W$ and a form $\beta \in \wedge^{p-1}W$ so that 
\[
\boxed{
\nu = \frac{1}{(-1)^{p}m_{k}}\cdot\left[d\log\wbar\wedge\beta -  
\frac{1}{\wbar}\cdot \psi\right].
}
\]
This shows that the subsheaf $\wedge^{p}V\cap
\Omega^{p}_{\Ybar}\left(\log D_{\Ybar},\wbar\right)$ in
$\wedge^{p}V$ consisting of $\wbar$-adapted forms is given by
\[
\wedge^{p}V\cap
\Omega^{p}_{\Ybar}\left(\log D_{\Ybar},\wbar\right) =
\frac{1}{\wbar}\wedge^{p}W\oplus d\log\wbar \wedge \left(\wedge^{p-1}W\right).
\]
In particular $\wedge^{p}V\cap
\Omega^{p}_{\Ybar}\left(\log D_{\Ybar},\wbar\right)$ is locally free
and hence $\Omega^{a}_{\Ybar}\left(\log D_{\Ybar},\wbar\right)$ is
locally free. Explicitly 
\begin{equation}
\label{eq:adapted.forms}
\boxed{
\Omega^{a}_{\Ybar}\left(\log D_{\Ybar},\wbar\right) = 
\bigoplus_{p = 0}^{a} \left[ \frac{1}{\wbar}\wedge^{p}W\bigoplus d\log\wbar
  \wedge \left(\wedge^{p-1}W\right)\right]\bigotimes \wedge^{a-p}R.
}
\end{equation}
This completes the proof of part (a) of the lemma. Part (b) follows
immediately from the formula \eqref{eq:adapted.forms}, the description
of the exceptional divisor of $\varepsilon : \widehat{\Ybar} \to \Ybar$ as a
projectivized normal bundle, and the Euler sequence of this projective
bundle. 
\ \hfill $\Box$

\

\begin{rem} \label{rem:complexes}
The $\wbar$-adapted logarithmic forms are equipped with two natural
differentials of degree one: 
\begin{itemize}
\item the de Rham differential $d : \Omega_{\Ybar}^{a}\left(\log
  D_{\Ybar},\wbar\right) \to \Omega_{\Ybar}^{a+1}\left(\log
  D_{\Ybar},\wbar\right)$, and
\item the differential $d\wbar\wedge : \Omega_{\Ybar}^{a}\left(\log
  D_{\Ybar},\wbar\right) \to \Omega_{\Ybar}^{a+1}\left(\log
  D_{\Ybar},\wbar\right)$. 
\end{itemize}
Note that by definition the differential $d\wbar\wedge$ is
$\mathcal{O}_{\Ybar}$-linear, while the de Rham differential satisfies
the Leibnitz rule as usual. Note also that for any complex numbers
$c_{1}$ and $c_{2}$  the linear combination $c_{1}d + c_{2}d\wbar\wedge$ is
also a differential and so we get a family of complexes of
$\wbar$-adapted logarithmic forms
\begin{equation} \label{eq:2param}
\left(\Omega_{\Ybar}^{\bullet}\left(\log
  D_{\Ybar},\wbar\right), c_{1}d + c_{2}d\wbar\wedge\right),
\end{equation}
parametrized by $(c_{1},c_{2}) \in
\mathbb{C}^{2}$. 
\end{rem}

\

\noindent
The previous discussion connects directly to the $L_{\infty}$-algebra
$R\Gamma(\Ybar,\dtg^{\bullet})$ since in the Calabi-Yau case
we can use the holomorphic volume form to convert $\wbar$-adapted
logarithmic forms to poly-vector fields. Suppose
$\left((\Ybar,\wbar),D_{\Ybar},\vol_{\Ybar}\right)$ is a compactified
Landau-Ginzburg model. The contraction with the meromorphic volume form
gives a map of $\mathcal{O}_{\Ybar}$-modules 
\begin{equation} \label{eq:contraction}
\xymatrix@R-2pc{
\con_{\vol_{\Ybar}}(\bullet): \hspace{-1pc} &  \wedge^{a} T_{\Ybar} \ar[r] &
\jmath_{Y*}\Omega^{n-a}_{Y} \\
& \xi \ar@{|->}[r] & \con_{\vol_{\Ybar}}(\xi).
}
\end{equation}
The preimage of $\Omega_{\Ybar}^{n-a}\left(\log
D_{\Ybar},\wbar\right)$ under the map \eqref{eq:contraction} will be a
coherent subsheaf in $\wedge^{a} T_{\Ybar}$. Furthermore when
$\left((\Ybar,\wbar),D_{\Ybar},\vol_{\Ybar}\right)$ satisfies the
tameness condition \eqref{eq:tame}, the explicit description of the
local frames of $\Omega_{\Ybar}^{n-a}\left(\log
D_{\Ybar},\wbar\right)$ above gives that 
\[
\left(\wedge^{a}T_{\Ybar}\right)\left(-\log D_{\Ybar},\wbar\right) :=
\left(\con_{\vol_{\Ybar}}(\bullet)\right)^{-1}\left(\Omega_{\Ybar}^{n-a}\left(\log
D_{\Ybar},\wbar\right)\right)
\] 
is a locally free subsheaf of maximal rank in $\wedge^{a} T_{\Ybar}$, and
that $\con_{\vol_{\Ybar}}$ induces an isomorphism between
$\left(\wedge^{a}T_{\Ybar}\right)\left(-\log D_{\Ybar},\wbar\right)$ and
$\Omega_{\Ybar}^{n-a}\left(\log D_{\Ybar},\wbar\right)$.

\

\noindent
With this notation we now have the following:

\begin{lemma} \label{lemma:bracket}
Let $\left((\Ybar,\wbar),D_{\Ybar},\vol_{\Ybar}\right)$ be a
compactified Landau-Ginzburg model satisfying the tameness assumption
\eqref{eq:tame}. Then the subsheaf
\[
\left(\wedge^{\bullet} T_{\Ybar}\right)\left(- \log
D_{\Ybar},\wbar\right) \subset \wedge^{\bullet} T_{\Ybar}
\]
is closed under the Nijenhuis-Schouten bracket on $\wedge^{\bullet} T_{\Ybar}$.
\end{lemma}
{\bfseries Proof.} Recall that the Nijenhuis-Schouten bracket on
polyvector fields is a degree $(-1)$ bracket that extends the Lie
bracket, acts as a graded derivation for the wedge product,  and is
given on decomposable polyvector fields by  
\begin{equation} \label{eq:bracket}
\left|\left|
\begin{split}
\left[\text{{\sf g}},\xi_{1}\wedge  \cdots \wedge \xi_{a}\right]  & =
\con_{d\text{{\sf g}}} \left(\xi_{1}\wedge \cdots \wedge
\xi_{a}\right), \\[+1pc]
\left[\xi_{1}\wedge \cdots  \wedge \xi_{a},\eta_{1}\wedge\cdots \wedge
  \eta_{b}\right]
& = \sum_{i,j} (-1)^{i+j}
\left[\xi_{i},\eta_{j}\right]\wedge \left( \xi_{1}\wedge
\cdots\wedge\widehat{\xi_{i}}\wedge\cdots  \wedge
\xi_{a}\right) \\[-1pc] 
& \hspace{1.5in} \wedge \left(\eta_{1}\wedge\cdots \wedge\widehat{\eta_{j}}
\wedge\cdots \wedge
  \eta_{b}\right),
\end{split}
\right.\right.
\end{equation}
for all $\text{{\sf g}} \in \mathcal{O}_{\Ybar}$, and all $\xi_{i},
\eta_{j} \in T_{\Ybar}$. 

The statement of the lemma is local on $\Ybar$ and is obvious away
from the points of $\ver{D}$. Indeed, away from $\ver{D}$ we have that
$\Omega^{n-a}_{\Ybar}\left(\log D_{\Ybar},\wbar\right)$ is isomorphic
to $\Omega^{n-a}_{\Ybar}\left(\log D_{\Ybar}\right)$. Since
$\vol_{\Ybar}$ has first order poles along the components of
$D_{\Ybar}$ this implies that on $Y - \ver{D}$ we have an isomorphism
$\left(\wedge^{a} T_{\Ybar}\right)\left(-\log D_{\Ybar},\wbar\right)
\cong \wedge^{a}T_{\Ybar,D_{\Ybar}}$. Since the subsheaf
$T_{\Ybar,D_{\Ybar}} \subset T_{\Ybar}$ is preserved by the Lie
bracket we get that away from $\ver{D}$ the subsheaf $\left(\wedge^{a}
T_{\Ybar}\right)\left(-\log D_{\Ybar},\wbar\right)$ is preserved by
the Nijenhuis-Scouten bracket.

Suppose next $p \in \ver{D}$. As before we choose local coordinates
$z_{1}, \cdots, z_{n}$ centered at $p$  so that near $p$ we have:
\begin{itemize}
\item $\ver{D} : z_{1}\cdots z_{k} = 0$,\quad  $\hor{D} :
  z_{k+1}\cdots z_{k+l} = 0$; 
\item ${\displaystyle \wbar(z_{1}, \ldots, z_{n}) =
  \frac{1}{z_{1}z_{2}\cdots z_{k}}}$;
\item ${\displaystyle \vol_{\Ybar} = \frac{dz_{1}\wedge \cdots
    \wedge dz_{n}}{z_{1}\cdots z_{k+l}}}$. 
\end{itemize}

Using this formula for $\vol_{\Ybar}$ and the description
\eqref{eq:adapted.forms} of the sheaf of $\wbar$-adapted logarithmic
forms, it is straightforward to compute $\left( \wedge^{a}
T_{\Ybar}\right)\left( -\log D_{\Ybar},\wbar\right)$. Let
\[
\begin{split}
M & =
\op{Span}_{\mathcal{O}_{\Ybar}}\left(z_{1}\frac{\partial}{\partial
  z_{1}},\ldots, z_{k-1}\frac{\partial}{\partial
  z_{k-1}}\right) \subset T_{\Ybar,D_{\Ybar}},\\[+0.5pc]
N & = \op{Span}_{\mathcal{O}_{\Ybar}}\left(z_{k+1}\frac{\partial}{\partial
  z_{k+1}},\ldots, z_{k+l}\frac{\partial}{\partial
  z_{k+l}},\frac{\partial}{\partial
  z_{k+l+1}},\ldots, \frac{\partial}{\partial
  z_{n}}\right) \subset T_{\Ybar,D_{\Ybar}}.
\end{split}
\]
In terms of these sheaves we have 
\begin{equation}
\label{eq:adapted.polyvectors}
\boxed{
\begin{split}
\left(\wedge^{a}T_{\Ybar}\right)& \left(-\log D_{\Ybar},\wbar\right)
\\[+0.5pc]
& = 
\bigoplus_{p = 0}^{a} \left[ \left(\frac{1}{\wbar}\wedge^{p-1}M\wedge
  z_{k}\frac{\partial}{\partial z_{k}}\right) \bigoplus
  \left(\con_{d\log\wbar}\left(\wedge^{p}M\wedge
  z_{k}\frac{\partial}{\partial z_{k}}\right)\right)
\right]\bigotimes \wedge^{a-p}N.
\end{split}
}
\end{equation}

From the formulas  \eqref{eq:bracket} it is now
immediate that $\left( \wedge^{\bullet}
T_{\Ybar}\right)\left( -\log D_{\Ybar},\wbar\right)$ is preserved by
the Nijenhuis-Schouten bracket. \ \hfill $\Box$

\

\begin{ex} \label{ex:dim1} It is instructive to examine more carefully
  the simplest case of a one dimensional compactified Landau-Ginzburg
  model. Near a point $p$ of $\ver{D}$ we can choose a local coordinate
  $z$ on $\Ybar$ so that
  $\wbar(z) = z^{-1}$, $\vol_{\Ybar} = dz/z$. Then locally near $p$ we
  get
\[
\begin{split}
\Omega^{\bullet}_{\Ybar}\left(\log D_{\Ybar}\right) & =
\mathcal{O}_{\Ybar}\cdot 1 \oplus \mathcal{O}_{\Ybar}\cdot
\frac{dz}{z}; 
\\[+0.5pc]
\Omega^{\bullet}_{\Ybar}\left(\log D_{\Ybar},\wbar\right) & =
\mathcal{O}_{\Ybar}\cdot z\oplus \mathcal{O}_{\Ybar}\cdot
\frac{dz}{z}; 
\\[+0.5pc]
\wedge^{\bullet} T_{\Ybar,D_{\Ybar}} & = \mathcal{O}_{\Ybar}\cdot 1
\oplus \mathcal{O}_{\Ybar}\cdot z \frac{\partial}{\partial z}; 
\\[+0.5pc]
\left(\wedge^{\bullet} T_{\Ybar}\right)\left(-\log
D_{\Ybar},\wbar\right) & = \mathcal{O}_{\Ybar}\cdot 1\oplus
\mathcal{O}_{\Ybar}\cdot z^{2}\frac{\partial }{\partial z}.
\end{split}
\]
\end{ex}

\

\noindent
Using Lemma~\ref{lemma:bracket} we can organize the $\wbar$-adapted
polyvector fields into a sheaf of dg Lie algebras. For any $1-n \leq b
\leq 1$ set 
\[
\pvg^{b} := \left( \wedge^{-b+1}T_{\Ybar}
\right)\left(-\log D_{\Ybar},\wbar\right).
\]
The sheaves  $\pvg^{b}$ 
fit together with the Nijenhuis-Schouten bracket $\left[ \bullet, \bullet
  \right]$  and the differential $\left[ \wbar, \bullet\right] =
\con_{d\wbar}$ into a sheaf of dg Lie algebras
\[
\left(\pvg^{\bullet}, [\wbar,\bullet]\right):= 
\xymatrix@R-2.5pc@M+0.5pc{
\left[\pvg^{1-n} \ar[r]^-{[\wbar,\bullet]}\right. & 
\pvg^{2-n} \ar[r]^-{[\wbar,\bullet]} &
\cdots \ar[r]^-{[\wbar,\bullet]}& 
\pvg^{0} \ar[r]^-{[\wbar,\bullet]} & 
\left. \pvg^{1} \right] \\
(1-n) & (2-n) & \cdots & 0 & 1 
}
\] 
This sheaf of dg Lie algebras is directly related to our
unobstructedness problem. Indeed, note that any stupid truncation of
$\left(\pvg^{\bullet}, [\wbar,\bullet]\right)$
will be a subsheaf of dg Lie algebras. In particular we have a
subsheaf of dg Lie algebras
\[
\left[ 
\xymatrix@1{\pvg^{0}
    \ar[r]^-{[\wbar,\bullet]} & \pvg^{1}} 
\right] = \sigma_{\geq 0}\left(\pvg^{\bullet},
     [\wbar,\bullet]\right)  \ 
 \hookrightarrow  \ \left(\pvg^{\bullet},
     [\wbar,\bullet]\right). 
\]
On the other hand this subsheaf maps naturally to the sheaf of dg Lie
algebras \linebreak $\dtg^{\bullet} = \left[\xymatrix@1{\dtg^{0}
    \ar[r]^-{d\wbar} & \dtg^{1}}\right]$ that controls our
deformation problem. So we get a diagram
\begin{equation} \label{eq:defo.roof}
\xymatrix{
& \sigma_{\geq 0} \left(\pvg^{\bullet}, 
     [\wbar,\bullet]\right)
\ar[dl] \ar[dr] & \\
\dtg^{\bullet}
 & &
  \left(\pvg^{\bullet}, 
     [\wbar,\bullet]\right) 
}
\end{equation}
of sheaves of dg Lie algebras. In fact \eqref{eq:defo.roof}
is a roof diagram. More precisely we have the following

\begin{prop} \label{prop:qis} The natural map of sheaves of dg Lie
  algebras 
\begin{equation} \label{eq:Gtog}
\left[\xymatrix@1{\pvg^{0}
    \ar[r]^-{[\wbar,\bullet]} & \pvg^{1}}\right]
\to
\left[\xymatrix@1{\dtg^{0}
    \ar[r]^-{d\wbar} & \dtg^{1}}\right] 
\end{equation}
is an $L_{\infty}$ quasi-isomorphism. 
\end{prop}
{\bfseries Proof.} The question is local on $\Ybar$. From the
definition of $\pvg^{\bullet}$ it is clear that the map
\eqref{eq:Gtog} is actually an isomorphism away from $\ver{D}$. Thus
it only remains to check the statement locally near a point 
$p \in \ver{D}$. 

Choose local coordinates $z_{1},\ldots,z_{n}$ as in the proof
of Lemma~\ref{lemma:bracket}. In terms of these coordinates we can
describe our dg Lie algebras explicitly. For the sheaves of 
$\wbar$-adapted  poly vector fields we have
\[
\left|\left| \quad
\begin{split}
\pvg^{0}  & = \left\{ \ \sum_{i=1}^{k} \sg_{i}
\frac{z_{i}\partial}{\partial z_{i}} +  \sum_{j=k+1}^{k+l} \sg'_{j}
\frac{z_{j}\partial}{\partial z_{j}} + \sum_{s=k+l}^{n} \sg''_{s}
\frac{\partial}{\partial z_{s}} \ \left| \ \quad
\text{
\begin{minipage}[c]{1.45in}
${\displaystyle \sg_{i}, \sg'_{j}, \sg''_{s} \in \mathcal{O}_{\Ybar}}$, and 
${\displaystyle
\sum_{i=1}^{k} \sg_{i} \in z_{1}\cdots z_{k}\mathcal{O}_{\Ybar}}$
\end{minipage}
}
 \right.\right\}, \\[+1pc]
\pvg^{1} & = \mathcal{O}_{\Ybar}.
\end{split}
\right.\right.
\]
The differential  $[\wbar,\bullet] =
\con_{d\wbar} : \pvg^{0} \to \pvg^{1}$ is given explicitly by
the formula
\begin{equation} \label{eq:differential}
\xymatrix@R-1pc{
\pvg^{0} \ar[r]^-{[\wbar,\bullet]} & \pvg^{1} \\
(\sg,\sg',\sg'') \ar[r] & {\displaystyle \frac{\sg_{1} + \sg_{2} + \cdots +
  \sg_{k}}{\sg_{1}\cdot \sg_{2} \cdot \cdots \cdot \sg_{k}}}.
}
\end{equation}
We have an analogous local description of the deformation theory dg algebra:
\[
\left|\left| \quad
\begin{split}
\dtg^{0}  & = \left\{\left. \ \sum_{i=1}^{k} \sg_{i}
\frac{z_{i}\partial}{\partial z_{i}} +  \sum_{j=k+1}^{k+l} \sg'_{j}
\frac{z_{j}\partial}{\partial z_{j}} + \sum_{s=k+l}^{n} \sg''_{s}
\frac{\partial}{\partial z_{s}} \ \right| \ \quad
\text{
\begin{minipage}[c]{1in}
${\displaystyle \sg_{i}, \sg'_{j}, \sg''_{s} \in \mathcal{O}_{\Ybar}}$
\end{minipage}
}\right\}, \\[+1pc]
\dtg^{1} & = \frac{1}{z_{1}\cdot \cdots \cdot
  z_{k}}\mathcal{O}_{\Ybar} = \wbar^{*}T_{\mathbb{P}^{1},\infty} \cong
\wbar^{*}\mathcal{O}_{\mathbb{P}^{1}}(1).
\end{split}
\right.\right.
\]
The differential $\dtg^{0} \to \dtg^{1}$ is again given by the formula
\eqref{eq:differential} and the map of dg Lie algebras $\left[\pvg^{0} \to
\pvg^{1}\right] \to \left[ \dtg^{0} \to \dtg^{1}\right]$ is given by the
natural inclusions $\pvg^{0} \subset \dtg^{0}$, $\pvg^{1} \subset
\dtg^{1}$. Thus we get a short exact sequence of complexes
\[
\xymatrix{
0 \ar[r] & \pvg^{0} \ar[r] \ar[d]_-{[\wbar,\bullet]} & \dtg^{0}
\ar[r] \ar[d]_-{d\wbar} &
\mathcal{O}_{\Ybar}/\left(z_{1}\cdots z_{k}\mathcal{O}_{\Ybar}\right)
\ar[r] \ar[d]_-{\frac{1}{z_{1}\cdots z_{k}}}& 0 \\
0 \ar[r] & \pvg^{1} \ar[r] & \dtg^{1}  \ar[r] &
\left(\frac{1}{z_{1}\cdots
  z_{k}}\mathcal{O}_{\Ybar}\right)/\mathcal{O}_{\Ybar} 
\ar[r] & 0
}
\]
Since the last vertical map is clearly an isomorphism, this implies
that $\left[\pvg^{0} \to
\pvg^{1}\right] \to \left[ \dtg^{0} \to \dtg^{1}\right]$ is a
quasi-isomorphism. 
\ \hfill $\Box$

\

\noindent
Proposition~\ref{prop:qis} and the roof diagram \ref{eq:defo.roof}
suggest that the unobstructedness statement in Theorem~\ref{theo:ha}
is related to the unobstructedness of the $L_{\infty}$ algebra
$R\Gamma(\Ybar,\left(\pvg^{\bullet},[\wbar,\bullet]\right)$. In fact
the standard formality yoga for $L_{\infty}$ algebras allows us to
deduce Theorem~\ref{theo:ha} from a stronger double degenration
statement for the cohomology of a two parameter family of $L_{\infty}$
algebras. Specifically let $\ddiv_{\vol_{\Ybar}} =
\con_{\vol_{\Ybar}}^{-1}\circ d \circ \con_{\vol_{\Ybar}} : \pvg^{a}
\to \pvg^{a+1}$ denote the divergence operator associated with
$\vol_{\Ybar}$. Note that by definition the differentials
$[\wbar,\bullet]$ and $\ddiv_{\vol_{\Ybar}}$ anticommute and so for
any pair of complex numbers $(c_{1},c_{2})$ we will get a well defined
complex
$R\Gamma\left(\Ybar,\left(\pvg^{\bullet},c_{1}\ddiv_{\vol_{\Ybar}} +
c_{2} [\wbar,\bullet]\right)\right)$. With this notation we now have
the following

\

\begin{prop} \label{prop:2param.algebra} 
Suppose that 
\begin{equation} \label{eq:2param.flat}
\text{ For all $a$ } \dim_{\mathbb{C}}
\mathbb{H}^{a}\left(\Ybar,\left(\pvg^{\bullet},c_{1}\ddiv_{\vol_{\Ybar}} +
c_{2} [\wbar,\bullet]\right)\right) \text{ is independent of $(c_{1},c_{2})
\in \mathbb{C}^{2}$.}
\end{equation}
 Then $R\Gamma(\Ybar,\dtg^{\bullet})$ is homotopy
abelian.
\end{prop}
{\bfseries Proof.} By Proposition~\ref{prop:qis} we deduce that the
$L_{\infty}$ algebra $R\Gamma\left(\Ybar,\dtg^{\bullet}\right)$ is
homotopy abelian (i.e. unobstructed) if and only if
$R\Gamma\left(\Ybar,\left[\xymatrix@1{\pvg^{0}
    \ar[r]^-{[\wbar,\bullet]} & \pvg^{1}}\right]\right)$ is homotopy
abelian.  Now in view of \cite[Proposition~4.11(ii)]{kkp} this
reduces\footnote{In \cite{kkp} Proposition~4.11 is formulated and
  proven for d$(\mathbb{Z}/2)$ graded algebras. However the statement
  of the proposition and its proof transfer verbatim to the
  d$(\mathbb{Z})$graded case, and we use this d$(\mathbb{Z})$graded
  version here.}  the unobstructedness statement in
Theorem~\ref{theo:ha} to showing that
\begin{itemize}
\item[(1)] The $L_{\infty}$ algebra $R\Gamma\left(\Ybar,\left(
  \pvg^{\bullet},[\wbar,\bullet]\right)\right)$ is homotopy
  abelian;
\item[(2)] The induced map 
\[
R\Gamma\left(\Ybar,\sigma_{\geq 0}\left(
  \pvg^{\bullet},[\wbar,\bullet]\right)\right) \to 
R\Gamma\left(\Ybar,\left(
  \pvg^{\bullet},[\wbar,\bullet]\right)\right)
\]
is  injective on cohomology. 
\end{itemize}

First note that the stupid filtration $\sigma_{\geq \bullet}
(\pvg^{\bullet},[\wbar,\bullet))$ gives rise to a spectral sequence
  which abuts to the spaces $\mathbb{H}^{a}\left(\Ybar,\left(
  \pvg^{\bullet},[\wbar,\bullet]\right)\right)$. By assumption
  $\dim_{\mathbb{C}} \mathbb{H}^{a}\left(\Ybar,\left(
  \pvg^{\bullet},c[\wbar,\bullet]\right)\right)$ is independent of $c
  \in \mathbb{C}$ and thus this spectral sequence will degenerate at
  $E_{1}$. This implies that $\mathbb{H}^{a}\left(\Ybar,\sigma_{\geq
    k}\left( \pvg^{\bullet},[\wbar,\bullet]\right)\right) \to
  \mathbb{H}^{a}\left(\Ybar,\left(
  \pvg^{\bullet},[\wbar,\bullet]\right)\right)$ is injective for all
  $k$ and in particular property (2) holds.

To prove property (1) we consider the flat family of $L_{\infty}$
algebras over $\mathbb{C}\left[\left[\hbar\right]\right]$ given by 
$\mathfrak{k} := R\Gamma(\Ybar,
(\pvg^{\bullet}\left[\left[\hbar\right]\right],
\hbar\cdot \ddiv_{\vol_{\Ybar}} + [ \wbar, \bullet])$. 
According to \cite[Proposition~4.11(i)]{kkp} it suffices to check
that $\mathfrak{k}$ satisfies:
\begin{itemize}
\item[(A)]
  $\mathfrak{k}\otimes_{\mathbb{C}\left[\left[\hbar\right]\right]}
  \mathbb{C}\left(\left(\hbar\right)\right)$ is homotopy abelian over
  $\mathbb{C}\left(\left(\hbar\right)\right)$, and  
\item[(B)] $H^{\bullet}\left(\mathfrak{k},d_{\mathfrak{k}}\right)$ is
  a flat $\mathbb{C}\left[\left[\hbar\right]\right]$-module.
\end{itemize}
Condition (B) follows immediately from the flatness assumption
\eqref{eq:2param.flat}. To check condition (A) we will use an
observation from \cite{bk98}: the map
$\mathfrak{k}\otimes_{\mathbb{C}\left[\left[\hbar\right]\right]}
\mathbb{C}\left(\left(\hbar\right)\right) \to
\mathfrak{k}\otimes_{\mathbb{C}\left[\left[\hbar\right]\right]}
\mathbb{C}\left(\left(\hbar\right)\right)$, given by $\gamma \mapsto
\exp(\gamma/\hbar)$, is a quasi-isomorphism between
$\mathfrak{k}\otimes_{\mathbb{C}\left[\left[\hbar\right]\right]}
\mathbb{C}\left(\left(\hbar\right)\right)$ and an abelian dg algebra
over $\mathbb{C}\left(\left(\hbar\right)\right)$.  The proposition is
proven.  \ \hfill $\Box$

\

\noindent
Proposition~\ref{prop:2param.algebra} finishes the proofs of
Theorem~\ref{theo:ha} and Theorem~\ref{theo:mainTT} modulo the flatness
assumption \eqref{eq:2param.flat}. Converting back to $\wbar$-adapted
logarithmic forms  via 
$\con_{\vol_{\Ybar}}$ the assumption  \eqref{eq:2param.flat} is equivalent to 
the statement that the dimension of the hypercohomology
$\mathbb{H}^{\bullet}(\Ybar,(\Omega^{\bullet}_{\Ybar}(\log
D_{\Ybar},\wbar),c_{1}d + c_{2}d\wbar\wedge))$ is independent of
$(c_{1},c_{2}) \in \mathbb{C}^{2}$. We investigate this
Hodge theoretic statement in the next section. 

\subsection{The double degeneration
  property} \label{ssec:double.degeneration} 

In this section we complete the proof of the unobstructedness
Theorem~\ref{theo:mainTT} by establishing the double
degeneration property for the complex \eqref{eq:2param} of
$\wbar$-adapted logarithmic forms associated with a compactified tame
Landau-Ginzburg model which is not necessarily of log Calabi-Yau
type. An alternative proof and a generalization of this statement
can be found in the recent work of Esnault-Sabbah-Yu \cite{esy}. For
the convenience of the reader we give our original argument here.

\

\begin{theo} \label{theo:double.degeneration} Let
  $((\Ybar,\wbar),D_{\Ybar})$ be geometric datum where
\begin{itemize}
\item[(a)] $\Ybar$ is a smooth projective variety, and $\wbar : \Ybar
  \to \mathbb{P}^{1}$ is a flat projective morphism.
\item[(b)] $D_{\Ybar} = \left(\cup_{i} \hor{D}_{i}  \right)\cup \left(
  \cup_{j} \ver{D}_{i} \right)$ is a reduced normal crossing divisor,
  such that 
\begin{itemize}
\item $\ver{D} = \cup_{j} \ver{D}_{j}$ is the pole divisor of $\wbar$,
  i.e. $\ver{D} = \wbar^{-1}(\infty)$ is the scheme-theoretic fiber of
  $\wbar$ at $\infty \in \mathbb{P}^{1}$.
\item $\op{crit}(\wbar) \cap \hor{D}_{\Ybar} = \varnothing$.
\end{itemize}
\end{itemize}
Then the following flatness property holds
\begin{equation}
\label{eq:2param.forms.flat}
\boxed{
\text{
\begin{minipage}[c]{4in}
For all $a \geq 0$  
$\dim_{\mathbb{C}} \mathbb{H}^{a}\left(\Ybar,\left(
\Omega^{\bullet}_{\Ybar}(\log 
D_{\Ybar},\wbar),c_{1}d + c_{2}d\wbar\wedge \right)\right)$  is
  independent of $(c_{1},c_{2}) \in \mathbb{C}^{2}$.
\end{minipage}
}}
\end{equation}
\end{theo}
{\bfseries Proof.} We will obtain the proof by checking the constancy
of dimension of cohomology along various lines in $\mathbb{C}^{2}$.
First we have the following

\

\begin{lemma} \label{lemma:Hodge.degeneration}
For every $((\Ybar,\wbar),D_{\Ybar})$ satisfying the hypothesis if the
theorem and every $a \geq 0$ we have 
\begin{equation} \label{eq:Hodge.degeneration}
\begin{split}
\dim_{\mathbb{C}} \mathbb{H}^{a}(\Ybar,
(\Omega_{\Ybar}^{\bullet}(\log D_{\Ybar},\wbar),d)) & = \dim_{\mathbb{C}}
  \mathbb{H}^{a}(\Ybar,
(\Omega_{\Ybar}^{\bullet}(\log D_{\Ybar},\wbar),0)) \\[+0.5pc] 
& = \sum_{i+j = a}
\dim_{\mathbb{C}} H^{i}(\Ybar,\Omega_{\Ybar}^{j}(\log
D_{\Ybar},\wbar)).
\end{split}
\end{equation}
In particular, the spectral sequence corresponding to the stupid
filtration on $(\Omega_{\Ybar}^{\bullet}(\log D_{\Ybar},\wbar),d)$
degennerates at $E_{1}$. 
\end{lemma}
{\bfseries Proof.} We will use the method of Deligne-Illusie
\cite{deligne-illusie,esnault-viehweg,illusie}. Here we only sketch the
necessary modifications that make the method applicable to
$\wbar$-adapted logarithmic forms. More details can be found 
in the Esnault-Sabbah-Yu writeup in \cite[Appendix~D]{esy}. 

By the standard spreading-out argument of
\cite{deligne-illusie,esnault-viehweg,illusie} it suffices to check
the $E_{1}$ degeneration of the spectral sequence
\begin{equation} \label{eq:deg.over.k}
 H^{i}(\Ybar,\Omega_{\Ybar/\bk}^{j}(\log
D_{\Ybar},\wbar)) \Rightarrow \mathbb{H}^{i+j}(\Ybar,
(\Omega_{\Ybar/\bk}^{\bullet}(\log D_{\Ybar},\wbar),d))
\end{equation}
in the case when the geometric datum $((\Ybar,\wbar),D_{\Ybar})$
satisfying the hypotheses of the lemma is defined over a perfect field
$\bk$ of characteristic $p > \dim X$ and admits a smooth lift to
characteristic $0$ (or at least to the second Witt vectors $W_{2}(\bk)$
of $\bk$).  

Write $((\Ybar',\wbar'),D'_{\Ybar'})$ for the Frobenius twist of the datum
$((\Ybar,\wbar),D_{\Ybar})$. In other words
$((\Ybar',\wbar'),D'_{\Ybar'})$ is the the base change of
$((\Ybar,\wbar),D_{\Ybar})$ by the absolute Frobenius map \linebreak 
$\varphi : \op{Spec}
\bk \to \op{Spec} \bk$. Let $\Phi : ((\Ybar',\wbar'),D'_{\Ybar'}) \to
((\Ybar,\wbar),D_{\Ybar})$ be the base change map and let $\Fr :
((\Ybar,\wbar),D_{\Ybar}) \to ((\Ybar',\wbar'),D'_{\Ybar'})$ denote the
induced relative Frobenius morphism over $\bk$. 

The base change property for algebraic differential forms combined
with the fact that $\Fr$ is a homeomorphism, and with the local
description \eqref{eq:adapted.forms} of $\wbar$-adapted forms implies
that we have canonical isomorphisms $\Phi^{*}
H^{i}\left(\Ybar,\Omega^{j}_{\Ybar/\bk}(\log D_{\Ybar},\wbar)\right)
\cong H^{i}\left(\Ybar',\Omega^{j}_{\Ybar'/\bk}(\log
D_{\Ybar'},\wbar')\right)$ and
$H^{a}\left(\Ybar',\Fr_{*}\left(\Omega^{\bullet}_{\Ybar/\bk}(\log
D_{\Ybar},\wbar),d\right)\right) =
H^{a}\left(\Ybar,\left(\Omega^{\bullet}_{\Ybar/\bk}(\log 
D_{\Ybar},\wbar),d\right)\right)$. This gives equality of dimensions
of these matching  cohomology groups and so  the $E_{1}$ degeneration
of \eqref{eq:deg.over.k} will follow immediately (see
e.g. \cite[Section 4.8]{illusie}) if we can show that the complex 
$\Fr_{*}\left(\Omega^{\bullet}_{\Ybar/\bk}(\log
D_{\Ybar},\wbar),d \right)$ is formal as an object in the derived
category of quasi-coherent $\mathcal{O}_{\Ybar'}$-modules.

To that end, recall \cite{cartier,katz-nilpotent} that the (inverse)
Cartier map defined by $\car(\Phi^{*}z) = z^{p}$ and $\car(d\Phi^{*}z)
= \left[z^{p-1}dz\right]$ on a local function $z$ on $\Ybar$, extends
uniquely by multiplicativity and gives rise to an isomorphism
\[
\xymatrix@1@M+0.5pc{
\car : \hspace{-2.5pc} &  
\bigoplus_{a \geq 0} \Omega^{a}_{\Ybar'/\bk}(\log D'_{\Ybar'})
\ar[rr]^-{\cong} & & \bigoplus_{a \geq 0}
\mycal{H}^{a}\left(\Fr_{*}\left(\Omega^{\bullet}_{\Ybar/\bk}(\log
D_{\Ybar}),d\right)\right). 
}
\]
of sheaves of
super commutative algebras over $\mathcal{O}_{\Ybar'}$. 

Using the
explicit local description \eqref{eq:adapted.forms} of the
$\wbar$-adapted logarithmic forms one checks immediately 
that $\car$ also restricts to an
isomorphism 
\begin{equation} \label{eq:adapted.cartier}
\xymatrix@1@M+0.5pc{
\car : \hspace{-2.5pc} & 
\bigoplus_{a \geq 0} \Omega^{a}_{\Ybar'/\bk}(\log D'_{\Ybar'},\wbar)
\ar[rr]^-{\cong} & &  \bigoplus_{a \geq 0}
\mycal{H}^{a}\left(\Fr_{*}\left(\Omega^{\bullet}_{\Ybar/\bk}(\log
D_{\Ybar},\wbar),d\right)\right).
}
\end{equation}
of sheaves of
super commutative algebras over $\mathcal{O}_{\Ybar'}$.

In view of the isomorphism \eqref{eq:adapted.forms} the formality of
$\Fr_{*}\left(\Omega^{\bullet}_{\Ybar/\bk}(\log D_{\Ybar},\wbar),d
\right)$ as an object in $D(\Ybar')$ is equivalent to the existence of
a morphism in $D(\Ybar')$:
\[
\xymatrix@1@M+0.5pc{
\for : \hspace{-2.5pc} & \bigoplus_{a\geq 0}
\Omega^{a}_{\Ybar'/\bk}(\log D'_{\Ybar'},\wbar')[-a] \ar[rr] & & 
\Fr_{*}\left(\Omega^{\bullet}_{\Ybar/\bk}(\log D_{\Ybar},\wbar),d
\right)
}
\]
which induces $\car$ on cohomology sheaves. Following
\cite{deligne-illusie} the construction of $\for$ can be carried out
in three stages. Fix a lift $\left( (\mathfrak{Z},\mathfrak{f}),
D_{\mathfrak{Z}} \right)$ over $W_{2}(\bk)$. We abuse notation and
again write $\varphi : \op{Spec} W_{2}(\bk) \to \op{Spec} W_{2}(\bk)$
for the absolute Frobenius. Similarly we will write $\left(
(\mathfrak{Z}',\mathfrak{f}'), D'_{\mathfrak{Z}'} \right)$ for the
pullback of $\left( (\mathfrak{Z},\mathfrak{f}),
D_{\mathfrak{Z}} \right)$ via $\varphi$ and will write 
$\Phi : \mathfrak{Z}' \to \mathfrak{Z}$ for the base change map.

As a first step suppose that the relative Frobenius $\Fr : \Ybar \to
\Ybar'$ admits a global lifting to a morphism $\lFr : \mathfrak{Z} \to
\mathfrak{Z}'$ of $W_{2}(\bk)$-schemes which furthermore satisfies
$\lFr^{*}(\mathfrak{f}') = \mathfrak{f}^{p}$ and $\lFr^{*}
\mathcal{O}_{\mathfrak{Z}'}\left(D'_{\mathfrak{Z}'}\right) =
\mathcal{O}_{\mathfrak{Z}}\left(p\cdot D_{\mathfrak{Z}}\right)$.  With
such a lifting we associate a formality morphism $\for_{\lFr}$ as
follows:
\begin{itemize}
\item For $a = 0$ we set $\for^{0}_{\lFr} = \Fr^{*} :
\mathcal{O}_{\Ybar'} \to  \Fr_{*}\mathcal{O}_{\Ybar}$;
\item For $a = 1$ we set $\for^{1}_{\lFr} = \left( (1/p)\cdot \lFr^{*}
  \ \op{mod} \ p \right) \ :
  \Omega^{1}_{\Ybar'/\bk}\left(\log D'_{\Ybar'},\wbar'\right) \to
  \Fr_{*}\Omega^{1}_{\Ybar/\bk}\left(\log D_{\Ybar}\right)$;
\item For $a > 1$ we define $\for^{a}_{\lFr}$ to be the composition of 
$\wedge^{a} \for^{1}_{\lFr}$ with the product map $\wedge^{a} \Fr_{*}
  \Omega^{1}_{\Ybar/\bk}\left(\log D_{\Ybar}\right) \to
  \Omega^{a}_{\Ybar/\bk}\left(\log D_{\Ybar}\right)$. 
\end{itemize}
The key observation now is that for all $a$ the map $\for^{a}_{\lFr}$
sends $\Omega^{a}_{\Ybar'/\bk}\left(\log D'_{\Ybar'},\wbar'\right)$
to \linebreak $\Fr_{*}\Omega^{a}_{\Ybar/\bk}\left(\log
D_{\Ybar},\wbar\right)$. Once this is checked, the fact that
$\for_{\lFr}$ is a quasi-isomorphism inducing $\car$ on all
cohomology sheaves follows tautologically. To show that 
\[
\for^{a}_{\lFr}\left(\Omega^{a}_{\Ybar'/\bk}\left(\log
D'_{\Ybar'},\wbar'\right)\right) \subset \Fr_{*}\Omega^{a}_{\Ybar/\bk}\left(\log
D_{\Ybar},\wbar\right)
\] 
we argue locally on $\Ybar$. By \eqref{eq:adapted.forms} we know that 
$\Omega^{a}_{\Ybar/\bk}\left(\log
D_{\Ybar},\wbar\right)$ is a locally free sheaf which near 
$\ver{D}_{\Ybar}$ is equal to the sum of
$d\log\wbar \wedge \Omega^{a-1}_{\Ybar/\bk}(\log D_{\Ybar})$ and
$(1/\wbar)\Omega^{a}_{\Ybar/\bk}(\log D_{\Ybar})$ inside
$\Omega^{a}_{\Ybar/\bk}(\log D_{\Ybar})$.

Choosing an appropriate Zariski local etale map to an affine space we
obtain local coordinates $\bz = (z_{1}, \ldots z_{n})$ as in the proof
of Lemma~\ref{lemma:locally.free}. In particular we have that the
divisor $D_{\Ybar}$ is given by the union of $z_i=0$ for $i = 1,
\ldots, k+l$, and $\wbar = 1/(z_1\cdot\ldots\cdot z_k)$. In these
coordinates the lifted Frobenius $\lFr$ has the form
$\lFr^{*}(\Phi^{*}z_i) = z_i^p+p \cdot z_{i}^{p} v_i(\bz)$ for $i = 1,
\ldots, k+l$ and $\lFr^{*}(\Phi^{*}z_i) = z_i^p + p \cdot v_i(\bz)$
for $i = k+l+1, \ldots, n$. Furthermore we have $v_1(\bz) + \cdots +
v_{k}(\bz) =0$.

From these formulas we now see that for the forms in
$\Omega^a_{\Ybar/\bk}(\log D_{\Ybar},\wbar)$ of type $d\log\wbar
\wedge \Omega^{a-1}_{\Ybar/\bk}(\log D_{\Ybar})$ the pullback via
$(1/p^{a})\lFr^{*}$ modulo $p$ does not depend on the choice of
$v_{1}, \ldots, v_{k}$, and is therefore again a log form multiplied
by $d\log\wbar$. For the forms of second type, i.e. forms in
$(1/\wbar)\Omega^a_{\Ybar/\bk}(\log D_{\Ybar},\wbar) \subset
\Omega^a_{\Ybar/\bk}(\log D_{\Ybar},\wbar)$ we note that these forms
belong to the $\mathcal{O}$-module generated by products over all $i$
of either $z_{i}$ or $dz_i$.  Hence the $(1/p^{a})\lFr^{*}$ pullback
of such form modulo $p$ will belong to the $\mathcal{O}$ module
generated by products of either $z_{i}^{p}$ or $z_{i}^{p-1} dz_{i}$
and is therefore again of second type.

In the second step one notes that locally in the Zariski topology we
can always choose etale maps to an affine space and then use local
coordinates as above to construct a lift of the relative Frobenius
over $W_{2}(\bk)$. Thus we have to analyze the relation between the
formality isomorphisms associated to different local liftings of the
relative Frobenius. Following \cite{deligne-illusie} we want to show
that for any 
two liftings $\lFr_{1} : \mathfrak{Z}_{1} \to \mathfrak{Z}'$ and
$\lFr_{2} : \mathfrak{Z}_{2} \to \mathfrak{Z}'$  of $\Fr$, we
can find a canonical map of sheaves 
\[
h(\lFr_{1},\lFr_{2}) : \Omega^{1}_{\Ybar'/\bk}\left(\log
D'_{\Ybar'},\wbar'\right)  \to \Fr_{*}\mathcal{O}_{\Ybar}
\]
so that $\for^{1}_{\lFr_{1}} - \for^{1}_{\lFr_{2}} = d
h(\lFr_{1},\lFr_{2})$. Furthermore for  a third lifting $\lFr_{3} :
\mathfrak{Z}_{3} \to \mathfrak{Z}'$  these maps should satisfy the
cocycle condition $h(\lFr_{1},\lFr_{2}) + h(\lFr_{2},\lFr_{3}) =
h(\lFr_{1},\lFr_{3})$.

To construct the maps $h$ we repeat verbatim the reasoning in
\cite{deligne-illusie,esnault-viehweg}. To show that the corresponding
$d h(\lFr_{1},\lFr_{2})$ belongs again to
$\Fr_{*}\Omega^{1}_{\Ybar'/\bk}\left(\log D'_{\Ybar'},\wbar'\right)$
one notes that by construction $h(\lFr_{1},\lFr_{2})$ is given by
substitutions with vector fields of the form $\sum_{i = 1}^{k+l}
u_{i}(\bz)\cdot z_{i}^{p}\cdot (\partial/\partial z_{i}) + \sum_{i =
  k+l+1}^{n} u_{i}(\bz)\cdot (\partial/\partial z_{i})$ satisfying
$\sum_{i=1}^{k} u_{i} = 0$.  Such substitution vanish for forms
divisible by $d\log\wbar$. For forms $\alpha= (1/\wbar)\cdot \beta$
with $\beta$ being log form, the substitution of such vector field in
$\beta$ is again a log form, hence its pullback is a log form. Also
note that the pullback of $1/\mathfrak{f}'$ is $1/\mathfrak{f}^{p}$
which is divisible by $1/\mathfrak{f}$. Thus we again get a form of
the second type.  The same argument should work for a triple of
lifts. The key point here is that the forms $\Omega^a_{\Ybar/\bk}(\log
D_{\Ybar},\wbar)$ are closed under contractions with vector fields
preserving $\wbar$.

From this point on the argument proceeds exactly as in
\cite{deligne-illusie}. First we cover $\Ybar$ by Zariski open sets
$U_{i}$ on
which we can choose Frobenius lifts $\lFr_{i} : \mathfrak{U}_{i} \to
\mathfrak{U}'_{i}$ as above. Then on overlaps we use the maps
$h(\lFr_{i},\lFr_{j})$
on overlaps to glue the formality morphisms $\for^{1}_{\lFr_{i}}$
into a morphism in the derived category 
\[
\for^{1}_{\mathfrak{Z}} : \Omega^{1}_{\Ybar'/\bk}(\log
D'_{\Ybar'},\wbar')[-1] \to \Fr_{*}\Omega^{\bullet}_{\Ybar/\bk}(\log
D_{\Ybar},\wbar)
\]
which induces $\car$ on $\mathcal{H}^{1}$.

In the last step we use the condition $n = \dim \Ybar < p$ and
multiplicative structure on the de Rham 
complex to define a map 
\[
\for^{a}_{\mathfrak{Z}} : \Omega^{a}_{\Ybar'/\bk}(\log
D'_{\Ybar'},\wbar')[-1] \to \Fr_{*}\Omega^{\bullet}_{\Ybar/\bk}(\log
D_{\Ybar},\wbar)
\]
by composing 
$\left(\for^{1}_{\mathfrak{Z}}\right)^{\otimes a}$ with the
anti-symmetrization map $\Omega^{a}_{\Ybar'/\bk}(\log  D'_{\Ybar'})
\to \left(\Omega^{1}_{\Ybar'/\bk}(\log  D'_{\Ybar'})\right)^{\otimes
  a}$ given by $\alpha_{1}\otimes \cdots \otimes \alpha_{a} \mapsto 
\frac{1}{a!} \sum_{\sigma \in S_{a}} \op{sgn}(\sigma)
\alpha_{\sigma(1)}\otimes \cdots \otimes \alpha_{\sigma(a)}$. 

This completes the proof of the lemma. \ \hfill $\Box$

\

\begin{rem} \label{rem:saito} Morihiko Saito recently found
  \cite{saito.appendixE} a different analytic proof of this lemma. Saito's
  argument uses Hodge theory with degenerating coefficients and takes
  place entirely in characteristic zero.
\end{rem}

\

\noindent 
Lemma~\ref{lemma:Hodge.degeneration} implies that the dimension of the
hypercohohomology of the complex \linebreak
$\left(\Omega^{\bullet}_{\Ybar}(\log D_{\Ybar},\wbar), c_{1}d +
  c_{2}d\wbar\wedge \ \right)$ is constant on the line $\{ c_{2} = 0
\} \subset \mathbb{C}^{2}$. Next we will show that this
hypercohomology is also constant along the line $\{ c_{1} = c_{2} \}
\subset \mathbb{C}^{2}$. First we have the following topological
statement (see also \cite[Appendix~C]{esy} where a more general
statement allowing multiplicities is proven)

\

\begin{lemma} \label{lemma:cohomology.of.the.pair} Consider $\bw : Y
  \to \mathbb{C}$. Write $Y_{-\infty}$ for the fiber $\bw^{-1}(z)$
  over $z \in \mathbb{C}$ with $\op{Re} z \ll 0$. Then for every $a
  \geq 0$ we have 
\[
\dim_{\mathbb{C}} \mathbb{H}^{a}\left(\Ybar,\left(\Omega^{\bullet}_{\Ybar}(\log
  D_{\Ybar},\wbar), d\right)\right) = \dim_{\mathbb{C}}
H^{a}(Y,Y_{-\infty};\mathbb{C}). 
\]
\end{lemma} {\bfseries Proof.} Before we address the statement of the
lemma, it is instructive to look at the analogous statement in the
classical Hodge theory of smooth open varieties.  Let $\jmath_{Y} : Y
\hookrightarrow \Ybar$ be the natural inclusion viewed as a continuous
map in the analytic topology. The pushforward
$R\jmath_{Y*}\mathbb{C}_{Y}$ is a constructible complex of sheaves of
$\mathbb{C}$-vector spaces on $\Ybar$. Recall that the log de Rham
complex $\left(\Omega^{\bullet}_{\Ybar}(\log D_{\Ybar}), d\right)$ is
naturally quasi-isomorphic to this constructible complex. Indeed, a
direct local calculation \cite{griffiths-periods,deligne-h2} shows
that the natural map of complexes $\left(\Omega^{\bullet}_{\Ybar}(\log
  D_{\Ybar}), d\right) \to \jmath_{Y*}\mathcal{A}^{\bullet}_{Y}$ is a
quasi-isomorphism. Composing this map with the augmentation
quasi-isomorphism $\jmath_{Y*}\mathcal{A}^{\bullet}_{Y} \to
R\jmath_{Y*}\mathbb{C}_{Y}$ gives an identification of
$\left(\Omega^{\bullet}_{\Ybar}(\log D_{\Ybar}), d\right)$ and
$R\jmath_{Y*}\mathbb{C}_{Y}$ in
$D^{b}\left(\mathbb{C}_{\Ybar}\right)$. Since
$R\jmath_{Y*}\mathbb{C}_{Y}$ computes the Betti cohomology of the open
variety $Y$ this yields the classical statement that
$\mathbb{H}^{a}\left(\Ybar, \left(\Omega^{\bullet}_{\Ybar}(\log
    D_{\Ybar}), d\right)\right) \cong H^{a}(Y,\mathbb{C})$.

The idea is to modify this reasoning to take into account relative
cohomology and $\wbar$-adapted forms. To that end consider the real
oriented blow-up $\varepsilon : \widehat{\Ybar} \to \Ybar$ of $\Ybar$
along the reduced normal crossing divisor $D_{\Ybar}$, and the real
oriented blow-up $\pi : \widehat{\mathbb{P}}^{1} \to \mathbb{P}^{1}$
of $\mathbb{P}^{1}$ at $\infty \in \mathbb{P}^{1}$. The morphism
$\wbar : \Ybar \to \mathbb{P}^{1}$ lifts naturally to a real
semi-algebraic map $\hatwbar : \widehat{\Ybar} \to
\widehat{\mathbb{P}}^{1}$. The spaces $\widehat{\mathbb{P}}^{1}$ and
$\widehat{\Ybar}$ are manifolds with boundary, and $\partial
\widehat{\mathbb{P}}^{1} = \pi^{-1}(\infty) \cong S^{1}$ and $\partial
\widehat{\Ybar} = \varepsilon^{-1}\left(D_{\Ybar}\right) \supset
\varepsilon^{-1}\left(\ver{D}_{\Ybar}\right) =
\hatwbar^{-1}\left(\partial \widehat{\mathbb{P}}^{1}\right)$.

The boundary circle $\partial \widehat{\mathbb{P}}^{1} =
\pi^{-1}(\infty) \cong S^{1}$ is the circle of radial directions at
$\infty \in \mathbb{P}^{1}$. If as before $z$ denotes the affine
coordinate on $\mathbb{A}^{1} = \mathbb{P}^{1} - \{\infty\}$, then
this circle is parametrized by $\op{arg}(1/z)$. Choose a point
$\theta_{0} \in \partial \widehat{\mathbb{P}}^{1}$ for which
$\op{Re}(z) \geq 0$ and let $\widehat{\Ybar}_{\theta_{0}} =
\hatwbar^{-1}(\theta_{0})$. 

Consider the complex
$\left(\mathcal{A}_{\widehat{\Ybar},\widehat{\Ybar}_{\theta_{0}}}^{\bullet}\left(\log
D_{\Ybar}\right),d\right)$ of $C^{\infty}$ logarithmic forms on
$\widehat{\Ybar}$ that vanish along
$\widehat{\Ybar}_{\theta_{0}}$. Let
$\left(\mathcal{A}_{\widehat{\Ybar}_{\theta_{0}}}^{\bullet}\left(\log
D_{\Ybar}\right),d\right)$ denote the cone (= quotient complex) of the
natural map from
$\left(\mathcal{A}_{\widehat{\Ybar},\widehat{\Ybar}_{\theta_{0}}}^{\bullet}\left(\log
D_{\Ybar}\right),d\right)$ to
$\left(\mathcal{A}_{\widehat{\Ybar}}^{\bullet}\left(\log
D_{\Ybar}\right),d\right)$. Now from the explicit local description
(see Lemma~\ref{lemma:locally.free}) of $\left(
\Omega^{\bullet}_{\Ybar}\left(\log D_{\Ybar},\wbar\right), d\right)$
one checks immediately that near the boundary $\partial
\widehat{\Ybar}$ the quotient complex $\left(
\Omega^{\bullet}_{\Ybar}\left(\log
D\right)\left/\Omega^{\bullet}_{\Ybar}\left(\log
D_{\Ybar},\wbar\right)\right., d\right)$ maps to
$\varepsilon_{*}\left(\mathcal{A}_{\widehat{\Ybar}_{\theta_{0}}}^{\bullet}\left(\log
D_{\Ybar}\right),d\right)$ and that the map is a
quasi-isomorphism. This implies that the natural map from
$\Omega^{\bullet}_{\Ybar}\left(\log D_{\Ybar},\wbar\right)$
$\varepsilon_{*}
\left(\mathcal{A}_{\widehat{\Ybar},\widehat{\Ybar}_{\theta_{0}}}^{\bullet}\left(\log
D_{\Ybar}\right),d\right)$ is a quasi-isomorphism and  thus gives the
equality of dimensions claimed in the lemma.

 A detailed writeup of this argument and an explicit check of the
fact that the map of quotient complexes is a quasi-isomorphism can be
found in \cite[Appendix~C,Step~2]{esy}. Instead  of repeating this
calculation here we will give an alternative proof of the lemma
which is  be of independent interest. 

To simplify the discussion let us first assume that $\hor{D}_{\Ybar}$
is empty. Consider the de Rham cohomology of the pair $\left(Y,
\bw^{-1}(\rho)\right)$ where $\rho$ is real and $\rho \ll 0$. Using de
Rham's theorem and the Gauss-Manin parallel transport along the ray
$\rho \in \mathbb{R}_{< 0}$ we can identify the cohomology
$H^{a}\left(Y,Y_{-\infty};\mathbb{C}\right)$ with the limit of
$H^{a}_{DR}\left(Y,\bw^{-1}(\rho);\mathbb{C}\right)$ as $\rho \to
-\infty$. 

Therefore the statement of the lemma reduces to understanding the
limit \linebreak ${\displaystyle \lim_{\rho \to -\infty}\displaylimits
  H^{a}_{DR}\left(Y,\bw^{-1}(\rho);\mathbb{C}\right)}$ in terms of the
complex of $\wbar$-adapted logarithmic forms. The relative de Rham
cohomology $H^{a}_{DR}\left(Y,\bw^{-1}(\rho);\mathbb{C}\right)$ is
computed by the complex \linebreak $\left( \Omega^{\bullet}_{\Ybar}\left(\log
D_{\Ybar},\op{rel} \, \wbar^{-1}(\rho)\right),d\right)$ of holomorphic
forms on $\Ybar$ that have logarithmic poles along $D_{\Ybar} =
\wbar^{-1}(\infty)$ and vanish along the divisor $\wbar^{-1}(\rho) =
\bw^{-1}(\rho)$. We now have the following

\begin{claim} \label{claim:limit.complex} 
\begin{itemize}
\item[(a)] 
As
$\rho \to -\infty$ the complex $\left( \Omega^{\bullet}_{\Ybar}\left(\log
D_{\Ybar},\op{rel} \, \wbar^{-1}(\rho)\right),d\right)$ has a well defined
limit, namely the complex $\left( \Omega^{\bullet}_{\Ybar}\left(\log
D_{\Ybar},\wbar \right),d\right)$ of $\wbar$-adapted
logarithmic forms on $\Ybar$.
\item[(b)] The Gauss-Manin parallel transport along the ray $\rho \in
  \mathbb{R}_{< 0}$ is well defined at the limit $\rho \to -\infty$
  and identifies $H^{a}_{DR}\left(Y,\bw^{-1}(\rho);\mathbb{C}\right)$
  with
  $\mathbb{H}^{a}\left(\Ybar,\left(\Omega^{\bullet}_{\Ybar}\left(\log
  D_{\Ybar}, \wbar \right), d \right)\right)$.
\end{itemize}
\end{claim}
{\bfseries Proof.} The statement is local on $\Ybar$ and is obvious at
points of the open set $Y = \Ybar - D_{\Ybar}$. Suppose next $p \in
\ver{D}_{\Ybar} = D_{\Ybar} \subset \Ybar$. As in the proof of
Lemma~\ref{lemma:locally.free}  
we can choose a local coordinate system $z_{1}, \ldots, z_{n}$
centered at $p$ so that $D_{\Ybar}$ is given by the equation
$\prod_{i = 1}^{k} z_{i} = 0$ and $\wbar = \prod_{i = 1}^{k}
z_{i}^{-1}$. Now as in the proof of  Lemma~\ref{lemma:locally.free} 
we write $W \subset \Omega^{1}_{\Ybar}\left( \log D_{\Ybar}\right)$
for the sub $\mathcal{O}_{\Ybar}$-module spanned by $d\log z_{1}, \ldots,
d\log z_{k-1}$, and $R \subset \Omega^{1}_{\Ybar}\left( \log
D_{\Ybar}\right)$ for sub $\mathcal{O}_{\Ybar}$-module spanned by $dz_{k+1},
\ldots, dz_{n}$.

In these terms we have
\[
\Omega^{a}_{\Ybar}\left(\log D_{\Ybar}\right) = \bigoplus_{p = 0}^{a}
\left[ \wedge^{p}W \oplus d\log z_{k}\wedge \left(\wedge^{p-1}
W\right)\right]\bigotimes \wedge^{a-p} R. 
\]
Write $\beps = 1/\rho$, and let $Y_{\beps} =
\wbar^{-1}(\beps)$. Then for $\beps$ close to $0$  we can use
$z_{1}$, \ldots, $z_{k-1}$, $z_{k+1}$, \ldots, $z_{n}$ as coordinates
along the divisor $Y_{\beps}$. In particular the sheaf of
holomorphic forms $\Omega^{1}_{Y_{\beps}}$ is the
$\mathcal{O}_{Y_{\beps}}$-span of the forms $dz_{1}$, \ldots,
$dz_{k-1}$, $dz_{k+1}$, \ldots, $dz_{n}$, and so 
\[
\Omega^{a}_{Y_{\beps}} = \bigoplus_{p = 0}^{a}
\wedge^{p}W_{|Y_{\beps}}\bigotimes \wedge^{a-p} R_{|Y_{\beps}}. 
\]
From here we get 
\[
\begin{split}
\Omega^{a}_{\Ybar}\left(\log D_{\Ybar}, \op{rel} \,
Y_{\beps}\right) & = \ker\left[ \Omega^{a}_{\Ybar}\left(\log
  D_{\Ybar}\right) \to \imath_{Y_{\beps}*}
  \Omega^{a}_{Y_{\beps}}\right]  \\[+0.5pc]
& = \bigoplus_{p = 0}^{a}
\left[ (z_{1}\cdots z_{k} - \beps)\wedge^{p}W \oplus d\log
z_{k}\wedge \wedge^{p-1} 
W\right]\bigotimes \wedge^{a-p} R \\[+0.5pc] 
& = \bigoplus_{p = 0}^{a}
\left[ (z_{1}\cdots z_{k} - \beps)\wedge^{p}W + d\log
\wbar \wedge \left(\wedge^{p-1} 
W\right)\right]\bigotimes \wedge^{a-p} R,
\end{split}
\]
and so when $\beps \to 0$ this sheaf specializes to the sheaf of
$\wbar$-adapted logarithmic forms (see \eqref{eq:adapted.forms})
\[
\Omega^{a}_{\Ybar}\left(\log D_{\Ybar},\wbar\right) =
\bigoplus_{p = 0}^{a}
\left[ z_{1}\cdots z_{k}\cdot\wedge^{p}W \oplus d\log
\wbar \wedge \left(\wedge^{p-1} 
W\right)\right]\bigotimes \wedge^{a-p} R.
\]
This shows that as $\beps \to 0$ the complex ${\displaystyle
  \left(\Omega^{\bullet}_{\Ybar}\left(\log D_{\Ybar}, \op{rel} \,
      Y_{\beps}\right), d\right)}$ will converge to the complex
${\displaystyle \left(\Omega^{\bullet}_{\Ybar}\left(\log
      D_{\Ybar},\wbar\right), d\right)}$. In other words we have a
family of complexes on $\Ybar$ parametrized by a small complex number
$\beps$, where ${\displaystyle
  \left(\Omega^{\bullet}_{\Ybar}\left(\log D_{\Ybar}, \op{rel} \,
      Y_{\beps}\right), d\right)}$ is the complex corresponding to
$\beps \neq 0$, while at $\beps = 0$ we have the complex
${\displaystyle \left(\Omega^{\bullet}_{\Ybar}\left(\log
      D_{\Ybar},\wbar\right), d\right)}$. 

More invariantly, let
$\bDelta \subset \mathbb{P}^{1}$ be a small disk centered at $\infty$
with coordinate $\beps$. Let $\bZ := \Ybar\times \bDelta$, and let $p
: \bZ \to \bDelta$ be the natural projection. The proper family $p :
\bZ \to \bDelta$ is equipped with two relative divisors
\[
\begin{split}
D_{\bZ} & := D_{\Ybar}\times \bDelta, \\[+0.5pc]
\bGamma & := (p\times\wbar)^{-1}\left(\op{{\sf graph}}\left(\bDelta
    \hookrightarrow \mathbb{P}^{1}\right)\right).
\end{split}
\]
By construction $\bGamma$ is smooth, $D_{\bZ}$ has strict normal
crossings, both $\bGamma$ and $D_{\bZ}$ are flat over $\bDelta$, and
the union of  $\bGamma \cup D_{\bZ}$ also has strict normal crossings.
Write $D_{\bGamma}$ for the normal crossing divisor in $\bGamma$ given
by $D_{\bGamma} = D_{\bZ}\cap \Gamma$.

Consider now the sheaves of relative meromorphic forms, i.e. forms
along the fibers of $p$, having logarithmic poles along $D_{\bZ}$ and
vanishing along $\Gamma$:
\[
\Omega^{a}_{\bZ/\bDelta}\left(\log D_{\bZ}, \op{rel} \, \bGamma
\right) := \ker \left[ \Omega^{a}_{\bZ/\bDelta}\left(\log
    D_{\bZ}\right) \to \imath_{\bGamma*} \Omega^{a}_{\bGamma/\bDelta}\left(\log
    D_{\bGamma}\right)\right].
\]
By definition these are locally free sheaves of (certain) relative logarithmic
forms along the fibers  of $p : \bZ \to \bDelta$,  and the graded
subsheaf $\Omega^{\bullet}_{\bZ/\bDelta}\left(\log D_{\bZ}, \op{rel} \, \bGamma
\right) \subset \Omega^{\bullet}_{\bZ/\bDelta}\left(\log
  D_{\bZ}\right)$ is clearly preserved by the relative de Rham
differential. The calculation in local coordinates above shows that
the complex 
\[
\bE^{\bullet}_{\bZ/\bDelta} : = 
\left( \Omega^{\bullet}_{\bZ/\bDelta}\left(\log D_{\bZ}, \op{rel} \, \bGamma
\right), \, d\right) 
\]
interpolates between relative logarithmic forms vanishing on $Y_{\beps}$ and
$\wbar$-adapted relative  logarithmic forms. In other words we have 
\[
\begin{split}
\left(\bE^{\bullet}_{\bZ/\bDelta}\right)_{|\Ybar \times \{\beps \neq 0\}} &
= \left( \Omega^{\bullet}_{\Ybar}\left(\log D_{\Ybar}, \op{rel} \, Y_{\beps}
\right), \, d\right); \\[+0.5pc]
\left(\bE^{\bullet}_{\bZ/\bDelta}\right)_{|\Ybar \times \{0\}} &
= \left( \Omega^{\bullet}_{\Ybar}\left(\log D_{\Ybar}, \wbar \right),
  \, d\right). 
\end{split}
\]
This proves the first part of the claim.

The statement about the Gauss-Manin parallel transport follows form
the homological description \cite{katz-oda} of the Gauss-Manin
connection. To spell this out one needs to describe the local system
of relative cohomology via differential forms. To fix notation, we
will use a superscript $(\bullet)^{\times}$ to indicate the removal of
the fiber over $\beps = 0$ in the various geometric and
sheaf-theoretic objects we are dealing with. Thus we will write
$\bDelta^{\times} = \bDelta - \{ 0 \}$, $\bZ^{\times} = \bZ -
p^{-1}(0)$, $D_{\bZ^{\times}} = \Ybar \times \bDelta^{\times}$,
$\bGamma^{\times} = \bGamma\cap \bZ^{\times}$, and
\[
\bE^{\bullet}_{\bZ^{\times}/\bDelta^{\times}} = \left(
  \Omega^{\bullet}_{\bZ^{\times}/\bDelta^{\times}}\left( \log
    D_{\bZ^{\times}}, \op{rel} \, \bGamma^{\times}\right), d\right).
\] 
Also we set $\bY = Y\times \bDelta = \bZ - D_{\bZ}$ and $\bY^{\times}
= Y\times \bDelta^{\times}$. 
 
Let $\mycal{E}^{a}_{B}$ denote the local system of $\mathbb{C}$-vector
spaces on $\bDelta^{\times}$ whose fiber over $\beps \in
\bDelta^{\times}$ is the relative Betti cohomology
$H^{a}(Y,Y_{\beps};\mathbb{C})$.  The underlying coherent sheaf
$\mycal{E}^{a}_{B}\otimes_{\mathbb{C}} \mathcal{O}_{\bDelta^{\times}}$
can be identified with the sheaf $\mathcal{H}^{a}_{DR}
\left(\bY^{\times}/\bDelta^{\times},\bGamma^{\times}/\bDelta^{\times};
  \mathbb{C}\right)$ of relative de Rham cohomology and is
thus computed as the hyperderived image
\[
\mycal{E}^{a}_{B}\otimes_{\mathbb{C}} \mathcal{O}_{\bDelta^{\times}}
\; \cong \;
\mathcal{H}^{a}_{DR}
\left(\bY^{\times}/\bDelta^{\times},\bGamma^{\times}/\bDelta^{\times};
  \mathbb{C}\right) \; = \; \mathbb{R}^{a}p_{*}
\bE^{\bullet}_{\bZ^{\times}/\bDelta^{\times}}.
\]
Note that the hyper-derived image $\mathbb{R}^{a}p^{\times}_{*}
\bE^{\bullet}_{\bZ^{\times}/\bDelta^{\times}}$ is naturally an
$\mathcal{O}_{\bDelta^{\times}}$-module since the de Rham differential
on relative forms is linear over
$p^{-1}\mathcal{O}_{\bDelta^{\times}}$. In these terms the Gauss-Manin connection
is given by a $\mathbb{C}$-linear map of sheaves
\[
\nabla^{\op{GM}} :
\mathbb{R}^{a}p_{*}\bE^{\bullet}_{\bZ^{\times}/\bDelta^{\times}}
  \longrightarrow 
\mathbb{R}^{a}p_{*}\bE^{\bullet}_{\bZ^{\times}/\bDelta^{\times}}
\otimes_{\mathcal{O}_{\bDelta^{\times}}}
  \Omega^{1}_{\bDelta^{\times}},
\]
satisfying the Leibnitz rule. The analysis of \cite{katz-oda} applies
verbatim to this setting and identifies $\nabla^{\op{GM}}$ with the
connecting homomorphism in the long exact sequence of hyperderived
direct images associated with the short exact sequence of complexes
\begin{equation} \label{eq:GMses}
\xymatrix@1{
0 \ar[r] &
*\txt{
\begin{minipage}[c]{2in}
\vspace{1pc}
$\bE^{\bullet}_{\bZ^{\times}/\bDelta^{\times}}[-1]
\underset{p^{-1}\mathcal{O}_{\bDelta^{\times}}}{\otimes}
p^{-1}\Omega^{1}_{\bDelta^{\times}}
$
\end{minipage}
} 
\ar[r] & 
\left(\Omega^{\bullet}_{\bZ^{\times}}\left( \log D_{\bZ^{\times}},
    \op{rel} \, \bGamma^{\times} \right), d\right)
\ar[r] & 
*\txt{
\begin{minipage}[c]{0.63in}
\vspace{0.4pc}
$\bE^{\bullet}_{\bZ^{\times}/\bDelta^{\times}}$
\end{minipage}
}
\ar[r] &
0.}
\end{equation}
In order to check that the parallel transport with respect to
$\nabla^{\op{GM}}$ has a well defined limit when $\beps \to 0$ it
suffices to show that the complex
$\left(\Omega^{\bullet}_{\bZ^{\times}}\left( \log D_{\bZ^{\times}},
    \op{rel} \, \bGamma^{\times} \right), d\right)$ extends to a well
defined subcomplex $\bE^{\bullet}_{\bZ}$ in
$\left(\Omega^{\bullet}_{\bZ}\left( \log D_{\bZ}\right), d\right)$, so
that  $\bE^{\bullet}_{\bZ}$ is defined on all of $\bZ$ and
fits in a short exact sequence of complexes  
\begin{equation} \label{eq:GMses.extended}
\xymatrix@1{
0 \ar[r] &
*\txt{
\begin{minipage}[c]{1.7in}
\vspace{1pc}
$\bE^{\bullet}_{\bZ/\bDelta}[-1]
\underset{p^{-1}\mathcal{O}_{\bDelta}}{\otimes}
p^{-1}\Omega^{1}_{\bDelta}
$
\end{minipage}
} 
\ar[r] & 
\bE^{\bullet}_{\bZ}
\ar[r] & 
*\txt{
\begin{minipage}[c]{0.48in}
\vspace{0.4pc}
$\bE^{\bullet}_{\bZ/\bDelta}$
\end{minipage}
}
\ar[r] &
0}
\end{equation}
extending 
\eqref{eq:GMses} to  all of $\bZ$.

The naive guess of taking $\bE^{\bullet}_{\bZ}$ to be
$\left(\Omega^{\bullet}_{\bZ}\left( \log D_{\bZ}, \op{rel} \,
\bGamma\right), d\right)$ will not work since the natural maps
$\Omega^{a}_{\bZ}\left( \log D_{\bZ}, \op{rel} \, \bGamma\right) \to
\Omega^{a}_{\bZ/\bDelta}\left( \log D_{\bZ}, \op{rel} \,
\bGamma\right)$ are not surjective. Because of this we will have to
work with the logarithmic de Rham complexes directly.  Consider the
short exact sequence
\begin{equation} \label{eq:logarithmic.ses}
\xymatrix@1{
0 \ar[r] & \Omega^{\bullet}_{\bZ/\bDelta}\left(\log D_{\bZ}\right)[-1]
\otimes_{\mathcal{O}_{\bZ}} p^{*}\Omega^{1}_{\bDelta} \ar[r] & 
\Omega^{\bullet}_{\bZ}\left(\log D_{\bZ}\right) \ar[r] &
\Omega^{\bullet}_{\bZ/\bDelta}\left(\log D_{\bZ}\right) \ar[r] &
0. 
}
\end{equation}
of logarithmic de Rham complexes\footnote{Here all terms are equipped
with the obvious absolute or relative de Rham differentials, so we
have omitted them from the notation. The differential in the first
term is defined via the identification
$\Omega^{\bullet}_{\bZ/\bDelta}\left(\log D_{\bZ}\right)[-1]
\otimes_{\mathcal{O}_{\bZ}} p^{*}\Omega^{1}_{\bDelta} \cong
\Omega^{\bullet}_{\bZ/\bDelta}\left(\log D_{\bZ}\right)[-1]
\otimes_{p^{-1}\mathcal{O}_{\bDelta}} p^{-1}\Omega^{1}_{\bDelta}$ and
the fact that de Rham differential on relative forms along the fibers
of $p : \bZ \to \bDelta$ is $p^{-1}\mathcal{O}_{\bDelta}$-linear.}  on
$\bZ$.  

View \eqref{eq:logarithmic.ses} as a morphism $\xi_{\bZ/\bDelta} :
\Omega^{\bullet}_{\bZ/\bDelta}\left(\log D_{\bZ}\right) \to
\Omega^{\bullet}_{\bZ/\bDelta}\left(\log D_{\bZ}\right) 
\otimes_{\mathcal{O}_{\bZ}} p^{*}\Omega^{1}_{\bDelta}$ in the derived
category of sheaves of $\mathbb{C}$-vector spaces on $\bZ$. Write
$\imath_{\bGamma} : \bGamma \hookrightarrow \bZ$ for the 
inclusion of the divisor $\bGamma$ in $\bZ$, and let $i$ and $q$
denote the maps in the short exact sequence of complexes defining
$\bE^{\bullet}_{\bZ/\bDelta}$:
\[
\xymatrix@1{
0 \ar[r] & \bE^{\bullet}_{\bZ/\bDelta} \ar[r]^-{i} &
\Omega^{\bullet}_{\bZ/\bDelta}\left(\log D_{\bZ}\right) \ar[r]^-{q} & 
\imath_{\bGamma*} \Omega^{\bullet}_{\bGamma/\bDelta}\left(\log
D_{\bGamma}\right) \ar[r] &
0.
}
\]
We claim that the composition 
\[
(q\otimes 1)\circ \xi_{\bZ/\bDelta}\circ i :
\bE^{\bullet}_{\bZ/\bDelta} \to \imath_{\bGamma*} \left[
\Omega^{\bullet}_{\bGamma/\bDelta}\left(\log 
D_{\bGamma}\right)\otimes_{\mathcal{O}_{\Gamma}} p^{*}\Omega^{1}_{\bDelta}\right]
\]
is the zero morphism in $D^{b}\left(\mathbb{C}_{\bZ}\right)$. This
follows immediately by noting that $\xi_{\bZ/\bDelta}$ fits in a commutative
diagram in $D^{b}\left(\mathbb{C}_{\bZ}\right)$:
\begin{equation} \label{eq:restriction.diagram}
\xymatrix@C+1pc{
\bE^{\bullet}_{\bZ/\bDelta} \ar[d]_-{i} 
& \bE^{\bullet}_{\bZ/\bDelta}\otimes_{\mathcal{O}_{\bZ}}
p^{*}\Omega^{1}_{\bDelta}  \ar[d]_-{i\otimes 1} \\
\Omega^{\bullet}_{\bZ/\bDelta}\left(\log D_{\bZ}\right) 
\ar[d]_-{q}
\ar[r]^-{\xi_{\bZ/\bDelta}}  &
\Omega^{\bullet}_{\bZ/\bDelta}\left(\log D_{\bZ}\right)\otimes_{\mathcal{O}_{\bZ}}
p^{*}\Omega^{1}_{\bDelta} \ar[d]_-{q\otimes 1} \\
\imath_{\bGamma*} \Omega^{\bullet}_{\bGamma/\bDelta}\left(\log
D_{\bGamma}\right) 
\ar[r]^-{\imath_{\bGamma*}\xi_{\bGamma/\bDelta}}  &
\imath_{\bGamma*}\left[ \Omega^{\bullet}_{\bGamma/\bDelta}\left(\log
D_{\bGamma}\right)\otimes_{\mathcal{O}_{\bGamma}}
p^{*}\Omega^{1}_{\bDelta}\right]
}
\end{equation}
in which the columns are parts of exact triangles and
\[
\xi_{\bGamma/\bDelta} : \Omega^{\bullet}_{\bGamma/\bDelta}\left(\log
D_{\bGamma}\right) \to \Omega^{\bullet}_{\bGamma/\bDelta}\left(\log
D_{\bGamma}\right)\otimes_{\mathcal{O}_{\bGamma}}
p^{*}\Omega^{1}_{\bDelta}
\] 
is the map in
$D^{b}\left(\mathbb{C}_{\bGamma}\right)$ corresponding to the short
exact sequence of complexes
\[
\xymatrix@1{
0 \ar[r] & \Omega^{\bullet}_{\bGamma/\bDelta}\left(\log D_{\bGamma}\right)[-1]
\otimes_{\mathcal{O}_{\bGamma}} p^{*}\Omega^{1}_{\bDelta} \ar[r] & 
\Omega^{\bullet}_{\bGamma}\left(\log D_{\bGamma}\right) \ar[r] &
\Omega^{\bullet}_{\bGamma/\bDelta}\left(\log D_{\bGamma}\right) \ar[r] &
0 
}
\]
of logarithmic forms on $\bGamma$. The vanishing of $(q\otimes 1)\circ
\xi_{\bZ/\bDelta}\circ i$ follows immediately now since from 
\eqref{eq:restriction.diagram} we see that $(q\otimes 1)\circ
\xi_{\bZ/\bDelta}\circ i =
\left(\imath_{\bGamma*}\xi_{\bGamma/\bDelta}\right)\circ q\circ i =
0$. This in turn implies that $\xi_{\bZ/\bDelta}$ comes from a
morphism $\op{cone}(q) \to \op{cone}(q\otimes 1)$ in
$D^{b}\left(\mathbb{C}_{\bZ}\right)$. In other words we can find a map
$\xi_{\bE} : \bE^{\bullet}_{\bZ/\bDelta} \to
\bE^{\bullet}_{\bZ/\bDelta}\otimes_{\mathcal{O}_{\bZ}}
  p^{*}\Omega^{1}_{\bDelta}$ so that $(i\otimes 1)\circ \xi_{\bE} =
  \xi_{\bZ/\bDelta}\circ i$. Since $i$ is an isomorphism over the open
  $\bZ^{\times} \subset \bZ$ it now follows that over $\bZ^{\times}$
  the map 
  $\xi_{\bE}$ coincides with the map $\xi_{\bE}^{\times} : 
  \bE^{\bullet}_{\bZ^{\times}/\bDelta^{\times}} \to 
\bE^{\bullet}_{\bZ^{\times}/\bDelta^{\times}}\otimes_{p^{-1}\mathcal{O}_{\bDelta^{\times}}}
p^{-1}\Omega^{1}_{\bDelta^{\times}}$ corresponding to the short exact
    sequence of complexes \eqref{eq:GMses}. 
Since $\nabla^{\op{GM}} = \mathbb{R}^{a}p_{*}(\xi_{\bE}^{\times})$ it
follows that $\nabla^{\op{GM}}$ extends to a holomorphic connection
$\mathbb{R}^{a}p_{*}(\xi_{\bE}) :
\mathbb{R}^{a}p_{*}\bE^{\bullet}_{\bZ/\bDelta} \to
\mathbb{R}^{a}p_{*}\bE^{\bullet}_{\bZ/\bDelta}\otimes_{\mathcal{O}_{\bDelta}}
\Omega^{1}_{\bDelta}$. This completes the proof of the claim. \ \hfill $\Box$

\

\noindent
The statement of Claim~\ref{claim:limit.complex} proves
Lemma~\ref{lemma:cohomology.of.the.pair} in the case when the divisor
$D_{\Ybar}$ does not have a horizontal part. In fact we can
incorporate the horizontal divisor into the proof of
Claim~\ref{claim:limit.complex} without any modification. The local
calculation for the limit, and the extension argument repeat verbatim,
only the notation becomes more cumbersome. We will not spell this out
here and leave it to the interested reader to fill in the details. \
\hfill $\Box$

\

\noindent
We can now complete the proof of the double degeneration property by
combining  Lemma~\ref{lemma:cohomology.of.the.pair} with the following
well known facts:

\begin{lemma} \label{lemma:barannikov.kontsevich} For every $a \geq 0$
  we have 
\[
\begin{split}
\dim_{\mathbb{C}} H^{a}\left(Y,Y_{-\infty}; \,
\mathbb{C}\right) & =
\dim_{\mathbb{C}}\mathbb{H}^{a}\left(Y_{\op{Zar}}, 
\left(
\Omega^{\bullet}_{Y}, d + d\bw\wedge \right)\right) \\
& = \dim_{\mathbb{C}}\mathbb{H}^{a}\left(Y_{\op{Zar}}, \left(
\Omega^{\bullet}_{Y}, d\bw\wedge \right)\right) 
\end{split}
\]
\end{lemma}
{\bfseries Proof.} The first equality is the usual identification of
de Rham nearby cycles with twisted de Rham cohomology via
the Fourier transform for regular holonomic D-modules on the affine
line. The second is the degeneration theorem for twisted de Rham
complexes proven in the work of Barannikov and Kontsevich
(unpublished), Sabbah \cite{sabbah-twisted}, or Ogus-Vologodsky
\cite{ogus.vologodsky}. \ \hfill $\Box$

\

\begin{lemma} \label{lemma:all.poles} For every $a \geq 0$
  we have 
\[
\dim_{\mathbb{C}}\mathbb{H}^{a}\left(Y_{\op{Zar}}, \left(
\Omega^{\bullet}_{Y}, d + d\bw\wedge \right)\right) 
= \dim_{\mathbb{C}}\mathbb{H}^{a}\left(\Ybar_{\op{Zar}}, \left(
\Omega^{\bullet}_{\Ybar}(\log D_{\Ybar}, \wbar), d + d\wbar\wedge 
\right)\right) 
\]
\end{lemma}
{\bfseries Proof.} This follows from the usual Grothendieck argument
\cite{grothendieck-deRham}. The local calculation comparing
logarithhmic forms with meromorphic forms transfers immediately to the
$\wbar$-adapted complex and combined with the local description of
adapted forms given in Lemma~\ref{lemma:locally.free} implies that the
natural inclusion of complexes 
\[
\left( \Omega^{\bullet}_{\Ybar}(\log
D_{\Ybar}, \wbar), d + d\wbar\wedge \right) \hookrightarrow \left(
\Omega^{\bullet}_{\Ybar}(* D_{\Ybar}), d + d\wbar\wedge \right) =
R\jmath_{Y*}\left( \Omega^{\bullet}_{Y}, d + d\bw\wedge \right)
\] 
is a
quasi-isomorphism. \ \hfill $\Box$

\

\noindent
Taken together Lemmas~\ref{lemma:cohomology.of.the.pair},
\ref{lemma:barannikov.kontsevich}, and 
\ref{lemma:all.poles} imply that the hypercohomology of the complex 
$\left(
\Omega^{\bullet}_{\Ybar}(\log D_{\Ybar}, \wbar), c_{1}\cdot d +
c_{2}\cdot  d\wbar\wedge 
\right)$ is constant on the line $c_{1} = c_{2}$. This completes the
proof of Theporem~\ref{theo:double.degeneration}. \ \hfill $\Box$

\section{Invariants of \nc \ 
Hodge structures of geometric
  origin} \label{s-variations}

In this section we use the deformation theory developed in
Section~\ref{sec:moduliLG} to elucidate the motivic and Hodge
theoretic data naturally present on the cohomology of a compactifiable
Landau-Ginzburg model. Using considerations from mirror symmetry we
propose various new refined invariants of \nc \ Hodge structures of
Landau-Ginzburg type, and discuss in particular the subtleties
involved in understaning  Hodge numbers and decorations.

\subsection{Hodge numbers of Landau-Ginzburg
  models} \label{ssec:hodge.numbers} 

Suppose $\left( (\Ybar,\wbar), D_{\Ybar}, \vol_{\Ybar}\right)$ is a
tame compactified Landau-Ginzburg model in the sense of
Definition~\ref{defi:compactLG} and
assumption~\eqref{eq:tame}. Conjecturally (see \cite{kkp}) the
cohomology
\linebreak 
$H^{\bullet}\left(Y,Y_{-\infty}; \, \mathbb{C}\right)$ of the associated
quasi-projective Landau-Ginzburg model $\bw : Y \to \mathbb{A}^{1}$
carries a $B$-model pure \nc \ Hodge structure. The de Rham data for
this \nc \ Hodge structure is described in \cite[Section~3.2]{kkp}, where
it is also argued that this data satisfies the \nc \ Hodge filtration
axiom. The much trickier opposedness axiom has only been verified for
models $((Y,\bw),\vol_{Y})$ which mirror a general symplectic toric
weak Fano manifold \cite{reichelt-sevenheck}.

A somewhat dissapointing feature of \nc \ Hodge structures in general
is that their complexity is not readily captured by simple numerical
invariants. The absence of easy to compute linear-algebraic invariants
in this setting is a reflection of the nature of the \nc \ de Rham
datum. The \nc \ Hodge filtration is encoded in a connection with
irregular singularities, and the Stokes structures characterizing this
connection can not be encoded in simple linear algebraic quantities.
The special nature of the Landau-Ginzburg context however allows one
to discern additional sophisticated linear-algebraic data compatible
with the \nc \ Hodge structure on $H^{\bullet}\left(Y,Y_{-\infty}; \,
  \mathbb{C}\right)$. Moreover, this data posseses computable
numerical invariants, such as weights, level, amplitude, and Hodge
numbers. Most naturally this additional data arises from the concept
of an irregular Hodge filtration that can be associated with a
Landau-Ginzburg potential. There are two variants of such irregular
Hodge filtrations - the version of Deligne and Sabbah
\cite{deligne-irregular,sabbah-irregular}, and the version of J.-D. Yu
\cite{yu-irregular,esy}. In \cite{esy} these two variants of the
irregular Hodge filtration are generalized, ultimately identified with
each other, and under the assumption \eqref{eq:tame} identified with
the Hodge filtration on the complex of $\wbar$-adapted logarithmic
forms. Here we will not discuss this identification but rather will
look more closely at the resulting Hodge numbers and will compare
those to other more classical definitions of Hodge numbers arising
from vanishing cohomology and mirror data.

Given a Calabi-Yau Landau-Ginzburg model $\bw : Y \to \mathbb{A}^{1}$
which admits a tame compactification $\left( (\Ybar,\wbar), D_{\Ybar},
  \vol_{\Ybar}\right)$, we will define geometrically three sets of
Hodge numbers $i^{p,q}(Y,\bw)$, $h^{p,q}(Y,\bw)$, and
$f^{p,q}(Y,\bw)$, each of which adds up to the rank of the
algebraic de Rham cohomology $H^{a}_{DR}\left( (Y_{\op{Zar}}, \bw); \,
  \mathbb{C} \right) = \mathbb{H}^{a}\left(Y_{\op{Zar}}, \left(
    \Omega^{\bullet}_{Y}, d + d\bw\wedge\right)\right)$ of the
Landau-Ginzburg model. Since by
Lemma~\ref{lemma:barannikov.kontsevich} we have 
$\dim_{\mathbb{C}} H^{a}_{DR}\left( Y_{\op{Zar}}, \bw; \, \mathbb{C}
\right) = \dim_{\mathbb{C}} H^{a}\left( Y, Y_{-\infty}; \, \mathbb{C}
\right)$ this implies
\[
\boxed{
\dim_{\mathbb{C}} H^{a}\left( Y, Y_{-\infty}; \, \mathbb{C}
  \right) =  \sum_{p+q = a} i^{p,q}(Y,\bw) = \sum_{p+q = a}
h^{p,q}(Y,\bw) = \sum_{p+q = a} f^{p,q}(Y,\bw).
}
\]
Each of these sets of Hodge numbers has a different origin. The
numbers $i^{p,q}(Y,\bw)$ come from ordinary mixed Hodge theory, the
numbers $h^{p,q}(Y,\bw)$ come from mirror considerations, and the
numbers $f^{p,q}(Y,\bw)$ come from the sheaf cohomology of the
$\wbar$-adapted logarithmic forms. The specific definitions are as
follows.

\

\smallskip

\noindent
\punkt \ {\bfseries The numbers $f^{p,q}(Y,\bw)$.} \ Let $\left(
(\Ybar,\wbar), D_{\Ybar}, \vol_{\Ybar}\right)$ be a tame
compactification of $\bw : Y \to \mathbb{A}^{1}$.

\begin{defi} \label{defi:hpq} The 
{\bfseries Landau-Ginzburg Hodge numbers $f^{p,q}(Y,\bw)$} are defined
by
\[
\boxed{
f^{p,q}(Y,\bw) = \dim_{\mathbb{C}} H^{p}(\Ybar,\Omega^{q}_{\Ybar}(\log
D_{\Ybar}, \wbar)).
}
\]
\end{defi}

\

\noindent
The fact that $\dim_{\mathbb{C}} H^{a}\left( Y, Y_{-\infty}; \, \mathbb{C}
  \right) =  \sum_{p+q = a} f^{p,q}(Y,\bw)$ follows from
  Theorem~\ref{theo:mainTT} and 
  Lemma~\ref{lemma:all.poles}. 

\

\smallskip

\noindent
\punkt \ {\bfseries The numbers $h^{p,q}(Y,\bw)$.} \ Before we explain
the definition we need to recall a basic construction from linear
algebra. Let $V$ be a finite dimensional complex vector space,  $N
: V \to V$ be a nilpotent linear operator satisfying $N^{m+1} = 0$ for
some $m \geq 0$. The
{\bfseries (monodromy) weight filtration of $N$ centered at $m$} is
the unique increasing filtration $W = W_{\bullet}(N,m)$ of $V$:
\[
0 \subset W_{0}(N,m) \subset W_{1}(N,m) \subset \cdots \subset
W_{2m-1}(N,m) \subset 
W_{2m}(N,m) = V
\]
 with the properties 
\begin{itemize}
\item $N (W_{i}(N,m)) \subset W_{i-2}(N,m)$;
\item the map ${\displaystyle N^{\ell} :  \op{gr}^{W}_{m+\ell} V \to
    \op{gr}^{W}_{m-\ell} V}$ 
is an isomorphism for all $\ell \geq 0$.
\end{itemize}
The existence and uniquencess of this filtration can be deduced from
the representation theory of $\mathfrak{s}\mathfrak{l}_{2}$-triples
and the Jacobson-Morozov theorem. A direct elementary proof can also
be found in \cite[Lemma~6.4]{schmid}.  Explicitly the monodromy weight
filtration $W_{\bullet}(N,m)$ is defined as follows.  Choose a Jordan
basis for the nilpotent endomorphism $N : V \to V$ and assign integer
weights to the basis vectors so that $N$ lowers weights by $2$, and so
that the weights of each Jordan block are arranged symmetrically about
$m$. Note that even though the Jordan canonical form is not canonical,
the monodromy weight filtration will be canonical since $W_{k}(N,m)$ is the
span of the basis vectors of weights less than or equal to $k$.

Now let $c_{0} \in \mathbb{A}^{1}$ be a regular value of $\bw$ near
infinity. Consider the monodromy transformation 
 $T : H^{\bullet}(Y,Y_{c_{0}};\, \mathbb{C}) \to
H^{\bullet}(Y,Y_{c_{0}};\, \mathbb{C})$ 
corresponding to moving the smooth fiber $Y_{c_{0}}$ once around
infinity. By assumption $(Y,\bw)$ admits a  tame compactification and
so, as explained in Remark~\ref{rem:tameness}{\bfseries (ii)} the
operator $T$ is unipotent. Set $N = \log \, T$. With this notation we
have the following

\begin{defi} \label{defi:hpq} The 
{\bfseries Landau-Ginzburg Hodge numbers $h^{p,q}(Y,\bw)$} are defined
by 
\[
\boxed{
h^{p,q}(Y,\bw) := \dim_{\mathbb{C}} \op{gr}_{p}^{W(N,p+q)}
H^{p+q}(Y,Y_{c_{0}}; \, \mathbb{C}). 
}
\]
\end{defi}

\smallskip

\

\noindent
The rationale behind this definition is the geometric mirror symmetry
prediction explained in  Remark~\ref{rem:tameness}{\bfseries
  (ii)}. Specifically, if $(Y,\bw)$ is part of a mirror pair 
\[
\left(X,\omega_{X},s_{X}\right) \bmid
\left((Y,\bw),\omega_{Y},\vol_{Y}\right)
\]
of Fano type with $\dim_{\mathbb{C}} X = \dim_{\mathbb{C}} Y = n$,
then the homological mirror equivalence 
\[
\sD^{b}(X) \cong \FS((Y,\bw),\omega_{Y},\vol_{Y}) 
\]
induces an isomorphism on period cyclic and on Hochschild homologies
of these categories. In particular we expect a mirror isomorphism 
\begin{equation} \label{eq:mirror.hh}
HH_{a}(\sD^{b}(X)) \cong HH_{a}(\FS((Y,\bw),\omega_{Y},\vol_{Y}))
\end{equation}
for all $a$ where the Hochschild homology can possibly be non-zero,
i.e. for all $a$ such that $-n \leq a \leq n$. It is also expected
that these categorical
homology groups have  geometric incarnations:
\begin{equation} \label{eq:hh=geom}
\begin{split}
HH_{a}(\sD^{b}(X)) & \cong \bigoplus_{p - q = a} H^{p}(X,\Omega^{q}_{X}),
\\[+0.3pc]
HH_{a}(\FS((Y,\bw),\omega_{Y},\vol_{Y})) & \cong H^{a+n}(Y,
Y_{-\infty}; \, \mathbb{C}).
\end{split}
\end{equation}
The first of the above identification follows from the work 
\cite{weibel-HC} of Weibel, while the second has been conjectured in
general, and proven in special cases in the works of Seidel (see
e.g. \cite{seidel-book,seidel-sh.as.hh}).  

Thus combining the conjectural mirror isomorphism \eqref{eq:mirror.hh} with this
geometric interpretation of Hochschild homology we will get a
conjectural isomorphism 
\begin{equation} \label{eq:mirror.hh.geometric}
H^{a+n}(Y,
Y_{-\infty}; \, \mathbb{C}) \cong \bigoplus_{p - q = a} H^{p}(X,\Omega^{q}_{X}) 
\end{equation}

\

\smallskip

\begin{rem} \label{rem:HPmirror} Since the mirror identification
  \eqref{eq:mirror.hh} comes from the mirror equivalence of
  categories, it is clear that a similar equivalence can also be
  formulated for the periodic cyclic homologies of the d$({\mathbb Z})$g
  categories $\sD^{b}(X)$ and
  $\FS((Y,\bw),\omega_{Y},\vol_{Y})$. Respectively the mirror
  identification \eqref{eq:mirror.hh.geometric} can be formulated for
  the de Rham cohomologies of $X$ and $(Y,\bw)$. In these cases we
  have natural \nc \ Hodge filtrations on each group. In the
  categorical setting the Hodge filtrations are encoded in the
  negative cyclic homologies $HC_{\bullet}^{-}(\sD^{b}(X))$ and
  $HC_{\bullet}^{-}((\FS((Y,\bw),\omega_{Y},\vol_{Y}))$ viewed as
  modules over $\mathbb{C}[[u]]$ (see
  \cite[Section~2.2.3]{kkp}). Homological mirror symmetry implies the
  existence of an isomorphism of $\mathbb{C}[[u]]$-modules
\[
HC_{\bullet}^{-}(\sD^{b}(X)) \cong  HC_{\bullet}^{-}(\FS((Y,\bw),\omega_{Y},\vol_{Y}))
\]
which after tensoring with $\mathbb{C}((u))$ induces an isomorphism of 
 $\mathbb{C}((u))$-vector spaces
\[
HP_{\bullet}(\sD^{b}(X)) \cong  HP_{\bullet}(\FS((Y,\bw),\omega_{Y},\vol_{Y})).
\]
From this point of view the isomorphism \eqref{eq:mirror.hh} is
recovered as the
induced isomorphism of specializations 
\[
HC_{\bullet}^{-}(\sD^{b}(X))/uHC_{\bullet}^{-}(\sD^{b}(X)) \cong 
HP_{\bullet}(\FS((Y,\bw),\omega_{Y},\vol_{Y}))/uHP_{\bullet}(
\FS((Y,\bw),\omega_{Y},\vol_{Y})).
\]
These mirror isomorphisms translate readily into the geometric
language. Recall that similarly to \eqref{eq:hh=geom} we have
identifications
\begin{equation} \label{eq:hp=geom}
\begin{split}
HP_{\bullet}(\sD^{b}(X)) & = H_{DR}^{\bullet}(X,\mathbb{C})\otimes \mathbb{C}((u)) \\
HP_{\bullet}(\FS((Y,\bw),\omega_{Y},\vol_{Y})), & =
H_{DR}^{\bullet}(Y, Y_{-\infty};\, \mathbb{C})\otimes \mathbb{C}((u)). \\
& 
\end{split}
\end{equation}
Furthermore, in geometric terms the $\mathbb{C}[[u]]$-module
$HC_{\bullet}^{-}(\sD^{b}(X))$ is identified with the Rees module of
the filtration 
\[
F^{a}_{\text{\nc}}H^{\bullet}_{DR}(X,\mathbb{C}) = \bigoplus_{p-q \geq a}
H^{p}(X,\Omega^{q}_{X})
\]
on the complex vector space $H^{\bullet}_{DR}(X,\mathbb{C})$, while the
$\mathbb{C}[[u]]$-module 
$HC_{\bullet}^{-}(\FS((Y,\bw),\omega_{Y},\vol_{Y}))$ is identified
with the Rees module of the filtration
\[
F^{a}_{\text{\nc}}H^{\bullet}_{DR}(Y,Y_{-\infty}; \mathbb{C}) = \bigoplus_{b \geq n+a}
H^{b}(Y,Y_{-\infty}; \mathbb{C})
\] 
on the complex vector space $H^{\bullet}(Y,Y_{-\infty}; \mathbb{C})$. 

The de Rham version of the Dolbeault mirror statement
\eqref{eq:mirror.hh.geometric} then becomes the statement that mirror
symmetry induces an isomorphism of the filtered complex vector spaces
$F^{\bullet}_{\text{\nc}}H^{\bullet}_{DR}(X,\mathbb{C})$ and
$F^{\bullet}_{\text{\nc}}H^{\bullet}_{DR}(Y,Y_{-\infty};
\mathbb{C})$. In fact,  more should be true. the induced isomorpism of
the algebraic vector bundles on $\mathbb{A}^{1}$ associated with the
respective  Rees modules should also intertwine the irregular
meromorphic connections describing the  \nc \ Hodge structures on both
sides. 
\end{rem}

\

\smallskip

\noindent
Going back to the mirror isomorphism \eqref{eq:mirror.hh.geometric}, we
are faced with the usual conundrum: simply identifying Hochschild (or
cyclic) homologies of the two mirror categories does not give us
matching of Hodge numbers. The comparison
\eqref{eq:mirror.hh.geometric} identifies the homological de Rham  grading on
the Landau-Ginzburg side with the $(p-q)$-folding of the Dolbeault
bigrading on the Fano side. 

The key to reconstructing the bigradings and thus extracting Hodge
numbers on both sides lies in the observation that, in the Fano case,
the Dolbeault bigrading also has a categorical interpretation. Indeed,
the nilpotent operator $c_{1}(K_{X})\cup (\bullet)$ acts on each
$HH_{a}(\sD^{b}(X)) = \oplus_{p-q} H^{p}(X,\Omega^{q}_{X})$, and so
induces a monodromy weight filtration centered at $a$. Since the
canonical class is anti-ample this filtration is given by the forms of
degree $\leq (p+q)$. In particular the dimensions of the graded pieces
for this monodromy weight filtration are precisely the Hodge numbers
$h^{p,q}(X)$ of the Fano variety $X$. Up to a sign, the nilpotent
operator $c_{1}(K_{X})\cup (\bullet) : H^{\bullet}_{DR}(X,\mathbb{C})
\to H^{\bullet}_{DR}(X,\mathbb{C})$ is just the logarithm of the
action of the Serre functor $S_{\sD^{b}(X)}$ on
$HH_{\bullet}(\sD^{b}(X)) \cong
H^{\bullet}_{DR}(X,\mathbb{C})$. Thus, this monodromy weight
filtration has a  categorical interpretation. But,  as we noted in 
Remark~\ref{rem:tameness}{\bfseries (ii)}, the Serre functor of
$\FS((Y,\bw),\omega_{Y},\vol_{Y})$ can be identified with the inverse
of the monodromy autoequivalence $T$. The logarithm of the action of
$T$ on $HH_{\bullet}(\FS((Y,\bw),\omega_{Y},\vol_{Y})) \cong
H^{\bullet+n}_{DR}(Y,Y_{c_{0}}; \mathbb{C})$ is just the nilpotent 
operator $N$ we considered above. Therefore the monodromy weight
filtration corresponding to $N$ is expected to have categorical
origin, and homological mirror symmetry predicts, the mirror matching of
Hodge numbers:
\begin{equation} \label{eq:match.hpq}
\boxed{
h^{p,q}(Y,\bw) = h^{p,n-q}(X)
}
\end{equation}
for all $p, q$. This prediction is still conjectural in general but the
case $(p,q) = (1,1)$ was recently proven by Przyjalkowski and Shramov
\cite{przyjalkowski.shramov} for all smooth Fano varieties.

\

\noindent
\punkt \label{sssec:ipq} {\bfseries The numbers $i^{p,q}(Y,\bw)$.}  
To simplify the discussion we will first assume that $\hor{D}_{\Ybar} =
\varnothing$, i.e. that $\bw : Y \to \mathbb{A}^{1}$ is proper. 

It is well known (see e.g. \cite{sabbah-twisted}) that the dimension
of the (Zariski) hypercohomology of the $\bw$-twisted de Rham complex
on $Y$ can be computed from the dimensions of the vanishing cohomology
for $\bw$:
\[
\dim_{\mathbb{C}} \mathbb{H}^{a}\left(Y_{\op{Zar}}, \left( \Omega^{\bullet}_{Y}, d +
    d\bw\wedge\right)\right) = \sum_{\lambda \in \mathbb{A}^{1}}
    \dim_{\mathbb{C}} \mathbb{H}^{a-1}\left(Y_{\lambda, \op{an}},
    \bphi_{\bw - \lambda} \mathbb{C}_{Y}\right), 
\]
where as usual $\bphi_{\bw - \lambda} \mathbb{C}_{Y}$ denotes the
perverse sheaf of vanishing cocycles for the fiber $Y_{\lambda}$.
From the works of Schmid and Steenbrink (see e.g. \cite{schmid} 
\cite[Section~11.2]{peters.steenbrink})  and Saito
\cite{saito-mhm} it is classically known that the constructible
complex $\bphi_{\bw - \lambda} \mathbb{C}_{Y}$ carries a structure of
a mixed Hodge module and so its cohomology is furnished with a
functorial mixed Hodge structure. 

Given a mixed Hodge structure $\mathcal{V}$ we will
write $i^{p,q}\mathcal{V}$ for the  $(p,q)$ Hodge number of the
$p+q$ weight graded piece $\op{gr}_{p+q}^{W}\mathcal{V}$. We now have
the following

\begin{defi} \label{defi:ipq} For a proper potential 
$\bw : Y \to \mathbb{A}^{1}$ on a quasi-projective variety $Y$ the
  {\bfseries Landau-Ginzburg Hodge numbers $i^{p,q}(Y,\bw)$} are
  defined by
\[
\boxed{
i^{p,q}(Y,\bw) := \sum_{\lambda \in \mathbb{A}^{1}}  \sum_{k}
i^{p,q+k}\mathbb{H}^{p+q-1}\left(Y_{\lambda}, 
    \bphi_{\bw - \lambda} \mathbb{C}_{Y}\right),
}
\]
where each vanishing cohomology $\mathbb{H}^{a}\left(Y_{\lambda}, 
    \bphi_{\bw - \lambda} \mathbb{C}_{Y}\right)$  is taken with its
    Schmid-Steenbrink mixed Hodge structyure. 
\end{defi}

\

\begin{rem} \label{rem:ipq}
{\bfseries (i)} \ The combination of Hodge numbers of different weight
pieces in this definition is motivated by mirror symmetry. In the
paper \cite{gkr} it was argued that for a Landau-Ginzburg mirror
of a general type complete-intersection $S$ in a toric variety, the
above definition of Hodge numbers reproduces the rotated Hodge diamond
of $S$.

\

\noindent
{\bfseries (ii)} \ The assumption that $\hor{D}_{\Ybar} = \varnothing$ above
was introduced solely for technical convenience and is not really
needed. If $\hor{D}_{\Ybar} \neq \varnothing$, we can still define
$i^{p,q}(Y,\bw)$ by setting 
\[
i^{p,q}(Y,\bw) := \sum_{\lambda \in \mathbb{A}^{1}}  \sum_{k}
i^{p,q+k}\mathbb{H}^{p+q-1}\left(\Ybar_{\lambda}, 
    \bphi_{\wbar - \lambda} R\jmath_{Y*}\mathbb{C}_{Y}\right)
\]
where $\jmath_{Y} : Y \hookrightarrow \Ybar$ is the natural
inclusion. 

\

\noindent
{\bfseries (iii)} \ It is very interesting to try and understand the
categorical meaning of the numbers $i^{p,q}(Y,\bw)$. At a first
glance, the definition of $i^{p,q}(Y,\bw)$ relies heavily on the
geometry since the information of the variety $Y$ and the potential
$\bw$ enter in an essential way in the construction of the pertinent
mixed Hodge structures. On the other hand, from the works of Shklyarov
\cite{shklyarov-nchodge} and Efimov \cite{efimov-cyclic} it is known
that the space $H^{\bullet}_{DR}(Y,\bw; \, \mathbb{C})$ together with
its \nc \ Hodge filtration admits a purely categorical
interpretation. Specifically, in \cite{efimov-cyclic} it is shown that
$H^{\bullet}_{DR}\left(Y_{\op{Zar}}, \left( \Omega^{\bullet}_{Y}((u)),
ud - d\bw\wedge \right)\right)$ is isomorphic to the periodic cyclic
homology $HP_{\bullet}(\MF(Y,\bw))$ of the d$({\mathbb Z}/2)$g category
of matrix factorizations of $\bw$, and that this isomorphism can be
chosen so that the irregular connection $\nabla^{\op{DR}}_{d/du} =
d/du + u^{-1}\sGr + u^{-2}\bw\cdot (\bullet)$ codifying the \nc
\ Hodge filtration on $H^{\bullet}_{DR}(Y,\bw; \, \mathbb{C})$ gets
identified with the connection $\nabla^{\op{cat}}_{d/du}$ from
\cite[Section~2.2.5]{kkp} used to define the categorical \nc \ Hodge
filtration on $HP_{\bullet}(\MF(Y,\bw))$.

In other words, the \nc \ Hodge filtration of a Landau-Ginzburg model
$(Y,\bw)$ admits a purely categorical interpretation.  In the case
when $\nabla^{\op{cat}}$ satisfies the \nc-opposedness axiom of
\cite{kkp} we can hope for more. In this case we expect that the pure
complex \nc \ Hodge structure on $HP_{\bullet}(\MF(Y,\bw))$ is
polarizable and that it admits a natural limit mixed twistor structure
(in the sense of \cite{sabbah-twistor}) which in turn is isomorphic to
the $\mathbb{Z}/2$-folding of the ordinary mixed Hodge structure on
the vanishing cohomology $\oplus_{\lambda \in \mathbb{A}^{1}}
\mathbb{H}^{\bullet-1}\left(Y_{\lambda}, \bphi_{\bw - \lambda}
  \mathbb{C}_{Y}\right)$. Concretely, we have a one parameter
deformation $\{\mathcal{A}_{t}\}_{t \in \mathbb{A}^{1}}$ of the
d$(\mathbb{Z}/2)$g category $\mathcal{A}_{1} := \MF(Y,\bw)$, where
$\mathcal{A}_{t}$ has the same objects and hom sets as $\MF(Y,\bw)$
but the composition, differential, and units in $\mathcal{A}_{t}$ are
scaled as $m_{\mathcal{A}_{t}} = t\cdot m_{\MF(Y,\bw)}$,
$d_{\mathcal{A}_{t}} = t\cdot d_{\MF(Y,\bw)}$, $1_{\mathcal{A}_{t}} =
t^{-1}\cdot 1_{\MF(Y,\bw)}$. The periodic cyclic homology of these
categories equipped with the connection $\nabla^{\op{cat}}$ in the
$u$-direction and with the Getzler-Gauss-Manin connection
\cite{getzler-gm} in the $t$-direction is a variation of twistor
$D$-modules. When the opposedness and polarizability properties hold,
e.g. for Landau-Ginzburg mirrors of toric Fano varieties, see
\cite{reichelt-sevenheck}, we can form the limit mixed twistor
 $D$-module for $t \to \infty$ and we conjecture that this mixed twistor
$D$-module is the one corresponding to the ordinary mixed Hodge
structure $\oplus_{\lambda \in \mathbb{A}^{1}}
\mathbb{H}^{\bullet-1}\left(Y_{\lambda}, \bphi_{\bw - \lambda}
  \mathbb{C}_{Y}\right)$. In the case of potentials given by tame
Laurent polynomials this conjecture is verified in \cite{sabbah}.  The
conjecture gives a categorical interpretation of the mixed Hodge
structure on vanishing cohomology (modulo Tate twists) and as a
consequence gives a categorical interpretation of the
$\mathbb{Z}/2$-folding of the numbers $i^{p,q}(Y,\bw)$.
\end{rem}

\

\smallskip

\noindent
\punkt \ {\bfseries Comparison conjectures} \ Because of their similar
behavior under the mirror correspondence we expect that the various
Landau-Ginzburg Hodge numbers are equal to each other:

\begin{conn} \label{conn:Hodge.equal} If $\bw : Y \to \mathbb{A}^{1}$
is an $n$-dimensional Landau-Ginzburg model which 
  admits a tame compactification, then 
\[
\boxed{
f^{p,q}(Y,\bw) = h^{p,q}(Y,\bw) = i^{p,q}(Y,\bw).
}
\]
\end{conn}

\

\smallskip

\noindent
Combined with the mirror matching \eqref{eq:match.hpq} the previous
conjecture predicts

\

\begin{conn} \label{conn:match.fpq} If $\left(X,\omega_{X},s_{X}\right) \bmid
\left((Y,\bw),\omega_{Y},\vol_{Y}\right)$ is a mirror pair of Fano
type, and if $\left((\Ybar, \wbar), D_{\Ybar}, \vol_{\Ybar}\right)$ is
a tame compactification of $\left((Y,\bw),\omega_{Y},\vol_{Y}\right)$,
then we have 
\[
\boxed{
f^{p,q}(Y,\bw) = h^{p,n-q}(X),
}
\]
for all $p, q$.
\end{conn}



\subsection{Mirrors of compactified Landau-Ginzburg
  models} \label{ssec:mirrorsZf} 

In this section we look more closely at the role that compactified
Landau-Ginzburg models play in mirror symmetry. In the setting where a
complex Landau-Ginzburg model $(Y,\bw)$ is the mirror of a symplectic
Fano variety $(X,\omega_{X})$, we give a mirror {\sf A}-model
interpretation of the Hodge information encoded in a compactification
$((\Ybar,\wbar),D_{\Ybar})$. This suggests that the mirror symmetry
bewteen $(Y,\bw)$ and $(X,\omega_{X})$ can be extended to a mirror
symmetry between $((\Ybar,\wbar),D_{\Ybar})$ and a one parameter
symplectic deformation of $(X,\omega_{X})$ which interpolates between
the Fukaya category of the symplectic Fano variety $(X,\omega_{X})$ and the
Fukaya category of the symplectic non-compact Calabi-Yau datum
$\left(X-D_{X},\omega_{X|X-D_{X}},\vol_{X-D_{X}}\right)$. We discuss
such an extension and give some evidence for its validity. This
picture is not new and has already been proposed and analyzed in one
form or another in the works of Seidel
\cite{seidel-book,seidel-genus2,seidel-AinftyI,seidel-AinftyII,seidel-ihes}
and Abouzaid et al \cite{abouzaid.seidel-viterbo,abouzaid-viterbo,
  abouzaidetal-punctured,abouzaidetal-punctured}. Our main
contribution here is to formulate a new procedure for reconstructing
the Hodge theory of $\wbar$-adapted logarithmic forms from the \nc
\ Hodge structure on the cohomology of the Landau-Ginzburg model
$(Y,\bw)$ or, in the mirror picture, from the {\sf A}-model \nc
\ Hodge structure on the cohomology of the Fano variety $X$.

\

\medskip

\punkt \ {\bfseries One parameter families of symplectic Fano
  varieties} \label{sss-symplectic.Fano} \ Let $(X,\omega_{X})$ be a
symplectic manifold underlying a smooth compact Fano variety of
$\dim_{\mathbb{C}} X = n$. Let $\boldsymbol{k}_{X}$ be a closed
$2$-form representing the canonical class $K_{X}$ and let $\kappa_{X}
\in H^{2}(X,\mathbb{Z})$ denote the first Chern class of $K_{X}$,
i.e. $\kappa_{X} = \left[\boldsymbol{k}_{X}\right] = c_{1}(K_{X})$.
Consider the (multivalued) family $\left\{\omega_{q}\right\}_{q \in
  \mathbb{C}}$ of complex $2$-forms $\omega_{q} := \omega_{X} +
\log(q) \boldsymbol{k}_{X}$ on $X$. In the regime when $|q| \to 1$ these are
complexified K\"{a}hler forms.

This is an affine-linear one-parameter family of symplectic structures
on $X$ which gives rise to a one-parameter variation of pure \nc \ Hodge
structures parametrized by the $q$-line. As discussed in
\cite[Section~3.1]{kkp} the de Rham part of such
variation is
encoded 
in a pair $(\Amod{H},\Amod{\nabla})$, where
\begin{itemize}
\item $\Amod{H} := H^{\bullet}(X,\mathbb{C})\otimes_{\mathbb{C}}
  \mathcal{O}_{\mathbb{A}^{2}}$ is a trivial algebraic
  $\mathbb{Z}/2$-graded vector bundle on the affine plane with
  coordinates $(u,q)$, with $\mathbb{Z}/2$-grading given by
\[
\begin{aligned}
\Amod{H}^{0} & = \left( \bigoplus_{\substack{k = n \\[+0.3pc] \mod \, 2}} 
H^{k}(X,\mathbb{C})
\right)\otimes_{\mathbb{C}} 
  \mathcal{O}_{\mathbb{A}^{2}},  \\[+0.5pc]
\Amod{H}^{1} & = \left( \bigoplus_{\substack{k = n+1 \\[+0.3pc] \mod \, 2}} 
H^{k}(X,\mathbb{C})
\right)\otimes_{\mathbb{C}}
  \mathcal{O}_{\mathbb{A}^{2}}.
\end{aligned}
\]
\item $\Amod{\nabla}$ is a meromorphic connection on $\Amod{H}$, with
  poles along the divisor $uq =0$, given by
\begin{equation} \label{eq:dubrovin.KX}
\left| \;
\begin{aligned}
\Amod{\nabla}_{\frac{\partial}{\partial u}} & := \frac{\partial}{\partial u}
+ u^{-2}\left(\kappa_{X}*_{q} \bullet\right) + u^{-1}\sGr \\[+0.5pc]
\Amod{\nabla}_{\frac{\partial}{\partial q}} & := \frac{\partial}{\partial q}
- q^{-1}u^{-1}\left(\kappa_{X}*_{q} \bullet\right),
\end{aligned}
\right.
\end{equation}
where  
\begin{description}
\item[\qquad $*_{q}$] denotes the quantum product corresponding to
  $\omega_{q}$, and 
\item[\qquad $\sGr : \Amod{H} \to \Amod{H}$] is the grading operator
  defined to be  
$\sGr_{|H^{k}(X,\mathbb{C})} := \frac{k-n}{2}
\op{id}_{H^{k}(X,\mathbb{C})}$. 
\end{description}
\end{itemize}

\

\begin{rem} \label{rem:affine.lines} 
More generally we have a variation of \nc \ Hodge
  structures over the whole complexified K\"{a}hler cone. 
The meromorphic connection defining the de Rham part of the variation
is the Dubrovin first structure connection
\cite{dubrovin-berlin,manin-frobenius}  which on each affine line 
$\omega + \log(q)\alpha$ is given by the formula
\[
\left| \; 
\begin{aligned}
\Amod{\nabla}_{\frac{\partial}{\partial u}} & := \frac{\partial}{\partial u}
+ u^{-2}\left(\kappa_{X}*_{q} \bullet\right) + u^{-1}\sGr \\[+0.5pc]
\Amod{\nabla}_{\frac{\partial}{\partial q}} & := \frac{\partial}{\partial q}
- q^{-1}u^{-1}\left([\alpha]*_{q} \bullet\right),
\end{aligned}
\right.
\]
with $*_{q}$ being the quantum product corresponding to $\omega + \log
(q)\alpha$. Traditionally in mirror symmetry one works with the line
of slope $\omega$ 
passing through the large volume limit, i.e. the line $0 +
\log(q)\omega$. 
This is the situation considered in \cite{fooo,fooo2} and in
\cite[Section~3.1]{kkp}. In contrast, here we need to work with a line
of slope $\boldsymbol{k}_{X}$, i.e. the line $\omega +
\log(q)\boldsymbol{k}_{X}$ which leads to the formula
\eqref{eq:dubrovin.KX}.
\end{rem}

\

\noindent
The particular affine linear deformation of the symplectic structure
$(X,\omega + \log(q)\boldsymbol{k}_{X})$ that we are considering has
many special properties even when compared to other affine linear
families. For instance, it is expected that the affine
one parameter deformation of the symplectic structure $(X,\omega +
\log(q)\boldsymbol{k}_{X})$ does not change the Fukaya category. From
the point of view of \nc \ Hodge theory this family is significant
because of the following simple observation 

\begin{lemma} \label{lem:restriction.to.lines} The restriction of
  $\left(\Amod{H},\Amod{\nabla}\right)$ to any non-vertical line $L$ through the
  origin in $\mathbb{A}^{2}$ is a meromorphic connection on the
  trivial $\mathbb{Z}/2$-graded vector bundle
  $H^{\bullet}(X,\mathbb{C})\otimes \mathcal{O}_{L}$, which has a
  first order pole at $0$, and monodromy around $0$ equal to
  $(-1)^{k}$ on the graded piece
  $H^{k}(X,\mathbb{C})\otimes \mathcal{O}$.
\end{lemma}
{\bfseries Proof.} From the formulas \eqref{eq:dubrovin.KX} we see that 
the $\Amod{\nabla}$-covariant derivative in the direction of the Euler
vector field ${\displaystyle \frac{u\partial}{\partial
    u}+\frac{q\partial}{\partial q}}$ is given by
\[
\Amod{\nabla}_{\frac{u\partial}{\partial u}+\frac{q\partial}{\partial q}} =
\frac{u\partial}{\partial u}+\frac{q\partial}{\partial q} + \sGr.
\]
Since  ${\displaystyle \frac{u\partial}{\partial
    u}+\frac{q\partial}{\partial q}}$ is tangent to any line through
the origin and is equal to the Euler vector field on any such line, we
get the statement of the lemma.

To explicate, choose a slope $v \neq 0$ and let $L_{v} \subset
\mathbb{A}^{2}$ be the line given by $u = vq$. The variable $q$ is the
natural coordinate on $L_{v}$, and so on $L_{v}$ we have $du =
vdq$. To shorten the notation, write $\sM := \kappa_{X}*_{q}(\bullet)$
for the operator of quantum multiplication by $\kappa_{X}$. Then
on $\mathbb{A}^{2}$ we have 
\[
\Amod{\nabla} = d + \left(u^{-2}\sM + u^{-1}\sGr\right)du  +
\left(-u^{-1}q^{-1} \sM\right)dq,
\]
and
\[
\begin{aligned}
\Amod{\nabla}_{|L_{v}} & = d + \left(v^{-2}q^{-2}\sM +
  v^{-1}q^{-1}\sGr\right)\cdot vdq  +
\left(-v^{-1}q^{-2}\sM\right)\cdot dq \\[+0.5pc]
& = d + \frac{\sGr}{q} dq.
\end{aligned}
\]
Thus $\Amod{\nabla}_{|L_{v}}$ is logarithmic at $0$ and has half
integer residues. This completes the proof of the lemma. \ \hfill
$\Box$

\

\smallskip

\noindent
Recall from \cite[Section~2.2.7]{kkp} that when viewed as \nc \ Hodge
structures, ordinary pure Hodge structures are given by meromorphic
connections on algebraic vector bundles over $\mathbb{A}^{1}$ that
have a first order pole at zero and monodromy $\pm 1$ on graded
pieces. Thus the pair $(\Amod{H},\Amod{\nabla})$ can be viewed as a
family of ordinary pure complex Hodge structures parametrized by $v
\in \mathbb{A}^{1}-\{0\}$. But this is exactly the type of data that our
Theorem~\ref{theo:double.degeneration} associates with a compactified
tame complex Landau-Ginzburg model.

\

\medskip

\punkt \ {\bfseries One parameter families of complex
  Landau-Ginzburg models} \label{sss-one.parameter.LG} \ Let
$\left(\left(\Ybar,\wbar),D_{\Ybar}\right)\right)$ be a compactified
tame complex Landau-Ginzburg model.  By
Theorem~\ref{theo:double.degeneration} the one parameter family  of
potentials $\left(\left(\Ybar,q\cdot\wbar),D_{\Ybar}\right)\right)$
gives rise to a variation of complex pure Hodge structures
parametrized by $q \in \mathbb{A}^{1}$. The de
Rham part of this variation is given by a pair
$\left(\Bmod{H},\Bmod{\nabla}\right)$, where 

\begin{itemize}
\item $\Bmod{H}$ is the coherent sheaf over $\mathbb{A}^{2}$
corresponding to the $\mathbb{C}[u,q]$-module
\[
\mathbb{H}^{\bullet}\left(\Ybar,
\left(\Omega_{\Ybar}^{\bullet}(\log D_{\Ybar},\wbar)[u,q], ud +
qd\wbar\right)\right).
\] 
Since by
Theorem~\ref{theo:double.degeneration}  the cohomology
$\mathbb{H}^{\bullet}\left(\Ybar, 
\left(\Omega_{\Ybar}^{\bullet}(\log D_{\Ybar},\wbar), ud +
qd\wbar\right)\right)$ has constant dimension for all $(u,q) \in
\mathbb{A}^{2}$ the sheaf  $\Bmod{H}$ is locally free. 
\item $\Bmod{\nabla}$ is the Gauss-Manin connection for the family of
  complexes of $\wbar$-adapted logarithmic forms. This is an algebraic
  meromorphic connection. For $u, q \neq 0$ the locally constant
  sections for $\Bmod{\nabla}$ are identified with the topological
  cohomology $H^{\bullet}(Y,Y_{-\infty};\mathbb{C})$ via the
  identifications from Lemma~\ref{lemma:cohomology.of.the.pair} and
  Lemma~\ref{lemma:barannikov.kontsevich}.
\end{itemize}

\

\noindent
By construction, the restriction of
$\left(\Bmod{H},\Bmod{\nabla}\right)$ on a line of the form $q = c$ is
the $u$-connection describing the Tate twist folding of the pure Hodge
structure on the vector space \linebreak
$\mathbb{H}^{\bullet}\left(\Ybar, \left(\Omega_{\Ybar}^{\bullet}(\log
D_{\Ybar},\wbar), ud + c\cdot d\wbar\right)\right)$. As explained in
\cite[Section~2.2.7]{kkp}, this means that the vector bundle
$\Bmod{H}_{|q = c}$ is the Rees module associated with the Hodge
filtration \linebreak $\mathbb{H}^{\bullet}\left(\Ybar,
\left(\Omega_{\Ybar}^{\bullet\geq p}(\log D_{\Ybar},\wbar), ud +
c\cdot d\wbar\right)\right)$ on $\mathbb{H}^{\bullet}\left(\Ybar,
\left(\Omega_{\Ybar}^{\bullet}(\log D_{\Ybar},\wbar), ud + c\cdot
d\wbar\right)\right)$, and that $\Bmod{\nabla}$ is a connection with
logarithmic singularity at $u = 0$ and monodromy $\pm 1$ on graded
pieces. In particular $\left(\Bmod{H},\Bmod{\nabla}\right)$ has a
logarithmic pole at $u = 0$.

\

\begin{rem} \label{rem:explicit.nabla} It will be useful to have an
  explicit formula for the Gauss-Manin connection $\Bmod{\nabla}$,
  similar to the formula \eqref{eq:dubrovin.KX}. We can write such a
  formula for the Gauss-Manin connection acting on the complex
  $\left(\Omega^{\bullet}_{\Ybar}\left(*D_{\Ybar}\right)[u,q], ud -
  qd\wbar\wedge \right)$:
\begin{equation} \label{eq:dubrovin.f}
\left| \;
\begin{aligned}
\Bmod{\nabla}_{\frac{\partial}{\partial u}} & := \frac{\partial}{\partial u}
+ u^{-2}\left(\wbar\cdot(\bullet)\right) + u^{-1}\bG \\[+0.5pc]
\Bmod{\nabla}_{\frac{\partial}{\partial q}} & := \frac{\partial}{\partial q}
- q^{-1}u^{-1}\left(\wbar\cdot(\bullet)\right),
\end{aligned}
\right.
\end{equation}
where  
$\bG$ is the grading operator
  defined to be  
$\bG := -\frac{p}{2}$
on $\Omega^{p}_{\Ybar}\left(*D_{\Ybar}\right)$.

This is the same formula that appears in the works of Shklyarov
\cite{shklyarov-nchodge} and Efimov \cite{efimov-cyclic}. Note however
that this formula does not preserve the subcomplex of $\wbar$-adapted
logarithmic forms since if $\alpha \in
\Omega^{\bullet}_{\Ybar}\left(\log D_{\Ybar},\wbar\right)$, the form
$f\alpha$ will not necessarily be in
$\Omega^{\bullet}_{\Ybar}\left(\log
  D_{\Ybar},\wbar\right)$. Therefore, we can not use these formulas
directly to describe the action of $\Bmod{\nabla}$ on $\Bmod{H}$. This
latter action is  a combination of the formulas \eqref{eq:dubrovin.f}
and the complicated limiting quasi-isomorphism in the proof of
Lemma~\ref{lemma:cohomology.of.the.pair}. 
\end{rem}

\

\medskip

\punkt \ {\bfseries Mirror symmetry for one parameter families} 
\label{sss-one.parameter.ms}  The formal similarity between the two
connections $\left(\Amod{H},\Amod{\nabla}\right)$ and
$\left(\Bmod{H},\Bmod{\nabla}\right)$ is very suggestive. We expect
that when the geometric data defining these connections is part of a
mirror pair, we should be able to go beyond a mere similarity and
identify the pairs $\left(\Amod{H},\Amod{\nabla}\right)$ and
$\left(\Bmod{H},\Bmod{\nabla}\right)$. More precisely we propose the
following conjecture.

\begin{conn} \label{conn:one.parameter.ms} Suppose
$\left(X,\omega_{X},s_{X}\right) \bmid 
\left((Y,\bw),\omega_{Y},\vol_{Y}\right)$  is a mirror pair of Fano
type. Then
\begin{itemize}
\item[{\normalfont\bfseries (i)}] The one parameter symplectic family
  $(X,\omega_{X} + \log(q)\boldsymbol{k}_{X})$ is mirrored into a one
  parameter complex family
  $\left(Y,q\cdot\bw)\right)$ of deformations of 
  $(Y,\bw)$.
\item[{\normalfont\bfseries (ii)}] The homological mirror
  correspondence induces an isomorphism
\[
\left(\Amod{H},\Amod{\nabla}\right) \cong
\left(\Bmod{H},\Bmod{\nabla}\right)
\]
of meromorphic connections on $\mathbb{A}^{2}$.  
\end{itemize}
\end{conn}

\

\noindent
The attentive reader will notice that the part {\bfseries (i)} of this
conjecture relies on a geometric one parameter perturbation of a Fano
mirror pair but does not involve a compactification of the
Landau-Ginzburg side of the pair. On the other hand, at least the {\sf
  B}-side of Conjecture~\ref{conn:one.parameter.ms}{\bfseries
  (ii)}  depends on a tame compactification
$\left(\left(\Ybar,\wbar\right),D_{\Ybar}\right)$ of the
Landau-Ginzburg model $(Y,\bw)$. 

\

\noindent
Nevertheless, in Section~\ref{sss-construction} we will argue that
part {\bfseries (ii)} of Conjecture~\ref{conn:one.parameter.ms} is in
fact a consequence of the homological mirror symmetry conjecture for
the Fano pair itself.  In other words: the existence of the tame
compactification matters, while the choice of a particular
compactification is not important. Indeed, as explained in
\cite[Sections~2.2.2, 3.1, and 3.2]{kkp}, the pure \nc \ Hodge
structures on $H^{\bullet}(X,\mathbb{C})$ and
$H^{\bullet}(Y,Y_{\op{sm}};\mathbb{C})$ can be defined intrinsically
in terms of the categories $\Fuk(X,\omega_{X})$ and $\MF(Y,\bw)$
respectively. Since homological mirror symmetry identifies these two
categories, it follows that the ${\sf A}$-model \nc \ Hodge structure
on $H^{\bullet}(X,\mathbb{C})$ will be isomorphic to the ${\sf
  B}$-model \nc \ Hodge structure on
$H^{\bullet}(Y,Y_{\op{sm}};\mathbb{C})$. Via these identifications
Conjecture~\ref{conn:one.parameter.ms}{\bfseries (ii)} reduces to
checking that the two parameter meromorphic connections
$\Amod{\nabla}$ and $\Bmod{\nabla}$ can be reconstructed from the \nc
\ Hodge structures on $H^{\bullet}(X,\mathbb{C})$ and
$H^{\bullet}(Y,Y_{\op{sm}};\mathbb{C})$ respectively.  To that end, in
Section~\ref{sss-construction} we describe a general method for
constructing a two parameter meromorphic connection from a pure \nc
\ Hodge structure.

\

\medskip

\punkt \ {\bfseries Mirrors of tame compactifications of
  Landau-Ginzburg models.} \label{sss-nc.pencils} Before we proceed
with the construction in Section~\ref{sss-construction}, it is
instructive to examine more closely the apparent mismatch in the
information contained in the one parameter mirror symmetry
$\left(X,\omega_{X} + \log(q)\boldsymbol{k}_{X}\right) \bmid
\left(Y,q\cdot\bw\right)$ and in the tame compactification of the
Landau-Ginzburg model. This mismatch is ultimately a reflection of the
fact that the one parameter deformations $\left(X,\omega_{X} +
\log(q)\boldsymbol{k}_{X}\right)$ and $\left(Y,q\cdot\bw\right)$ only
perturb one direction of the mirror symmetry: going from the {\sf
  A}-model on the Fano side to the {\sf B}-model on the
Landau-Ginzburg side.

Since the choice of a tame compactification
$\left(\left(\Ybar,\wbar\right),D_{\Ybar}\right)$ is a choice
additional data on the Landau-Ginzburg side, its mirror partner will
necessarily depend on the choice of some additional data on the Fano
side. A clue of what this additional data should be, appears  in
the works of Seidel \cite{seidel-book,seidel-genus2} where the
one parameter deformation $\left(X,\omega_{X} +
\log(q)\boldsymbol{k}_{X}\right)$ is interpretted intrinsically in
categorical terms. The relevant key fact from
\cite{seidel-book,seidel-genus2}  is the statement that the 
family of Fukaya categories $\Fuk\left(X, \omega_{X} +
\log(q)\boldsymbol{k}_{X}\right)$ has a well defined limit as $q\to
0$, namely the Fukaya category of the symplectic manifold underlying
the non-compact Calabi-Yau
$X - D_{X}$.  

To simplify notation write $U := X - D_{X}$ for the complement of
$D_{X}$, $\omega_{U} := \omega_{X|U}$ for the restriction of the
symplectic structure to $U$, and $\vol_{U} = 1/s_{X}$ for the
holomorphic volume form corresponding to the anticanonical section
$s_{X}$. As explained in \cite{seidel-book,seidel-genus2} (see also
\cite{auroux-anticanonical}) the $\mathbb{Z}$-graded $A_{\infty}$
category $\mycal{F}_{0} = \Fuk(U,\omega_{U},\vol_{U})$ admits a
natural one-parameter deformation $\{\mycal{F}_{q}\}$ as a
$\mathbb{Z}/2$-graded $A_{\infty}$ category. By construction
$\mycal{F}_{q}$ has the same objects and morphisms as $\mycal{F}_{0}$
but the $A_{\infty}$ operations $m_{k}^{q}$ in $\mycal{F}_{q}$ are
$q$-perturbations of the $A_{\infty}$ operations $m_{k}^{0}$ in
$\mycal{F}_{0}$ where the correction term $q^{a}$ comes with a
coefficient counting not pseudo-holomorphic discs in $U$ but rather
disks in $X$ that intersect the boundary divisor $D_{X}$ at $a$
points\footnote{Making this precise is quite subtle (see
  \cite{wehrheim-woodward,seidel-genus2}) and requires a version of
  the Fukaya category which is linear over $\mathbb{C}$ (rather than a
  Novikov field).  In \cite{seidel-genus2} such a version is built out
  of balanced (rather than arbitrary) Lagrangians in $U$. We thank
  Denis Auroux for illuminating explanations of this subtlety.}.  Now,
a comparison with the standard construction \cite{fooo,fooo2} of the
Fukaya category identifies $\mycal{F}_{q}$ for $q \neq 0$ the with the
$\mathbb{Z}/2$ category $\Fuk(X,\omega_{X} +
\log(q)\boldsymbol{k}_{X})$.

Thus we get a streamlined categorical (or \nc \ geometric)
interpretation of the {\sf A}-model data associated with the Fano
geometry $\left(X,\omega_{X},s_{X}\right)$. In summary  
Seidel's analysis shows that:
\begin{itemize}
\item from the point of view of \nc \ geometry, the
primordial object is the $\mathbb{Z}$-graded Fukaya category
$\mycal{F}_{0} = \Fuk\left(U,\omega_{U},\vol_{U}\right)$; 
\item the data of a symplectic 
compactification $(U,\omega_{U}) \subset (X,\omega_{X})$ with
anti-canonical boundary $D_{X} = X - U$ corresponds to a $q$-deformation
$\mycal{F}_{q} = \Fuk(X,\omega_{X} +
\log(q)\boldsymbol{k}_{X})$ of   $\mycal{F}_{0}$ as a
$\mathbb{Z}/2$-graded  Calabi-Yau category. 
\end{itemize}
To put it differently, the symplectic anti-canonical compactification
$(U,\omega_{U}) \subset (X,\omega_{X})$ is encoded in a one parameter
degeneration of the Fukaya category $\mycal{F}_{1} =
\Fuk(X,\omega_{X})$ of the compact symplectic Fano $(X,\omega_{X})$ to
the Fukaya category $\mycal{F}_{0} =
\Fuk\left(U,\omega_{U},\vol_{U}\right)$ of the symplectic non-compact
Calabi-Yau $(U,\omega_{U})$. 

This categorical interpretation of the compactification of $U$ has a
natural mirror incarnation. The non-compact symplectic Calabi-Yau
$(U,\omega_{U},\vol_{U})$ has a complex non-compact Calabi-Yau mirror
$Y$, constructed say by the SYZ prescription as in
\cite{auroux-anticanonical}. Homological mirror symmetry predicts
that $\mycal{F}_{0}$ is equivalent to the category
$\sD^{b}_{c}(Y)$. The one parameter deformation $\mycal{F}_{q}$ of
$\mycal{F}_{0} = \sD^{b}_{c}(Y)$ then correpsonds to a class in the
Hochschild cohomology $HH^{\bullet}(\mycal{F}_{0}) =
HH^{\bullet}(\sD^{b}_{c}(Y))$.  Since $\mycal{F}_{q}$ is only a
$\mathbb{Z}/2$-graded deformation, this Hochschild cohomology class
will have a non-trivial component in $HH^{0}$, i.e. will give us a
well defined element $\bw \in H^{0}(Y,\mathcal{O}_{Y})$. If we assume
for symplicity that the boundary divisor $D_{X}$ is smooth, then the
Fano/Landau-Ginzburg homological mirror symmetry conjecture will
identify $\mycal{F}_{q}$ with $\MF(Y,q\cdot \bw)$ for $q \neq 0$. If
we interpret $\MF(Y,q\cdot \bw)$ as a coproduct of the derived categories
of singularities of  the singular fibers of $q\cdot \bw$, we see that
this identification will specialize correctly when $q \to 0$. The
category $\mycal{F}_{q}$ specializes to $\mycal{F}_{0}$ while
$\MF(Y,q\cdot \bw)$ specializes to the compactly supported derived
category of singularities of the derived fiber of the zero function on
$Y$, which is readily identified with $\sD^{b}_{c}(Y)$.

The upshot of the previous discussion is that the mirror one parameter
families 
\[
(X,\omega_{X} + \log(q) \boldsymbol{k}_{X}) \bmid (Y, q\cdot \bw)
\]
arising from the Fano mirror pair 
$(X,\omega_{X},s_{X}) \bmid ((Y,\bw),\omega_{Y},\vol_{Y})$ 
have a natural homolgical interpretation as families of (term by term
equivalent) categories 
\[
\{\mycal{F}_{q}\} = \{\MF(Y,q\cdot \bw)\}, 
\]
where the family on the left hand side is the Seidel
$\mathbb{Z}/2$-graded deformation of $\Fuk(U,\omega_{U},\vol_{U})$
corresponding to the compactification $(U,\omega_{U}) \subset
(X,\omega_{X})$. 

This interpretation allows us to reverse the process and identify
the mirror information corresponding to a tame compactification of
$Y$. If we choose a tame compactification
$\left((\Ybar,\wbar),D_{\Ybar}\right)$ and also choose an extension
$\omega_{\Ybar}$ of
the symplectic form $\omega_{Y}$, then we can apply Seidel's analysis
to the symplectic anti-canonical compactification $(Y,\omega_{Y})
\subset (\Ybar,\omega_{\Ybar})$. Since by the tameness assumption
$D_{\Ybar}$ is an anti-canonical divisor, the same reasoning shows that
this compactification is encoded in a one parameter deformation of the
$\mathbb{Z}$-graded category $\Fuk(Y,\omega_{Y},\vol_{Y})$ to the
$\mathbb{Z}/2$-graded category $\Fuk\left(\Ybar,\omega_{\Ybar} +
\log(r)\boldsymbol{k}_{\Ybar}\right)$. Again the degree zero piece of
the Hochschild cohomology class governing this deformation will give
us a holomorphic function $\boldsymbol{{\sf v}} : U \to
\mathbb{A}^{1}$. In fact the description of $\boldsymbol{{\sf v}}$ in
terms of the weighted disk counting on $\Ybar$ relative to the
boundary $D_{\Ybar}$ also predicts that  $\boldsymbol{{\sf
    v}}$ has first order poles along $D_{X}$ and so
$\boldsymbol{{\sf v}} =
s/s_{X}$ for some anti-canonical section $s \in H^{0}(X,K_{X}^{-1})$. 
This can be packaged in the following conjecture.

\begin{conn} \label{conn:tame.comp} Suppose
$\left(X,\omega_{X},s_{X}\right) \bmid 
\left((Y,\bw),\omega_{Y},\vol_{Y}\right)$  is a mirror pair of Fano
type. Then
\begin{itemize}
\item[{\normalfont\bfseries (i)}] a choice of a tame
  compactification
  $\left((\Ybar,\wbar),D_{\Ybar},\omega_{\Ybar}\right)$ of the
  Landau-Ginzburg side gives rise to a compactified Fano mirror pair
\[
\left(X,\omega_{X},s_{X},\boldsymbol{{\sf v}}\right) \bmid
\left((\Ybar,\wbar),D_{\Ybar},\omega_{\Ybar}\right),
\]
where $\boldsymbol{{\sf v}}$ is a meromorphic function on $X$ with a
first order pole along $D_{X}$.
\item[{\normalfont\bfseries (ii)}] The Fano/Landau-Ginzburg homological mirror
  correspondence induces equivalences 
\[
\begin{split}
\Fuk(X,\omega_{X} + \log(q) \boldsymbol{k}_{X}) & \cong
\sD^{b}_{c}(Y,q\cdot \bw) \\
\sD^{b}_{c}(U,r\cdot \boldsymbol{{\sf v}}) & \cong 
\Fuk(\Ybar,\omega_{\Ybar} + \log(r) \boldsymbol{k}_{\Ybar})
\end{split}
\]
of one parameter families of categories.
\end{itemize}
\end{conn}

\

\noindent
Note that the geometric ingredients of the compactified mirror pair 
\[
\left(X,\omega_{X},s_{X},\boldsymbol{{\sf v}}\right) \bmid
\left((\Ybar,\wbar),D_{\Ybar},\omega_{\Ybar}\right).
\]
from Conjecture~\ref{conn:tame.comp}{\bfseries
  (i)} appear now in a symmetric fashion in the two sides of the pair.
In particular we expect {\sf A}-model data on one side
to be mirror to the {\sf B}-model data on the other. 

It is clear also that the statement of
Conjecture~\ref{conn:tame.comp}{\bfseries (ii)} is but one facet of
the homological mirror correspondence one should associate with the
compactified mirror pair.  The full homological mirror conjecture will
involve various equivalences of categories generalizing the
equivalences described in Table~\ref{table:hms.fano} for an ordinary
(non-compactified) Fano mirror pair. It is possible to list all these
equivalences but the list is somewhat cumbersome.  Very recently
Seidel found \cite{seidel-ihes} a uniform conceptual way for capturing
the homological content of either side of the compactified mirror pair
and gave a clean formulation of the complete homological mirror
conjecture for the pair Conjecture~\ref{conn:tame.comp}{\bfseries (i)}
in terms of an equivalence of categories equipped with \nc
\ anti-canonical pencils.

It is very interesting to understand how the two meromorphic
connections from Conjecture~\ref{conn:one.parameter.ms}{\bfseries (ii)}
arise directly as Hodge theoretic data associted with these \nc
\ pencils but we will leave this for future investigations.

\

\medskip

\punkt \ {\bfseries Construction of meromorphic connections over
  ${\mathbb A}^{2}$.}  \label{sss-construction} \quad Suppose
$(H^{\bullet},\nabla)$ is the de Rham part of a pure \nc \ Hodge structure. In
this section we explain how, under some mild technical assumptions,
the pair $(H^{\bullet},\nabla)$ gives rise to a meromorphic connection over the
affine plane $\mathbb{A}^{2}$. 

To keep track of the various copies of the affine line and the affine
plane appearing in the construction, we will indicate the coordinates
on these lines and planes as subscripts. Thus $\mathbb{A}^{1}_{u}$
will denote the affine line with coordinate $u$,
$\mathbb{A}^{2}_{(u,q)}$ will denote the affine plane with coordinates
  $(u,q)$, etc. By definition (see
  \cite[Section~2.1.4]{kkp}) the pair $(H,\nabla)$ is the de Rham part
  of a pure \nc \ Hodge structure if it satisfies:  
\begin{description}
\item[\quad $H^{\bullet}$] is a $\mathbb{Z}/2$-graded algebraic vector bundle on
  the affine line $\mathbb{A}^{1}_{u}$, and 
\item[\quad $\nabla$] is a meromorphic
  connection on $H^{\bullet}$, which has at most a regular singularity at $u =
  \infty$, at most a second order pole at $u = 0$, and no other
  singularities in $\mathbb{A}^{1}_{u}$.
\end{description}

\

\noindent
View $\mathbb{A}^{1}_{u}$ as the $u$-axis in the plane
$\mathbb{A}^{2}_{(u,q)}$. Our goal is to extend $H^{\bullet}$ to a holomorphic
bundle $\stwo{H}^{\bullet}$ on all of $\mathbb{A}^{2}_{(u,q)}$, and $\nabla$ to a
meromorphic connection $\stwo{\nabla}$ on $\stwo{H}^{\bullet}$ over
$\mathbb{A}^{2}_{(u,q)}$ so that $\stwo{\nabla}$ has poles only at $uq = 0$
and has logarithmic singularities along $q = 0$. We will carry this
out in two steps:

\

\smallskip

\noindent {\bfseries Step 1.} \ Start with the connection
$(H^{\bullet},\nabla)$ on $\mathbb{A}^{1}_{u}$. Write $(H,\nabla)$ for
the underlying ungraded algebraic vector bundle with connection. Since
by assumption $\nabla$ has a regular singularity at $u = \infty$ we
can consider the Deligne extension $(\mathcal{H},\nabla)$ of
$(H,\nabla)$ (see e.g. \cite[Chapter~II.5]{deligne} or
\cite[Corollary~II.2.21]{sabbah-frobenius}). The bundle $\mathcal{H}$
is the algebraic vector bundle on \linebreak $\mathbb{P}^{1}_{u} =
\mathbb{A}^{1}_{u}\cup \{\infty\}$ which is uniquely characterized by
the properties that at $\infty$ the connection $\nabla$ has a
logarithmic pole on $\mathcal{H}$, and that the residue
$\op{Res}^{\mathcal{H}}_{\infty}(\nabla) :
\mathcal{H}_{\infty} \to \mathcal{H}_{\infty}$ is
a nilpotent, grading preserving endomorphism of the fiber of
$\mathcal{H}$ at $\infty$.

The bundle $\mathcal{H}$ decomposes into a direct sum of line bundles
$\mathcal{H} = \oplus_{k = 1}^{r} \mathcal{O}_{\mathbb{P}^{1}}(d_{k})$
and so it admits a natural decreasing biregular filtration by subbundles
\[
F^{i}\mathcal{H} = \bigoplus_{d_{k} \geq i}
\mathcal{O}_{\mathbb{P}^{1}}(d_{k}), \qquad i \in \mathbb{Z}.
\] 
The restrictions $F^{i}H := F^{i}\mathcal{H}_{|\mathbb{A}^{1}_{u}}$
give a $\mathbb{Z}$-labeled filtration of $H$ by holomorphic
subbundles. 

For any complex number $v \in \mathbb{C}$ consider the Rees bundle
$\xi(H_{v},F^{\bullet}H_{v}) \to \mathbb{A}^{1}_{q}$ associated with this
filtration, \cite{simpson-hf}.  The bundle
$\xi(H_{v},F^{\bullet}H_{v})$ is defined as the locally free sheaf on
$\mathbb{A}^{1}_{q}$ associated with the $\mathbb{C}[q]$-submodule
$\sum_{i} q^{-i}F^{i}H_{v} \subset H_{v}[q,q^{-1}]$. By construction
$\xi(H_{v},F^{\bullet}H_{v})$ is a $\mathbb{C}^{\times}$-equivariant
vector bundle on $\mathbb{A}^{2}_{q}$ for the scaling action of
$\mathbb{C}^{\times}$ on the $q$-line.  Allowing $v$ to vary we get a
Rees bundle $\xi(H,F^{\bullet}H) \to \mathbb{A}^{2}_{(v,q)}$ which is
algebraic and equivariant for the $\mathbb{C}^{\times}$-action
$\lambda\cdot (v,q) := (v,\lambda q)$. By construction we have
canonical identifications
\[
\begin{aligned}
\xi(H,F^{\bullet}H)_{(v,1)} & \cong H_{v},
\\[+0.5pc]
 \xi(H,F^{\bullet}H)_{(v,0)} & \cong
 \op{gr}_{F^{\bullet}H_{v}} H_{v}.
\end{aligned}
\]
Since the filtration $F^{\bullet}H$ was compatible with the grading on
$H$, we get a natural $\mathbb{Z}/2$-grading on
$\xi(H,F^{\bullet}H)$. Similarly, since $F^{\bullet}H$ arose from the
Deligne extension of $(H,\nabla)$, the meromorphic connection $\nabla$
on $H = \xi(H,F^{\bullet}H)_{|\mathbb{A}^{1}_{v}\times \{ 1 \}}$
extends to a well defined meromorphic connecion on
$\xi(H,F^{\bullet}H)$ over the plane $\mathbb{A}^{2}_{(v,q)}$ which
has poles on $vq = 0$ and on each line $\{v\}\times \mathbb{q}$, $v
\neq 0$ has monodromy $+1$ on the even graded piece of
$\xi(H_{v},F^{\bullet}H)$ and $-1$ on the odd graded piece. To
simplify notation we will write $\sone{H}^{\bullet}$ for the
$\mathbb{Z}/2$-graded $\xi(H,F^{\bullet}H)$ on
$\mathbb{A}^{2}_{(v,q)}$, and will write $\sone{\nabla}$ for the
extension of the connection $\nabla$.

\

\medskip

\noindent
{\bfseries Step 2.} \ In this step we will modify the
$\mathbb{Z}/2$-graded bundle with connection $\left(
\sone{H}, \sone{\nabla}\right)$ to ensure that its
monodromy is $\pm 1$ not on vertical lines but rather on lines through
the origin.

Consider an affine plane $\mathbb{A}^{2}_{(u,q)}$ with coordinates
$(u,q)$. Let $\mathfrak{s} : S \to \mathbb{A}^{2}_{(u,q)}$ be the
blow-up of $\mathbb{A}^{2}_{(u,q)}$ at the origin $(u,q) =
(0,0)$. The surface $S$ is glued out of two affine charts
$\mathbb{A}^{2}_{(v,q)}$ and $\mathbb{A}^{2}_{(u,w)}$ via the gluing
map $u = vq$, $w = 1/v$. In particular $\mathbb{A}^{2}_{(v,q)}$ embeds
as a Zariski open subset in $S$ and we have a commutative diagram of
surfaces
\[
\xymatrix{ \mathbb{A}^{2}_{(v,q)} \ar@{^{(}->}[r]^-{i}
  \ar[rd]_-{\mathfrak{p}} & S \ar[d]^-{\mathfrak{s}} \\
& \mathbb{A}^{2}_{(u,q)} 
}
\]
where $i : \mathbb{A}^{2}_{(v,q)} \hookrightarrow S$ denotes the
inclusion, and  $\mathfrak{p}$ is the map
\[
\mathfrak{p} : \mathbb{A}^{2}_{(v,q)} \to \mathbb{A}^{1}_{(u,q)},
\qquad  (u,q) = \mathfrak{p}(v,q) = (vq,q).
\] 
Note that $S - \mathbb{A}^{2}_{(v,q)} = \mathbb{A}^{1}_{u}$ and that a
point $u \neq 0 \in \mathbb{A}^{1}_{u} \subset S$ is a limiting point
in $S$ completing the hyperbola $\{ (v,q)
\ | vq = u \} \subset \mathbb{A}^{2}_{(v,q)}$  to a copy $C^{u}$ of 
$\mathbb{A}^{1}$ embedded in $S$. 

Now observe that if we restrict $(\sone{H},\sone{\nabla})$ to the
hyperbola $vq = u$, the restricted connection has a regular
singularity as $v \to \infty$. Therefore we have a canonical Deligne
extension of $\sone{H}_{|vq = u}$ to an algebraic vector bundle on
$C^{u}$. This process depends algebraically on $u$ and so gives an
extension of $\sone{H}$ to an algebraic vector bundle $\Xi$ on the
punctured surface $S - \{ x\}$, where $x \in \mathbb{A}^{2}_{(u,w)}$
is the point with coordinates $u = 0$, $w = 0$. 

Next observe that since the surface $S$ is smooth, any vector bundle
$V$ on $S- \{ x\}$ will extend to a (necessarily unique) vector bundle
on $S$. Indeed, choose a torsion free coherent sheaf $F$ on $S$ which
restricts to $V$ on $S- \{ x\}$. For instance, if $\jmath$ denotes the
inclusion of $S- \{ x\}$ in $S$, we can take $F$ to be the
intersection $\cap K$ of all coherent subsheaves $K \subset
\jmath_{*}V$ such that $V \subset \jmath^{*}K$. Since $V$ is locally
free, the double dual $F^{\vee\vee}$ will also restrict to $V$. Being
the dual of a coherent sheaf $F^{\vee\vee}$ is automatically
reflexive, and by Auslander-Buchsbaum theorem can only fail to be
locally free in codimension three. Thus $F^{\vee\vee}$ is a locally
free sheaf which extends $V$ to $S$. The uniqueness of the extension
follows again from the fact that $x$ is a smooth point and so the
local ring $\mathcal{O}_{S,x}$ satisfies the Serre condition $S_{2}$.

Let $\widetilde{\Xi}$ be the unique extension of $\Xi$ to $S$. To
complete the construction we will need to know that $\widetilde{\Xi}$
satisfies a descent property for the morphism $\mathfrak{s} : S
\to \mathbb{A}^{2}_{(u,q)}$. Let $E \subset S$ denote the exceptional
  $\mathbb{P}^{1}$ of the blow-up morphism $\mathfrak{s} : S
\to \mathbb{A}^{2}_{(u,q)}$. With this notation we have the following

\begin{defi} \label{defi:good.extension} We will say that an \nc
  \ Hodge filtration $(H^{\bullet},\nabla)$ is {\bfseries extendable}
  if the restriction of the algebraic vector bundle $\widetilde{\Xi}$
  to $E$ is holomorphically trivial.
\end{defi}

\

\noindent
Note that if $(H^{\bullet},\nabla)$ is extendable, then by the
projection formula this extension $\widetilde{\Xi}$ is canonically a
pullback of a vector bundle on $\mathbb{A}^{2}_{(u,q)}$, namely
$\mathfrak{s}_{*} \widetilde{\Xi}$ is a vector bundle and
$\widetilde{\Xi} \cong \mathfrak{s}^{*}\mathfrak{s}_{*}
\widetilde{\Xi}$.

In particular, if the \nc \ Hodge filtration $(H^{\bullet},\nabla)$ is
extendable, we get a well defined holomorphic bundle $\stwo{H} :=
\mathfrak{s}_{*}\widetilde{\Xi}$ on $\mathbb{A}^{2}_{(u,q)}$. The
meromorphic connection $\sone{\nabla}$ on $\sone{H}$ is holomorphic on
the open set $vq \neq 0$ and so can be viewed as a meromorphic
connection $\stwo{\nabla}$ on $\stwo{H}$ with poles on $uq = 0$.
Altogether we have proven the following

\

\begin{lemma} \label{lem:good.extension}
Let $(H^{\bullet},\nabla)$ be an extendable  \nc \ Hodge filtration, then
$(H^{\bullet},\nabla)$ gives rise to a $\mathbb{Z}/2$-graded
meromorphic connection $\left(\stwo{H}^{\bullet},\stwo{\nabla}\right)$
on $\mathbb{A}^{2}_{(u,q)}$, such that
\begin{itemize}
\item $\stwo{\nabla}$ is holomorphic away from $uq  = 0$;
\item $\stwo{\nabla}$ has at most a logarithmic pole along $u = 0$, and
  a pole of order $\leq 2$ along $q = 0$;
\item The restriction of $\left(\stwo{H}^{0},\stwo{\nabla}\right)$ to
  a line through the origin has trivial monodromy, while the
  restriction of $\left(\stwo{H}^{1},\stwo{\nabla}\right)$ to
  a line through the origin has monodromy $(-1)$.
\end{itemize}
\end{lemma}

\

\noindent
The discussion in Section~\ref{sss-symplectic.Fano} shows that the
extendability assumption in the previous lemma holds for the de Rham
part of the \nc \ Hodge structure associated with a symplectic Fano
variety:

\begin{cor} \label{cor:Amodel.is.good} Let $(X,\omega_{X})$ be a
  symplectic manifold underlying a smooth Fano variety of complex
  dimension $n$. Let $*_{1}$ denote the quantum product corresponding
  to the symplectic form $\omega_{X}$. Then
the {\sf A}-model \nc \ Hodge filtration
\[
\left(\Anc{H}^{\bullet},\Anc{\nabla}\right) := 
\left( H^{\bullet}(X,\mathbb{C})\otimes
\mathcal{O}_{\mathbb{A}^{1}_{u}}, \ d +
(u^{-2}(\kappa_{X}\otimes_{1}(\bullet)) +  u^{-1}\sGr)du\right) 
\]
for the \nc \ Hodge structure on the cohomology of $(X,\omega_{X})$ is
extendable and
$\left(\stwo{H}^{\bullet},\stwo{\nabla}\right)$ reconstructs the
standard $q$-variation of \nc \ Hodge structures for the symplectic
manifold $(X,\omega_{X})$. That is,  we have
a canonical identification
\[
\left(\stwo{\left(\Anc{H}\right)}^{\bullet},
\stwo{\left(\Anc{\nabla}\right)}\right)   
 =
  \left(\Amod{H}^{\bullet},\Amod{\nabla}\right),
\] 
where
  $\left(\Amod{H}^{\bullet},\Amod{\nabla}\right)$ is the 
  connection defined in \eqref{eq:dubrovin.KX}. 
\end{cor}
{\bfseries Proof.} \ Follows immediately from the two step
construction above and by 
 Lemma~\ref{lem:restriction.to.lines}
and Lemma~\ref{lem:good.extension}. \ \hfill $\Box$

\

\medskip

\noindent
In particular Corollary~\ref{cor:Amodel.is.good} shows that the one
parameter mirror symmetry
Conjecture~\ref{conn:one.parameter.ms}{\bfseries (ii)} is 
equivalent to the extendability property for the {\sf B}-model \nc \
Hodge filtration. More precisely, suppose $(Y,\bw)$ is a complex
Landau-Ginzburg model. Consider the {\sf B}-model \nc \
Hodge filtration for the cohomology of $(Y,\bw)$:
\[
\left(\Bnc{H}^{\bullet},\Bnc{\nabla}\right) =
\left(\mathbb{H}^{\bullet}\left(\Omega_{Y}^{\bullet}[u], ud -
    d\bw\wedge\right), \ d + \left( u^{-2}\left(\bw\cdot
      (\bullet)\right) + u^{-1}\bG\right)\right).
\]
Here $\bG$ is the grading operator of multiplication by $-p/2$ on
$\Omega^{p}_{Y}$. 
When $(Y,\bw)$ is the mirror of a symplectic Fano variety
$(X,\omega_{X})$, homological
mirror symmetry identifies the \nc \ Hodge filtrations on the
cohomology, i.e. gives an isomorphism 
\[
\left(\Anc{H}^{\bullet},\Anc{\nabla}\right) \cong
\left(\Bnc{H}^{\bullet},\Bnc{\nabla}\right). 
\]
Combined with  Corollary~\ref{cor:Amodel.is.good} this identification
reduces Conjecture~\ref{conn:one.parameter.ms}{\bfseries (ii)} to the
following purely algebro-geometric conjecture:

\begin{conn} \label{conn:Bmodel.is.good} Suppose $(Y,\bw)$ is a
  complex Landau-Ginzburg model which admits a tame compactification
  $(\Ybar,\wbar),D_{\Ybar})$ of log Calabi-Yau type. Then the
  associated \nc \ Hodge filtration
  $\left(\Bnc{H}^{\bullet},\Bnc{\nabla}\right)$ is extendable, and 
\[
\left(\stwo{\left(\Bnc{H}\right)}^{\bullet},
\stwo{\left(\Bnc{\nabla}\right)}\right) 
= 
\left(\Bmod{H}^{\bullet},\Bmod{\nabla}\right). 
\]
\end{conn}

\

\begin{rem} \label{rem:extendable} In a very interesting recent work
  Sabbah and Yu \cite{sabbah.yu} consider a different  but related notion
  of extendability arising from a nilpotent orbit for the pure complex
  Hodge structure attached to a compactification of a
  Landau-Ginzburg model. Moreover they prove that the scaling
  variation of this pure Hodge structure is polarizable and satisfies
  their extendability condition. This result seems 
  closely related to Conjecture~\ref{conn:Bmodel.is.good} but we have
  not investigated the precise relation between the two statemnts.
\end{rem}

\subsection{Canonical decorations} \label{ssec:decorations}

In this section we take a closer look at the data needed to write
special coordinates on the versal deformation space
$\mycal{M}$ of tame compactified
Landau-Ginzburg models $(\Ybar,\wbar),D_{\Ybar})$ of log Calabi-Yau
type. Recall from \cite{tian} and \cite{andrey}, that when $Y$ is a
smooth compact Calabi-Yau manifold of dimension $\dim_{\mathbb{C}} Y =
d$, then any choice of a splitting of the Hodge filtration on
$H^{d}_{DR}(Y,\mathbb{C})$ defines an analytic affine structure (= an
integrable torsion free connection on the tangent bundle on) on the
versal deformation space of $Y$. In
\cite[Section~4.1.3]{kkp} we analyzed the \nc \ counterpart of this
statement. In the \nc \ setting the splitting of the \nc \ Hodge
filtration is encodded in the notion of a {\it\bfseries decoration}
(see \cite[Definition~4.5]{kkp}) and in \cite[Claim~4.6]{kkp} we argue
that for decorated variations of pure \nc \ Hodge structures of
Calabi-Yau type there is a natural affine structure on the base of the
variation.

This analysis applies directly to
$\mycal{M}$ and the {\sf B}-model
variation of \nc \ Hodge structures over it. Concretely this variation
is given by a $\mathbb{Z}/2$-graded holomorphic bundle with connection 
$\left(\BBmod{H}^{\bullet},\BBmod{\nabla}\right)$ over $\mathbb{A}^{1}_{u}\times
\mycal{M}$, where the fiber of
$\BBmod{H}^{\bullet}$ over a point $\{ u \}\times
\{(\Ybar,\wbar),D_{\Ybar})\}$ is the hypercohomology 
$\mathbb{H}^{\bullet}\left(\Ybar; \left(\Omega^{\bullet}_{\Ybar}(\log
    D_{\Ybar},\wbar), ud + d\wbar\wedge \right)\right)$ and $\BBmod{\nabla}$
is the Gauss-Manin connection. Consider the projective line $\mathbb{P}^{1}_{u} =
\mathbb{A}^{1}_{u}\cup \{\infty \}$ compactifying
$\mathbb{A}^{1}_{u}$. As explained in \cite[Section~4.1.3]{kkp} the
special coordinates on $\mycal{M}$ arise from  decoration data for
$\left(\BBmod{H}^{\bullet},\BBmod{\nabla}\right)$. By definition a decoration is a pair
$\left(\BBmod{\widetilde{H}}^{\bullet},\psi\right)$, where 
\begin{itemize}
\item $\BBmod{\widetilde{H}}^{\bullet}$ is an extension of
$\BBmod{H}^{\bullet}$ to a $\mathbb{Z}/2$-graded holomorphic vector bundle on
$\mathbb{P}^{1}_{u}\times \mycal{M}$ for which $\BBmod{\nabla}$ has a
  logarithmic pole along  $\{\infty\} \times \mycal{M}$.
\item $\psi$ is a holomorphic section of
$\BBmod{\widetilde{H}}^{\bullet}_{|\{\infty\} \times \mycal{M}}$ which
is horizontal with respect to the holomorphic connection
${}^{\BBmod{\widetilde{H}}}\left(\BBmod{\nabla}\right)$ induced from
$\BBmod{\nabla}$ 
\end{itemize}

\

\begin{rem} \label{rem:multi.to.one}
The variation $\left(\BBmod{H}^{\bullet},\BBmod{\nabla}\right)$ is a
multi parameter variant of the one parameter variation
$\left(\Bmod{H}^{\bullet},\Bmod{\nabla}\right)$ we considered in
Section~\ref{sss-one.parameter.LG} and
Conjecture~\ref{conn:Bmodel.is.good}. In fact, the one parameter
variation $\left(\Bmod{H}^{\bullet},\Bmod{\nabla}\right)$ is the
restriction of the multi parameter variation
$\left(\BBmod{H}^{\bullet},\BBmod{\nabla}\right)$ to the straingt line
in $\mycal{M}$ given by the scaling of a fixed potential by a complex
number $q$.
\end{rem}

\

\noindent
In the remainder of this section we will describe a conjectural
construction which will produce a natural decoration in this
setting, i.e. will lead to canonical special coordinates that do not
depend on random choices. The construction is based on mirror symmetry
considerations and a description of decorations for the {\sf A}-model
variation of \nc \ Hodge structures. We begin by recalling the
relationship between filtrations and logarithmic extensions that we
used repeatedly in the previous section and in \cite[Section~4.1.3]{kkp}.

\

\medskip

\noindent
\punkt {\bfseries Extensions and filtrations} \label{sss:extensions} 
\ Let $\bD = \{ t \in
\mathbb{C} \ | \ |t| <  R \ll 1\}$ be a small one dimensional complex
disk centered at zero, and let $\bD^{\times} = \bD - \{0\}$ denote the
corresponding punctured disk. Let $(V,\nabla)$ be a holomorphic bundle
with holomorphic connection on $\bD^{\times}$, and suppose $\nabla$ is
meromorphic and has a regular singularity at $0$. By Deligne's
extension theorem \cite[Chapter~II.5]{deligne},
\cite[Corollary~II.2.21]{sabbah-frobenius} we can always find a
holomorphic bundle on $\bD$ which extends $V$, and on which $\nabla$
has a logarithmic pole. Fixing one such extension $\mathcal{V}$ as a
reference point we can use the Deligne-Malgrange classification
theorem \cite[Theorem~III.1.1]{sabbah-frobenius} to enumerate all
other logarithmic extensions of $(V,\nabla)$ by their relative
position to $\mathcal{V}$ (see
e.g. \cite[Chapter~III]{sabbah-frobenius} or
\cite[Section~4.1.3]{kkp}). In particular the choice of $\mathcal{V}$
gives a bijection
\[
\left( 
\text{
\begin{minipage}[c]{2.2in}
Holomorphic extensions of $V$ to $\bD$ on which $\nabla$ has
logarithmic pole at $0$
\end{minipage} 
}
\right)
\quad 
{\boldsymbol{\longleftrightarrow}}
\quad
\left(
\text{
\begin{minipage}[c]{2.5in}
Increasing biregular filtrations of $V$ by covariantly constant
holomorphic subbundles $V^{\leq i} \subset V$ on $\bD^{\times}$
\end{minipage} 
}
\right).
\]
If we
choose for concreteness $\mathcal{V}$ to be the unique Deligne
extension on which $\nabla$ has a logarithmic pole at $0$ and a
residue with eigenvalues whose real parts are in $(-1,0]$, then the
above bijection can be described explicitly as follows. Let
$\widetilde{V}$ be another
extension  of $V$ on which $\nabla$ has a logarithmic
pole. Fix an analytic trivialization of $\widetilde{V}$ near $t = 0$
and let $||\bullet ||$ denote the Hermitian norm of a section of $V$
computed in this trivialization.  
For any $t \in \bD$ and any $v \in V_{t}$ we have a well defined 
$\nabla$-horizontal section $s_{v}(r)$ of $V$ over the
segment $(0,1]\cdot t$ uniquely determined by the initial condition
  $s_{v}(1) = v$. With this notation we have
\[
V^{\leq i}_{t} = \left\{ v \in V_{t} \ \left|  \ \text{
\begin{minipage}[c]{2.2in} The $\nabla$-horizontal section $s_{v}(r)$
satisfies $||s_{v}(r)|| = \bO\left(r^{-i}\right)$
\end{minipage}
} \right.\right\}
\]
\

\begin{rem} $\bullet$ \  The growth condition defining $V^{\leq i}$ depends
  on the extension $\widetilde{V}$ but not on the choice of a local
  holomorphic frame of $\widetilde{V}$ near $0$.

$\bullet$ \  In \cite[Section~4.1.3]{kkp}) we discussed the classification
of logarithmic extensions of $(V,\nabla)$ in terms of biregular
decreasing filtartions of $V$. The above description of $V^{\leq i}$
is just a relabeling of the filtrations described in
\cite[Section~4.1.3]{kkp}). 
\end{rem}

\

\medskip

\noindent
\punkt {\bfseries The {\sf A}-model decoration.} \ Let
$(X,\omega_{X})$ be a compact symplectic manifold of real dimension
$2n$. Under the convergence assumption from 
\cite[Section~3.1]{kkp} for the quantum multiplication $*_{q}$, the \nc
\ Hodge filtration on the de Rham cohomology of $X$ is encoded in the
meromorphic connection on $\mathbb{A}^{1}_{u}$:
\[
\left(\Anc{H},\Anc{\nabla}\right) := \left( H^{\bullet}(X,\mathbb{C})\otimes
\mathcal{O}_{\mathbb{A}^{1}_{u}}, \ d +
(u^{-2}(\kappa_{X}\otimes_{1}(\bullet)) +  u^{-1}\sGr)du\right), 
\]
where $*_{1}$ denotes the quantum product for $\omega_{X}$. 
Note that by definition $\Anc{\nabla}$ has a regular singularity at $u
= \infty$. 

\begin{rem} \label{rem:convergence}
Conjecturally the convergence assumption on $*_{q}$ is closely related
to the properties of the \nc \ geometry attached to the
pair $(X,\omega_{X})$. In particular if convergence for $q = 1$ holds
we expect that 
\begin{itemize}
\item[(i)] the Fukaya category $\Fuk(X,\omega_{X})$ is smooth and
compact;
\item[(ii)] the geometrically  defined \nc \ Hodge filtration 
$\left(\Anc{H},\Anc{\nabla}\right)$ coincides with the
\nc \ Hodge filtration on
$HP_{\bullet}(\Fuk(X,\omega_{X}))$ defined in
\cite[Section~2.2.5]{kkp};
\item[(iii)] the monodromy of  $\Anc{\nabla}$
  around $u = \infty$ is unipotent and conjugate to the classical
  multiplication $\kappa_{X}\wedge (\bullet)$ by the canonical class.
\end{itemize}
Trough a combination of various results from
\cite{abouzaid-toric,reichelt-sevenheck,shklyarov-nchodge,efimov-cyclic}
properties (i)-(iii) are known to hold when $(X,\omega_{X})$ underlies
a smooth toric Fano variety.
\end{rem}

\

\noindent
As we explained  in Section~\ref{sss:extensions}, logaritmic extensions of
$\left(\Anc{H},\Anc{\nabla}\right)$ across $u = \infty$ will
correspond to a $\Anc{\nabla}$-horizontal filtrations of $\Anc{H}$.  In
particular from the definition of $\sGr$ we see that extensions of
$\Anc{H}$ as a trivial bundle over $\mathbb{P}^{1}_{u} =
\mathbb{A}^{1}_{u}\cup \{\infty\}$, will correspond to filtrations
$\Anc{H}^{\leq \bullet}$ whose associated graded is isomorphic to
$H^{\bullet + n}(X,\mathbb{C})\otimes
\mathcal{O}_{\mathbb{A}^{1}_{u}}$ as a $\mathbb{Z}$-graded bundle on
$\mathbb{A}^{1}_{u}$. Therefore,  in order to get a decoration for the {\sf
  A}-model data $\left(\Anc{H},\Anc{\nabla}\right)$ we need to specify
a canonical $\Anc{\nabla}$-covariant filtration on $\Anc{H}$ whose
associated graded pieces have dimensions $h^{\bullet +
  n}(X,\mathbb{C})$. If in addition this filtration depends
holomorphically on $\omega_{X}$, then it will automatically give a
decoration not only  for the fixed \nc \ Hodge filtration
$\left(\Anc{H},\Anc{\nabla}\right)$ but also for the universal
variation $\left(\AAmod{H},\AAmod{\nabla}\right)$ over the cone of
complexified  symplectic structures.

Such a canonical decoration arises naturally in the Fano case. Indeed,
suppose that $(X,\omega_{X})$ underlies a complex Fano manifold of
complex dimension $n$, and that property (iii) from
Remark~\ref{rem:convergence} holds for $(X,\omega_{X})$. In this case
the operator $\kappa_{X}\wedge(\bullet)$ satisfies the Lefschetz
property on $H^{\bullet}(X,\mathbb{C})$ and in particular has Jordan
blocks which are symmetrically situated around the middle dimension
$n$. In particular the Lefschetz filtartion (= the monodromy weight
filtration for the nilpotent operator $\kappa_{X}\wedge(\bullet)$)
will have associated graded pieces with dimensions $h^{\bullet +
  n}(X,\mathbb{C})$. Thus the extension of $\Anc{H}$ across $u =
\infty$ corresponding to this filtration will be holomorphically
trivial on $\mathbb{P}^{1}_{u}$. This shows that for a symplectic Fano
the universal Calabi-Yau variation of \nc \ Hodge structures
$\left(\AAmod{H},\AAmod{\nabla}\right) \to \mathbb{A}^{1}_{u}\times
\mycal{K}$ over the complexified K\"{a}hler cone $\mycal{K}$ will have
a canonical decoration data: 
\begin{description}
\item[$\AAmod{\widetilde{H}}$] is the holomorphic extension of
$\AAmod{H}$ to $\mathbb{P}^{1}_{u}\times \mycal{K}$ which ocrresponds 
to the monodromy weight filtration for the monodromy around  $u =
\infty$. 
\item[$\psi$] is the covariantly constant section of
  $\AAmod{\widetilde{H}}_{|\{\infty\}\times \mycal{K}}$ defined
  by 
$\psi(\infty,\beta) = s(\infty)$, where $s \in
  \Gamma\left(\mathbb{P}^{1}_{u}\times\{\beta\},
  \AAmod{\widetilde{H}}\right)$ is the unique holomorphic section in
  the trivial bundle
  $\AAmod{\widetilde{H}}_{|\mathbb{P}^{1}_{u}\times\{\beta\}} \cong
  H^{\bullet}(X,\mathbb{C})\otimes \mathcal{O}$ whose value at
  $(0,\beta)$ is $1 \in H^{0}(X,\mathbb{C})$.
\end{description}

\

\noindent
For ease of reference it will be useful to introduce terminology that
describes this extendability behavior. Again fix a small complex disk
$\bD$ and a meromorphic connection $(V,\nabla)$ on $\bD^{\times}$ with
a regular singularity and unipotent monoromy around zero. Fix the
unique Deligne extension $\mycal{V} \to \bD$ of $V$ on which $\nabla$
has a logarithmic pole with nilpotent residue. As we saw in 
Section~\ref{sss:extensions} this data establishes a $1$-to-$1$
correspondence between logarithmic extensions of $V$ to
$\mathbb{P}^{1}_{u}$ and covariantly constant biregular increasing
filtrations of $V \to \mathbb{A}^{1}_{u}$.

\begin{defi} \label{defi:skewed} $\bullet$ \ The {\bfseries skewed canonical
    extension} of $V$ is the holomorphic vector bundle $\widetilde{V}$
  which corresponds to the monodromy weight filtration for the
  monodromy operator around $u = \infty$. 

$\bullet$ \ An abstract \nc \ Hodge filtration $(H^{\bullet}\to
  \mathbb{A}^{1}_{u},\nabla)$ will be called {\bfseries special} if
  $\nabla$ has unipotent monodromy around $u = \infty$ and the
  corresponding skewed canonical extension $\widetilde{H}$ of $H$ is
  holomorphically trivial.
\end{defi}

\

\noindent
In these terms the  discussion about the extendability behaviour
of the {\sf A}-model \nc \ Hodge structure above can be rephrased as the
statement that the {\sf A}-model \nc \ Hodge filtration
$\left(\Anc{H},\Anc{\nabla}\right)$ associated with a symplectic Fano
variety is special. Speciality is the main property needed to
define the canonical decoration for the universal {\sf A}-model
variation $\left(\AAmod{H},\AAmod{\nabla}\right)$. 

Since the skewed extension and the speciality property are
intrinsically determined by the monodromy, it is straightforward to
transfer them through the mirror correspondence and to formulate the
{\sf B}-model extendability that will give rise to a canonical
decoration of $\left(\BBmod{H},\BBmod{\nabla}\right)$ and canonical
special coordinates on the moduli $\mycal{M}$ of compactified
Landau-Ginzburg models.

\

\medskip

\noindent
\punkt {\bfseries The {\sf B}-model decoration.} \ The {\sf A}-model
picture in the previous section and the mirror identification of the
{\sf A} and {\sf B} model universal variations of \nc \ Hodge
structures suggest that a canonical decoration for
$\left(\BBmod{H},\BBmod{\nabla}\right)$ and canonical special
coordinates on the moduli $\mycal{M}$ arise from the skewed extension
of the {\sf B}-model \nc \ Hodge structure. Specifically we get the
following purely algebro-geometric conjecture:

\begin{conn} \label{conn:Bmodel.can} 
\begin{itemize}
\item[(a)] Let $(Y,\bw)$ be a complex Landau-Ginzburg model, and let 
\[
\left(\Bnc{H}^{\bullet},\Bnc{\nabla}\right) =
\left(\mathbb{H}^{\bullet}\left(\Omega_{Y}^{\bullet}[u], ud -
d\bw\wedge\right), \ d + \left( u^{-2}\left(\bw\cdot (\bullet)\right)
+ u^{-1}\bG\right)\right)
\]
be the \nc \ Hodge filtration on the de Rham cohomology of
$(Y,\bw)$. If $(Y,\bw)$ admits a tame compactification
$((\Ybar,\wbar),D_{\Ybar})$ of Calabi-Yau type, then
$\left(\Bnc{H}^{\bullet},\Bnc{\nabla}\right)$ is special.
\item[(b)] Let $\mycal{M}$ be the versal deformation space of
  $((\Ybar,\wbar),D_{\Ybar})$. The universal
  {\sf B}-model variation 
  $\left(\BBmod{H}^{\bullet},\BBmod{\nabla}\right)$  over
  $\mathbb{A}^{1}_{u}\times \mycal{M}$ has a canonical decoration data:
\begin{description}
\item[$\BBmod{\widetilde{H}}$] is the skewed extension of
$\BBmod{H}$ to $\mathbb{P}^{1}_{u}\times \mycal{M}$. 
\item[$\psi$] is the covariantly constant section of
  $\BBmod{\widetilde{H}}_{|\{\infty\}\times \mycal{M}}$ defined by
  $\psi(\infty,((\Ybar,\wbar),D_{\Ybar})) = s(\infty)$, where $s \in
  \Gamma\left(\mathbb{P}^{1}_{u}\times\{((\Ybar,\wbar),D_{\Ybar})\},
    \BBmod{\widetilde{H}}\right)$ is the unique holomorphic section in
  the trivial bundle
  $\BBmod{\widetilde{H}}_{|\mathbb{P}^{1}_{u}\times\{((\Ybar,\wbar),D_{\Ybar})\}}
  \cong H^{\bullet}_{DR}((Y,\bw);\mathbb{C})\otimes \mathcal{O}$ whose
  value at $(0,((\Ybar,\wbar),D_{\Ybar}))$ is $1 \in
  H^{0}(\Ybar,\Omega^{n}_{\Ybar}(\log D_{\Ybar},\wbar))$.
\end{description}
\end{itemize}
\end{conn}


\begin{thebibliography}{FOOO09b}

\bibitem[AAE{\etalchar{+}}13]{abouzaidetal-punctured}
M.~Abouzaid, D.~Auroux, A.~Efimov, L.~Katzarkov, and D.~Orlov.
\newblock Homological mirror symmetry for punctured spheres.
\newblock {\em J. Amer. Math. Soc.}, 26(4):1051--1083, 2013.

\bibitem[ABC{\etalchar{+}}09]{mirrorbookD}
P.~Aspinwall, T.~Bridgeland, A.~Craw, M.~Douglas, M.~Gross, A.~Kapustin,
  G.~Moore, G.~Segal, B.~Szendr{\H{o}}i, and P.M.H. Wilson.
\newblock {\em Dirichlet branes and mirror symmetry}, volume~4 of {\em Clay
  Mathematics Monographs}.
\newblock American Mathematical Society, Providence, RI, 2009.

\bibitem[Abo09]{abouzaid-toric}
M.~Abouzaid.
\newblock Morse homology, tropical geometry, and homological mirror symmetry
  for toric varieties.
\newblock {\em Selecta Math. (N.S.)}, 15(2):189--270, 2009.

\bibitem[Abo12]{abouzaid-wrapped}
M.~Abouzaid.
\newblock On the wrapped {F}ukaya category and based loops.
\newblock {\em J. Symplectic Geom.}, 10(1):27--79, 2012.

\bibitem[{Abo}13]{abouzaid-viterbo}
M.~{Abouzaid}.
\newblock {Symplectic cohomology and Viterbo's theorem}.
\newblock {\em ArXiv e-prints}, page 182 pp, 2013, 1312.3354.

\bibitem[AKO06]{auroux-katzarkov-orlov2}
D.~Auroux, L.~Katzarkov, and D.~Orlov.
\newblock Mirror symmetry for del {P}ezzo surfaces: vanishing cycles and
  coherent sheaves.
\newblock {\em Invent. Math.}, 166(3):537--582, 2006.

\bibitem[AKO08]{auroux-katzarkov-orlov}
D.~Auroux, L.~Katzarkov, and D.~Orlov.
\newblock {Mirror symmetry for weighted projective planes and their
  noncommutative deformations}.
\newblock {\em Ann. of Math. (2)}, 167(3):867--943, 2008.

\bibitem[AS10]{abouzaid.seidel-viterbo}
M.~Abouzaid and P.~Seidel.
\newblock An open string analogue of {V}iterbo functoriality.
\newblock {\em Geom. Topol.}, 14(2):627--718, 2010.

\bibitem[Aur07]{auroux-anticanonical}
D.~Auroux.
\newblock Mirror symmetry and {$T$}-duality in the complement of an
  anticanonical divisor.
\newblock {\em J. G\"okova Geom. Topol. GGT}, 1:51--91, 2007.

\bibitem[BK98]{bk98}
S.~Barannikov and M.~Kontsevich.
\newblock Frobenius manifolds and formality of {L}ie algebras of polyvector
  fields.
\newblock {\em Internat. Math. Res. Notices}, 4:201--215, 1998.

\bibitem[Blo10]{block}
J.~Block.
\newblock Duality and equivalence of module categories in noncommutative
  geometry.
\newblock In {\em A celebration of the mathematical legacy of {R}aoul {B}ott},
  volume~50 of {\em CRM Proc. Lecture Notes}, pages 311--339. Amer. Math. Soc.,
  Providence, RI, 2010.

\bibitem[Bog79]{fedya-unobstructed}
F.~Bogomolov.
\newblock Hamiltonian {K}\"{a}hlerian manifolds.
\newblock {\em Soviet Math. Dokl}, 19:1462–1465, 1979.

\bibitem[Bog81]{fedya-ihes}
F.~Bogomolov.
\newblock K\"{a}hler varieties with trivial canonical class.
\newblock IHES preprint M/81/10, 1981.

\bibitem[Car57]{cartier}
P.~Cartier.
\newblock Une nouvelle op\'{e}ration sur les formes diff\'{e}rentielles.
\newblock {\em C.R.Ac.Sc.Paris}, 244:426--428, 1957.

\bibitem[Del70]{deligne}
P.~Deligne.
\newblock {\em \'{E}quations diff\'erentielles \`a points singuliers
  r\'eguliers}.
\newblock Springer-Verlag, Berlin, 1970.
\newblock Lecture Notes in Mathematics, Vol. 163.

\bibitem[Del71]{deligne-h2}
P.~Deligne.
\newblock Th\'eorie de {H}odge. {II}.
\newblock {\em Inst. Hautes \'Etudes Sci. Publ. Math.}, 40:5--57, 1971.

\bibitem[DI87]{deligne-illusie}
P.~Deligne and L.~Illusie.
\newblock Rel\`evements modulo {$p^2$} et d\'ecomposition du complexe de de
  {R}ham.
\newblock {\em Invent. Math.}, 89(2):247--270, 1987.

\bibitem[DMR07]{deligne-irregular}
P.~Deligne, B.~Malgrange, and J.-P. Ramis.
\newblock {\em Singularit\'es irr\'eguli\`eres}.
\newblock Documents Math\'ematiques (Paris) [Mathematical Documents (Paris)],
  5. Soci\'et\'e Math\'ematique de France, Paris, 2007.
\newblock Correspondance et documents. [Correspondence and documents].

\bibitem[Dub98]{dubrovin-berlin}
B.~Dubrovin.
\newblock Geometry and analytic theory of {F}robenius manifolds.
\newblock In {\em Proceedings of the International Congress of Mathematicians,
  Vol. II (Berlin, 1998)}, Extra Vol. II, pages 315--326 (electronic), 1998.

\bibitem[{Efi}12]{efimov-cyclic}
A.~{Efimov}.
\newblock {Cyclic homology of categories of matrix factorizations}.
\newblock {\em ArXiv e-prints}, December 2012, 1212.2859.

\bibitem[ESY13]{esy}
H.~{Esnault}, C.~{Sabbah}, and J.-D. {Yu}.
\newblock {'$E_{1}$-degeneration of the irregular {H}odge filtration (with an
  appendix by {M}orihiko {S}aito)}.
\newblock {\em ArXiv e-prints}, February 2013, 1302.4537.
\newblock 55 pp.

\bibitem[EV92]{esnault-viehweg}
H.~Esnault and E.~Viehweg.
\newblock {\em Lectures on vanishing theorems}, volume~20 of {\em DMV Seminar}.
\newblock Birkh\"auser Verlag, Basel, 1992.

\bibitem[FOOO09a]{fooo}
K.~Fukaya, Y-G. Oh, H.~Ohta, and K.~Ono.
\newblock {\em Lagrangian intersection {F}loer theory: anomaly and obstruction.
  {P}art {I}}, volume~46 of {\em AMS/IP Studies in Advanced Mathematics}.
\newblock American Mathematical Society, Providence, RI, 2009.

\bibitem[FOOO09b]{fooo2}
K.~Fukaya, Y.-G. Oh, H.~Ohta, and K.~Ono.
\newblock {\em Lagrangian intersection {F}loer theory: anomaly and obstruction.
  {P}art {II}}, volume~46 of {\em AMS/IP Studies in Advanced Mathematics}.
\newblock American Mathematical Society, Providence, RI, 2009.

\bibitem[Get93]{getzler-gm}
E.~Getzler.
\newblock Cartan homotopy formulas and the {G}auss-{M}anin connection in cyclic
  homology.
\newblock In {\em Quantum deformations of algebras and their representations
  (Ramat-Gan, 1991/1992; Rehovot, 1991/1992)}, volume~7 of {\em Israel Math.
  Conf. Proc.}, pages 65--78. Bar-Ilan Univ., Ramat Gan, 1993.

\bibitem[Giv98]{givental}
A.~Givental.
\newblock A mirror theorem for toric complete intersections.
\newblock In {\em Topological field theory, primitive forms and related topics
  (Kyoto, 1996)}, volume 160 of {\em Progr. Math.}, pages 141--175.
  Birkh\"auser Boston, Boston, MA, 1998.

\bibitem[GKR12]{gkr}
M.~{Gross}, L.~{Katzarkov}, and H.~{Ruddat}.
\newblock {Towards Mirror Symmetry for Varieties of General Type}.
\newblock {\em ArXiv e-prints}, February 2012, 1202.4042.

\bibitem[GM90]{goldman-millson}
W.~Goldman and J.~Millson.
\newblock The homotopy invariance of the {K}uranishi space.
\newblock {\em Illinois J. Math.}, 34(2):337--367, 1990.

\bibitem[Gri69]{griffiths-periods}
P.~Griffiths.
\newblock On the periods of certain rational integrals. {I}, {II}.
\newblock {\em Ann. of Math. (2) 90 (1969), 460-495; ibid. (2)}, 90:496--541,
  1969.

\bibitem[Gro66]{grothendieck-deRham}
A.~Grothendieck.
\newblock On the de {R}ham cohomology of algebraic varieties.
\newblock {\em Inst. Hautes \'Etudes Sci. Publ. Math.}, 29:95--103, 1966.

\bibitem[Hin01]{hinich}
V.~Hinich.
\newblock D{G} coalgebras as formal stacks.
\newblock {\em J. Pure Appl. Algebra}, 162(2-3):209--250, 2001.

\bibitem[HKK{\etalchar{+}}03]{mirrorbook}
K.~Hori, S.~Katz, A.~Klemm, R.~Pandharipande, R.~Thomas, C.~Vafa, R.~Vakil, and
  E.~Zaslow.
\newblock {\em Mirror symmetry}, volume~1 of {\em Clay Mathematics Monographs}.
\newblock American Mathematical Society, Providence, RI, 2003.
\newblock With a preface by Vafa.

\bibitem[Hor74]{horikawa}
E.~Horikawa.
\newblock On deformations of holomorphic maps. {II}.
\newblock {\em J. Math. Soc. Japan}, 26:647--667, 1974.

\bibitem[HV00]{hori-vafa}
K.~Hori and C.~Vafa.
\newblock Mirror symmetry, 2000, hep-th/0002222.

\bibitem[{Iac}13]{iacono}
D.~{Iacono}.
\newblock {Deformations and obstructions of pairs {({X},{D})}}.
\newblock {\em ArXiv e-prints}, February 2013, 1302.1149.

\bibitem[Ill71]{illusie1}
L.~Illusie.
\newblock {\em Complexe cotangent et d\'eformations. {I}}.
\newblock Lecture Notes in Mathematics, Vol. 239. Springer-Verlag, Berlin,
  1971.

\bibitem[Ill72]{illusie2}
L.~Illusie.
\newblock {\em Complexe cotangent et d\'eformations. {II}}.
\newblock Lecture Notes in Mathematics, Vol. 283. Springer-Verlag, Berlin,
  1972.

\bibitem[Ill02]{illusie}
L.~Illusie.
\newblock Frobenius and {H}odge degeneration.
\newblock In {\em Introduction to {H}odge theory}, volume~8 of {\em SMF/AMS
  Texts and Monographs}, pages 96--145. American Mathematical Society, 2002.
\newblock Translated from the 1996 French original by James Lewis and Chris
  Peters.

\bibitem[IM13]{iacono-manetti}
D.~Iacono and M.~Manetti.
\newblock Semiregularity and obstructions of complete intersections.
\newblock {\em Adv. Math.}, 235:92--125, 2013.

\bibitem[Kat70]{katz-nilpotent}
N.~Katz.
\newblock Nilpotent connections and the monodromy theorem.
\newblock {\em Publ. Math. IHES}, 39:175--233, 1970.

\bibitem[Kaw92]{kawamata}
Y.~Kawamata.
\newblock Unobstructed deformations. {A} remark on a paper of {Z}. {R}an:
  ``{D}eformations of manifolds with torsion or negative canonical bundle''
  [{J}. {A}lgebraic {G}eom.\ {\bf 1} (1992), no.\ 2, 279--291; {MR}1144440
  (93e:14015)].
\newblock {\em J. Algebraic Geom.}, 1(2):183--190, 1992.

\bibitem[KKP08]{kkp}
L.~Katzarkov, M.~Kontsevich, and T.~Pantev.
\newblock Hodge theoretic aspects of mirror symmetry.
\newblock In {\em From {H}odge theory to integrability and {TQFT}
  tt*-geometry}, volume~78 of {\em Proc. Sympos. Pure Math.}, pages 87--174.
  Amer. Math. Soc., Providence, RI, 2008.

\bibitem[KO68]{katz-oda}
N.~Katz and T.~Oda.
\newblock On the differentiation of de {R}ham cohomology classes with respect
  to parameters.
\newblock {\em J. Math. Kyoto Univ.}, 8:199--213, 1968.

\bibitem[Kon95]{kontsevich-icm}
M.~Kontsevich.
\newblock Homological algebra of mirror symmetry.
\newblock In {\em Proceedings of the International Congress of Mathematicians,
  Vol.\ 1, 2 (Z\"urich, 1994)}, pages 120--139, Basel, 1995. Birkh\"auser.

\bibitem[KS05]{ks-defos}
M.~Kontsevich and Y.~Soibelman.
\newblock Deformation theory. {V}olume {I}.
\newblock Unpublished book available at
  http://www.math.ksu.edu/~soibel/Book-vol1.ps, 2005.

\bibitem[KS09]{ks-ncgeometry}
M.~Kontsevich and Y.~Soibelman.
\newblock Notes on {$A_\infty$}-algebras, {$A_\infty$}-categories and
  non-commutative geometry.
\newblock In {\em Homological mirror symmetry}, volume 757 of {\em Lecture
  Notes in Phys.}, pages 153--219. Springer, Berlin, 2009.

\bibitem[Lan73]{landman}
A.~Landman.
\newblock On the {P}icard-{L}efschetz transformation for algebraic manifolds
  acquiring general singularities.
\newblock {\em Trans. Amer. Math. Soc.}, 181:89--126, 1973.

\bibitem[LO10]{lunts-orlov}
V.~Lunts and D.~Orlov.
\newblock Uniqueness of enhancement for triangulated categories.
\newblock {\em J. Amer. Math. Soc.}, 23(3):853--908, 2010.

\bibitem[LP11]{lin-pomerleano}
K.~H. {Lin} and D.~{Pomerleano}.
\newblock {Global matrix factorizations}.
\newblock {\em ArXiv e-prints}, January 2011, 1101.5847.

\bibitem[Lur11]{lurie-defos}
J.~Lurie.
\newblock Derived algebraic geometry {X}: {F}ormal moduli problems.
\newblock available at http://www.math.harvard.edu/~lurie/papers/DAG-X.pdf,
  2011.

\bibitem[Man99]{manin-frobenius}
Yu. Manin.
\newblock {\em Frobenius manifolds, quantum cohomology, and moduli spaces},
  volume~47 of {\em American Mathematical Society Colloquium Publications}.
\newblock American Mathematical Society, Providence, RI, 1999.

\bibitem[Man04]{manetti}
M.~Manetti.
\newblock Lectures on deformations of complex manifolds (deformations from
  differential graded viewpoint).
\newblock {\em Rend. Mat. Appl. (7)}, 24(1):1--183, 2004.

\bibitem[Orl04]{orlov}
D.~Orlov.
\newblock Triangulated categories of singularities and {D}-branes in
  {L}andau-{G}inzburg models.
\newblock {\em Tr. Mat. Inst. Steklova}, 246(Algebr. Geom. Metody, Svyazi i
  Prilozh.):240--262, 2004.

\bibitem[Orl05]{orlov-lagrange}
D.~Orlov.
\newblock Triangulated categories of singularities and equivalences between
  {L}andau-{G}inzburg, 2005.
\newblock math.AG/0503630.

\bibitem[Orl12]{orlov-global}
D.~Orlov.
\newblock Matrix factorizations for nonaffine {LG}-models.
\newblock {\em Math. Ann.}, 353(1):95--108, 2012.

\bibitem[OV05]{ogus.vologodsky}
A.~Ogus and V.~Vologodsky.
\newblock Nonabelian {H}odge theory in characteristic p, 2005, preprint,
  arXiv.org:math/0507476.

\bibitem[{Pos}11]{positselski-mf}
L.~{Positselski}.
\newblock {Coherent analogues of matrix factorizations and relative singularity
  categories}.
\newblock {\em ArXiv e-prints}, February 2011, 1102.0261.

\bibitem[{Pre}11]{preygel}
A.~{Preygel}.
\newblock {Thom-Sebastiani Duality for Matrix Factorizations}.
\newblock {\em ArXiv e-prints}, January 2011, 1101.5834.

\bibitem[PS08]{peters.steenbrink}
C.~Peters and J.~Steenbrink.
\newblock {\em Mixed {H}odge structures}, volume~52 of {\em Ergebnisse der
  Mathematik und ihrer Grenzgebiete. 3. Folge. A Series of Modern Surveys in
  Mathematics [Results in Mathematics and Related Areas. 3rd Series. A Series
  of Modern Surveys in Mathematics]}.
\newblock Springer-Verlag, Berlin, 2008.

\bibitem[PS13]{przyjalkowski.shramov}
V.~{Przyjalkowski} and C.~{Shramov}.
\newblock {On {H}odge numbers of complete intersections and
  {L}andau--{G}inzburg models}.
\newblock {\em ArXiv e-prints}, page 15 pp, May 2013, 1305.4377.

\bibitem[Ran92]{ran}
Z.~Ran.
\newblock Deformations of manifolds with torsion or negative canonical bundle.
\newblock {\em J. Algebraic Geom.}, 1(2):279--291, 1992.

\bibitem[RS10]{reichelt-sevenheck}
T.~{Reichelt} and C.~{Sevenheck}.
\newblock {Logarithmic Frobenius manifolds, hypergeometric systems and quantum
  {D}-modules}.
\newblock {\em ArXiv e-prints}, October 2010, 1010.2118.

\bibitem[Sab99]{sabbah-twisted}
C.~Sabbah.
\newblock On a twisted de {R}ham complex.
\newblock {\em Tohoku Math. J. (2)}, 51(1):125--140, 1999.

\bibitem[Sab02]{sabbah-frobenius}
C.~Sabbah.
\newblock {\em D\'eformations isomonodromiques et vari\'et\'es de {F}robenius}.
\newblock Savoirs Actuels (Les Ulis). [Current Scholarship (Les Ulis)]. EDP
  Sciences, Les Ulis, 2002.
\newblock Math\'ematiques (Les Ulis). [Mathematics (Les Ulis)].

\bibitem[Sab05]{sabbah-twistor}
C.~Sabbah.
\newblock Polarizable twistor {$\mycal{D}$}-modules.
\newblock {\em Ast\'erisque}, 300:vi+208, 2005.

\bibitem[Sab10a]{sabbah-irregular}
C.~Sabbah.
\newblock Fourier-{L}aplace transform of a variation of polarized complex
  {H}odge structure, {II}.
\newblock In {\em New developments in algebraic geometry, integrable systems
  and mirror symmetry ({RIMS}, {K}yoto, 2008)}, volume~59 of {\em Adv. Stud.
  Pure Math.}, pages 289--347. Math. Soc. Japan, Tokyo, 2010.

\bibitem[Sab10b]{sabbah}
C.~Sabbah.
\newblock Fourier-{L}aplace transform of a variation of polarized complex
  {H}odge structure, {II}.
\newblock In {\em New developments in algebraic geometry, integrable systems
  and mirror symmetry ({RIMS}, {K}yoto, 2008)}, volume~59 of {\em Adv. Stud.
  Pure Math.}, pages 289--347. Math. Soc. Japan, Tokyo, 2010.

\bibitem[Sai90]{saito-mhm}
M.~Saito.
\newblock Mixed {H}odge modules.
\newblock {\em Publ. Res. Inst. Math. Sci.}, 26(2):221--333, 1990.

\bibitem[Sai13]{saito.appendixE}
M.~Saito.
\newblock On the {K}ontsevich-de {R}ham complexes and {B}eilinson's maximal
  extensions.
\newblock In {\em ArXiv e-prints\/} \cite{esy}, 1302.4537.
\newblock Appendix E, p. 44--55.

\bibitem[{San}13]{sano}
T.~{Sano}.
\newblock {Unobstructedness of deformations of weak {F}ano manifolds}.
\newblock {\em ArXiv e-prints}, February 2013, 1302.0705.

\bibitem[Sch73]{schmid}
W.~Schmid.
\newblock Variation of {H}odge structure: the singularities of the period
  mapping.
\newblock {\em Invent. Math.}, 22:211--319, 1973.

\bibitem[Sei08]{seidel-book}
P.~Seidel.
\newblock {\em Fukaya categories and {P}icard-{L}efschetz theory}.
\newblock Zurich Lectures in Advanced Mathematics. European Mathematical
  Society (EMS), Z\"urich, 2008.

\bibitem[Sei09]{seidel-sh.as.hh}
P.~Seidel.
\newblock Symplectic homology as {H}ochschild homology.
\newblock In {\em Algebraic geometry---{S}eattle 2005. {P}art 1}, volume~80 of
  {\em Proc. Sympos. Pure Math.}, pages 415--434. Amer. Math. Soc., Providence,
  RI, 2009.

\bibitem[Sei11]{seidel-genus2}
P.~Seidel.
\newblock Homological mirror symmetry for the genus two curve.
\newblock {\em J. Algebraic Geom.}, 20(4):727--769, 2011.

\bibitem[Sei12]{seidel-AinftyI}
P.~Seidel.
\newblock Fukaya {$A_\infty$}-structures associated to {L}efschetz fibrations.
  {I}.
\newblock {\em J. Symplectic Geom.}, 10(3):325--388, 2012.

\bibitem[Sei14a]{seidel-AinftyII}
P.~Seidel.
\newblock {Fukaya $A_{\infty}$-structures associated to {L}efschetz fibrations.
  {I}{I}.}
\newblock {\em ArXiv e-prints}, page 52 pp, 2014, 1404.1352.

\bibitem[Sei14b]{seidel-ihes}
P.~Seidel.
\newblock Non-commutative geometry of {L}efschetz pencils.
\newblock talk at "Algebra, Geometry and Physics: a conference in honour of
  {M}axim {K}ontsevich", IHES, June 2014, 2014.

\bibitem[Ser06]{sernesi}
E.~Sernesi.
\newblock {\em Deformations of algebraic schemes}, volume 334 of {\em
  Grundlehren der Mathematischen Wissenschaften [Fundamental Principles of
  Mathematical Sciences]}.
\newblock Springer-Verlag, Berlin, 2006.

\bibitem[{Shk}11]{shklyarov-nchodge}
D.~{Shklyarov}.
\newblock {Non-commutative Hodge structures: Towards matching categorical and
  geometric examples}.
\newblock {\em ArXiv e-prints}, July 2011, 1107.3156.

\bibitem[Sim97]{simpson-hf}
C.~Simpson.
\newblock The {H}odge filtration on nonabelian cohomology.
\newblock In {\em Algebraic geometry---Santa Cruz 1995}, volume~62 of {\em
  Proc. Sympos. Pure Math.}, pages 217--281. Amer. Math. Soc., Providence, RI,
  1997.

\bibitem[SY14]{sabbah.yu}
C.~{Sabbah} and J.-D. {Yu}.
\newblock {On the irregular {H}odge filtration of exponentially twisted mixed
  {H}odge modules}.
\newblock {\em ArXiv e-prints}, page 53 pp, 2014, 1406.1339.

\bibitem[Tia87]{tian}
G.~Tian.
\newblock Smoothness of the universal deformation space of compact
  {C}alabi-{Y}au manifolds and its {P}etersson-{W}eil metric.
\newblock In {\em Mathematical aspects of string theory (San Diego, Calif.,
  1986)}, volume~1 of {\em Adv. Ser. Math. Phys.}, pages 629--646. World Sci.
  Publishing, Singapore, 1987.

\bibitem[Tod89]{andrey}
A.~Todorov.
\newblock The {W}eil-{P}etersson geometry of the moduli space of {${\rm
  SU}(n\geq 3)$} ({C}alabi-{Y}au) manifolds. {I}.
\newblock {\em Comm. Math. Phys.}, 126(2):325--346, 1989.

\bibitem[Wei96]{weibel-HC}
C.~Weibel.
\newblock Cyclic homology for schemes.
\newblock {\em Proc. Amer. Math. Soc.}, 124(6):1655--1662, 1996.

\bibitem[WW10]{wehrheim-woodward}
K.~Wehrheim and C.~Woodward.
\newblock Functoriality for {L}agrangian correspondences in {F}loer theory.
\newblock {\em Quantum Topol.}, 1(2):129--170, 2010.

\bibitem[{Yu}12]{yu-irregular}
J.-D. {Yu}.
\newblock {Irregular Hodge filtration on twisted de Rham cohomology}.
\newblock {\em ArXiv e-prints}, March 2012, 1203.2338.

\end{thebibliography}

\newcommand{\etalchar}[1]{$^{#1}$}

\

\bigskip

\noindent
Ludmil Katzarkov, {\sc University of Vienna}, 
ludmil.katzarkov@univie.ac.at

\smallskip

\noindent
Maxim Kontsevich, {\sc IHES, and University of Miami}, maxim@ihes.fr

\smallskip

\noindent
Tony Pantev, {\sc University of Pennsylvania}, tpantev@math.upenn.edu

\end{document}